\documentclass{amsart}
\usepackage{amssymb}
\usepackage[hyperref]{degt}

\makeatletter
\def\HyPsd@CatcodeWarning#1{}
\makeatother

\def\barN{\bar N}

\def\bS{\bold S}
\def\bZ{\bold Z}
\def\CK{\Cal K}
\def\CL{\Cal L}
\def\CS{\Cal S}
\def\CV{\Cal V}
\def\tker{\Cal C}

\def\F{\Bbb F}
\def\kk{\Bbbk}

\let\Artin\sigma

\def\T{V}

\def\bbT{\bold{\T}}
\def\bT{\bar\bbT}
\def\bba{\bold a}
\def\bbc{\bold c}
\def\ba{\bar\bba}
\def\bc{\bar\bbc}

\def\TL{\bold T}

\def\bX{\bold X}
\def\bH{\bold H}
\def\be{\bold e}

\def\genus{-\frak S}

\def\Golay{\Cal C}

\def\IS{\Cal I}
\def\fQ{\frak Q}
\def\fC{\frak C}
\def\fA{\frak Q}

\def\red{\QOPNAME{red}}
\def\dense{\QOPNAME{dense}}
\def\rt{\QOPNAME{rt}}
\def\supp{\QOPNAME{supp}}
\def\stab{\QOPNAME{stab}}
\def\Stab{\QOPNAME{Stab}}
\def\rle{\QOPNAME{rle}}
\def\mds{\QOPNAME{mds}}
\def\type{\QOPNAME{tp}}
\def\Fn{\QOPNAME{Fn}}
\def\Fano{\Cal F}

\def\dual{^\vee}
\def\norm|#1|{\mathopen\|#1\mathclose\|}

\def\Sym{\SG{}}
\def\Fermat{\Phi_4}

\def\symb#1{(#1)\mapsto}
\def\bb#1{\setbox0\hbox{$-$}\symb{\hbox to2\wd0{#1}}}
\def\bbr#1{\bb{\hss\hbox to\wd0{\hss$#1$\hss}}}
\def\bbl#1{\bb{\hbox to\wd0{\hss$#1$\hss}\hss}}
\def\bbcc#1{\bb{\hss$#1$\hss}}
\def\bbccc#1{\setbox0\hbox{$-$}\symb{\hbox to3\wd0{\hss$#1$\hss}}}

\def\2{2}
\def\6{6}
\def\8{8}
\def\quadric{V}

\pdef\bis{{}^{\dagger}}

\makeatletter
\newif\ifcusp
\newif\iftacnode
\let\refbar\relax
\def\punct#1{\rlap{#1}}
\def\bref#1{{\let\refbar\bar\ref{#1}}}
\newcommand\diagram[2][\undefined]{%
 \ifx\undefined#1\else
  \ifx*#1\stepcounter{A}\tag{$\refbar\theA$}\else\tag{$\refbar#1$}\fi\fi
 \left[\ \vcenter{\offinterlineskip
 \def\\{$\egroup\kern\bk@\hbox\bgroup$}%
 \setbox0\hbox{$=$\rlap{$\bullet$}}%
 \def\bxi##1{\setbox1\hbox to\wd0{##1}\ht1=\ht0\dp1=\dp0\box1}%
 \def\bxii##1{\setbox1\hbox to\bw\wd0{##1}\ht1=\ht0\dp1=\dp0\box1}%
 \def\bx##1{\bxi{\hss$##1$\hss}}%
 \def\rb##1{\bxii{\def\ {\bx{}}\hss$##1$\hss}}
 \def\rr##1{\rb{\def\.{\mathrel\circ}
                \def\*{\mathrel\bullet}
                \def\-{\joinrel\relbar\joinrel}##1}}
 \def\bb##1##2{\rb{##1\joinrel##2}}%
 \def\bk@{0pt}%
 \ifcusp
  \let\bw\tw@
  \def\bk@{2pt}%
  \def\-{\bb\relbar\relbar}%
  \def\={\bb\Relbar\Relbar}%
  \def\ {\rb{}}%
 \else\iftacnode
  \let\bw\thr@@
  \def\bk@{2pt}%
  \def\-{\rb{\relbar\joinrel\relbar\joinrel\relbar}}%
  \def\={\rb{\Relbar\joinrel\Relbar\joinrel\Relbar}}%
  \def\+{\rb{\relbar\joinrel\mathrel+\joinrel\relbar}}%
  \def\ {\rb{}}%
 \else
  \let\bw\@ne
  \def\+{\bx+}%
  \def\-{\bx-}%
  \def\={\bx=}%
  \def\ {\bx{}}%
  \def\d{\bx\cdot}
 \fi\fi
 \def\*{\bx*}%
 \def\.{\bx\circ}%
 \def\n{\bx\bullet}%
 \def\>{\bb\relbar{\joinrel\to}}%
 \def\R{\bb\Relbar{\joinrel\Rightarrow}}%
 \def\<{\bb{\leftarrow\joinrel}\relbar}%
 \def\L{\bb{\Leftarrow\joinrel}\Relbar}%
 \def\d{\rb\cdot}%
 \def\r{\bxii{\hss\n}}%
 \def\l{\bxii{\n\hss}}%
 \def\C{\bxii{\hss\n\hss}}%
 \def\c{\rb\.}%
 \hbox\bgroup$#2$\egroup}\vcenter to\big@size{}\ \right]%
}
\makeatletter

\makeatletter
\def\newquartic#1#2#3{\expandafter\gdef\csname q@#1\endcsname{{$#2$}{#3}}}
\def\qref#1{\ifmmode\text{#1}\else#1\fi}
\def\quartic{\relax\@ifstar\quarticii\quartici}
\def\quarticii#1{{\let\refbar\bar\quartici{#1}}}
\def\quartici#1{\qref{\@ifundefined{q@#1}{\ref{#1}}
 {\edef\q@{\csname q@#1\endcsname}%
 \hyperref[\expandafter\@secondoftwo\q@]{\expandafter\@firstoftwo\q@}}}}
\def\qlabel#1#2#3{\@bsphack\protected@write\@auxout{}
 {\string\newquartic{#1}{#2}{#3}}\@esphack}

\def\newMukai#1#2#3{\expandafter\gdef\csname Mukai@#1\endcsname{{#2}{#3}}}
\def\Mukailabel#1#2{\@bsphack\protected@write\@auxout{}
 {\string\newMukai{#1}{#1}{#2}}\@esphack}
\def\Mukaigroup#1#2{\raise 10pt\hbox{\hypertarget{#1}{}}#2\Mukailabel{#1}{#2}}
\def\Mukai#1{{\edef\q@{\csname Mukai@#1\endcsname}%
 \expandafter\ifx\csname Mukai@#1\endcsname\relax#1\else
 \hyperlink{\expandafter\@firstoftwo\q@}{\expandafter\@secondoftwo\q@}\fi}}

\newcounter{A}
\def\theA{\mathrm{D}\arabic{A}}
\makeatother

\def\Ps{\operatorname{\frak{ps}}}

\makeatletter
\def\PS#1{\expandafter\@PSi#1\end@PS}
\def\@PSi#1[#2[#3]#4,#5]{(#3)^{#5}\futurelet\next\@PSii}
\def\@PSii{\ifx\next,\expandafter\@PSi\fi}
\let\end@PS\relax
\makeatother

\def\ealign{\crcr
 \noalign{\vspace{2pt}\hrule}\egroup\egroup}

\pdef\newpaper{\url{http://www.fen.bilkent.edu.tr/~degt/papers/hh.pdf}}

\title{Smooth models of singular $K3$-surfaces}

\author{Alex Degtyarev}

\address{%
Department of Mathematics\\
Bilkent University\\
06800 Ankara, TURKEY}

\email{
degt@fen.bilkent.edu.tr}

\thanks{%
The author was partially supported by the T\"{U}B\DOTaccent{I}TAK grant 114F325}

\keywords{%
$K3$-surface, smooth quartic surface, Fermat quartic,
sextic curve, sextic model, octic model, Niemeier lattice, Mukai group%
}

\subjclass[2000]{%
Primary: 14J28;
Secondary: 14N25%
}

\begin{document}

\begin{abstract}
We show that the classical Fermat quartic has exactly three smooth spatial
models. As a generalization, we give a classification of smooth spatial (as
well as some other) models of singular $K3$-surfaces of small discriminant.
As a by-product, we observe a correlation (up to a certain limit) between the
discriminant of a singular $K3$-surface and the number of lines in its
models.
We also construct a $K3$-quartic surface with $52$ lines and singular
points, as well as a few other examples with many lines or models.
\end{abstract}

\maketitle

\section{Introduction}\label{S.intro}

All algebraic varieties considered in the paper
(except \autoref{s.fields})
are over~$\C$.

\subsection{Fermat quartics}\label{s.Fermat}
The original motivation for this paper was the classical \emph{Fermat
quartic} $\Fermat=\quartic{X48}\subset\Cp3$ given by the equation
\[*
z_0^4+z_1^4+z_2^4+z_3^4=0.
\]
It is immediate from the equation that \quartic{X48} contains $48$ straight
lines, \viz.
\[*
z_a-\epsilon_1z_b=z_c-\epsilon_2z_d=0,
\]
where $\epsilon_1^4=\epsilon_2^4=-1$ and $\{\{a,b\},\{c,d\}\}$ is an
unordered partition of the index set $\{0,\ldots,3\}$ into two unordered
pairs.
The maximal possible number of lines in a smooth quartic surface is~$64$
(see \cite{Segre,rams.schuett}) and
there are but
ten (eight up to complex conjugation) quartics with more than $52$ lines
(see~\cite{DIS}).
When \cite{DIS} appeared,
it was immediately observed by T.~Shioda that one
of these extremal quartic, \viz. \quartic{X56} in the notation of~\cite{DIS},
is isomorphic, as an abstract $K3$-surface, to the classical Fermat
quartic~\quartic{X48}. This observation resulted in a beautiful
paper~\cite{Shimada:X56}, which provides explicit defining equations
for the surface~\quartic{X56} and isomorphism $\quartic{X48}\cong\quartic{X56}$ and
studies further geometric properties of~\quartic{X56}.
This explicit construction (first to my knowledge) is particularly
interesting  due to the fact
(see \cite[Theorem 1.8]{Oguiso} with a further
reference to~\cite{Matsumura})
that $(d,n)=(4,3)$ is the only pair
with $n\ge3$ for which
two smooth hypersurfaces of
the same
degree~$d$ in~$\Cp{n}$ may be isomorphic as
abstract algebraic varieties but not projectively equivalent.

Smooth quartics in~$\Cp3$ are $K3$-surfaces, and we define a
\emph{smooth spatial model} of a $K3$-surface~$X$ as an embedding
$X\into\Cp3$ defined by a very ample line bundle of degree~$4$. Two models
are \emph{projectively equivalent} if so are their images.
According to the authors of~\cite{Shimada:X56}, an extensive search for other
smooth spatial models of the Fermat quartic~$\Fermat$ did not produce any
results, suggesting that such models do not exist. This assertion is one of
the principal results of the present paper.

\theorem[see \autoref{proof.main}]\label{th.main}
Up to projective equivalence, there are three smooth models
$\Fermat\into\Cp3$ of the
Fermat quartic\rom: they are \quartic{X48}, \quartic{X56}, and \quartic*{X56}.
\endtheorem

A detailed proof of this theorem is given in
\autoref{S.Niemeier} and \autoref{proof.main}.

This statement is rather surprising, because
there must be several thousands of \emph{singular} spatial models of the Fermat
quartic (\cf. \autoref{rem.singular}; we do not attempt their
classification) and because the number of distinct spatial models
of a singular $K3$-surface
growth rather fast with the discriminant (\cf.
\autoref{rem.growth}
in \autoref{proof.infty}
or the discussion
of the smooth models of $X([2,1,82])$ in \autoref{s.Q}).

\subsection{Quartics of small discriminant}\label{s.det<80}
The approach used in the classification of the spatial
models of~$\Fermat$ applies
to other
$K3$-surfaces.
We confine ourselves to the most interesting case of the so-called
\emph{singular} $K3$-surfaces~$X$, \ie, those of
the maximal Picard rank $\rank\NS(X)=20$. Due to the global Torelli theorem,
such a
$K3$-surface~$X$ is determined
by its \emph{transcendental lattice}
$T:=\NS(X)^\perp\subset H_2(X)$, considered up to orientation preserving
isometry (see, \eg,~\cite{Shioda.Inose} and \autoref{s.K3}).
This is a positive definite even lattice of rank~$2$; we use the
notation $X:=X(T)$ and the single line (instead of a matrix)
notation $T:=[a,b,c]$ for
the lattice $T=\Z u+\Z v$, $u^2=a$, $v^2=c$, $u\cdot v=b$
(see \autoref{s.2x2} and \autoref{s.K3}).
To illustrate the approach, we outline the proof of the following theorem,
which
formally
incorporates \autoref{th.main}. For a technical reason
(see \autoref{rem.genus}),
we omit the
case $T=[4,0,16]$; conjecturally,
the corresponding $K3$-surface $X(T)$ has no
smooth spatial models.

\theorem[see \autoref{proof.others}]\label{th.others}
Let $T$ be a positive definite even lattice of rank~$2$, and assume that
$\det T\le80$ and $T\ne[4,0,16]$. Then, up to projective equivalence, any
smooth spatial model $X(T)\into\Cp3$ is one of those listed in
\autoref{tab.quartics}.
\table
\def\*{\rlap{$^*$}}
\def\+{\noalign{\kern2\lineskip}}
\def\-{\noalign{\kern-\lineskip}}
\rm
\caption{Nonsingular spatial models (see \autoref{th.others})}\label{tab.quartics}
\hbox to\hsize\bgroup\hss\vbox\bgroup
\halign\bgroup\strut\quad\hss$#$\hss\quad&\hss$#$\hss\quad&\hss$#$\hss\quad&
 $#$\hss\quad\cr
\noalign{\hrule\vspace{2pt}}%
\omit\strut\hidewidth$\det$\hidewidth&T&X&\text{Pencils, remarks}\cr
\noalign{\vspace{1pt}\hrule\vspace{2pt}}
48&[8,4,8]&\quartic{X64}&
 \PS{[ [ 6, 0 ], 16 ], [ [ 4, 6 ], 48 ]};\
 \text{Mukai group $\Mukai{T192}$}\cr\+
55&[4,1,14]\*&\quartic{X60''}&
 \PS{[ [ 4, 5 ], 60 ]}\cr\+
60&[4,2,16]&\quartic{X60'}&
 \PS{[ [ 6, 2 ], 10 ], [ [ 4, 4 ], 30 ], [ [ 3, 7 ], 20 ]}\cr\-
  &&\quartic{Q56}&
 \PS{[ [ 4, 4 ], 24 ], [ [ 3, 7 ], 32 ]}\cr\+
64&[2,0,32]&\quartic{Y56}&
 \PS{[ [ 4, 4 ], 32 ], [ [ 3, 7 ], 24 ]}\cr\+
64&[4,0,16]&??&
 \text{Conjecturally, none}\cr\+
64&[8,0,8]&\quartic{X56}\*&
 \PS{[ [ 4, 6 ], 8 ], [ [ 4, 4 ], 32 ], [ [ 2, 8 ], 16 ]}\
 \text{(see \autopageref{X56})}\cr\-
  &&\quartic{X48}&
 \PS{[ [ 2, 8 ], 48 ]}\
 \text{(see \autopageref{X48}); Mukai group $\Mukai{F384}$}\cr\+
75&[10,5,10]&\quartic{Q52'''}&
 \PS{[ [ 5, 0 ], 4 ], [ [ 3, 6 ], 48 ]}\cr\+
76&[2,0,38]&\quartic{Y52'}&
 \PS{[ [ 4, 6 ], 2 ], [ [ 4, 4 ], 16 ], [ [ 3, 5 ], 20 ], [ [ 2, 8 ], 14 ]}\cr\-
  &[8,2,10]\*&\quartic{Y48'}&
 \PS{[ [ 3, 5 ], 24 ], [ [ 2, 8 ], 24 ]}\cr\+
76&[4,2,20]&\quartic{Q54}&
 \PS{[ [ 4, 4 ], 24 ], [ [ 4, 3 ], 24 ], [ [ 0, 12 ], 6 ]};
    \ \text{see \eqref{X52''}}\cr\-
  &&\quartic{X52''}&
 \PS{[ [ 6, 0 ], 1 ], [ [ 4, 4 ], 9 ], [ [ 4, 3 ], 18 ], [ [ 3, 5 ], 18 ],
 [ [ 0, 12 ], 6 ]};\ \text{see \eqref{X52''}}\cr\+
79&[2,1,40]&\quartic{Y52''}&
 \PS{[ [ 4, 5 ], 8 ], [ [ 4, 3 ], 12 ], [ [ 3, 6 ], 16 ], [ [ 2, 7 ], 16 ]}\cr\-
  &[4,1,20]\*&\quartic{Y48''}&
 \PS{[ [ 4, 3 ], 6 ], [ [ 3, 6 ], 12 ], [ [ 2, 7 ], 30 ]}\cr\-
  &[8,1,10]\*\cr\+
80&[4,0,20]&\quartic{Z52}&
 \PS{[ [ 6, 0 ], 4 ], [ [ 4, 4 ], 12 ], [ [ 4, 2 ], 24 ], [ [ 2, 8 ], 12 ]};
    \ \text{see \eqref{Z48'}, \eqref{Z52}}\cr\-
  &&\quartic{Z50}\*&
 \PS{[ [ 4, 4 ], 10 ], [ [ 3, 5 ], 40 ]}\cr\-
  &&\quartic{Z48'}&
 \PS{[ [ 4, 2 ], 16 ], [ [ 2, 8 ], 32 ]};\ \text{see \eqref{Z48'}}\cr\-
  &&\quartic{Z48''}\*&
 \PS{[ [ 3, 5 ], 48 ]};\ \text{see \eqref{Z48'}, \eqref{Z52}}\cr\+
80&[8,4,12]&\quartic{X52'}&
 \PS{[ [ 6, 0 ], 1 ], [ [ 4, 4 ], 12 ], [ [ 4, 3 ], 12 ], [ [ 4, 2 ], 3 ],
   [ [ 3, 5 ], 18 ], [ [ 0, 12 ], 6 ]}\cr\-
  &&\quartic{Q52''}\*&
 \PS{[ [ 4, 4 ], 8 ], [ [ 4, 3 ], 32 ], [ [ 4, 2 ], 8 ], [ [ 0, 12 ], 4 ]}\cr\-
  &&\quartic{Q48}&
 \PS{[ [ 4, 2 ], 8 ], [ [ 3, 5 ], 32 ], [ [ 0, 12 ], 8 ]};
    \ \text{see \eqref{Q48}}\cr
\ealign\hss\egroup
\endtable
\endtheorem

In \autoref{tab.quartics},
the lattices~$T$ are grouped according to their genus
(or, equivalently, isomorphism class of the discriminant form, see
\autoref{s.lattice} for details). For each genus, we list the discriminant
$\det:=\det T$, the isomorphism classes of lattices~$T$
(line by line), the smooth
spatial models~$X$ (using, whenever possible, the notation
introduced
in~\cite{DIS} for the configurations of lines,
the subscript always indicating the number of lines;
all models apply to all lattices in the genus), and
the pencil structure $\Ps(X)$,
which can be used to identify quartics
(see \autoref{s.pencils}).
The number of models of $X(T)$ is constant within
each genus (see \autoref{rem.genus.constant}).
Lattices admitting no orientation reversing automorphism are marked with
a $^*$; the corresponding $K3$-surfaces are not real, and neither are their
models. The $^*$ next to a model designates the fact that, although the
$K3$-surface itself is defined over~$\R$, the model is not.

The
notation $X_*$, $Y_*$, \etc. refers to the diagrams found on
pages~\pageref{X64}--\pageref{Q52''-2} (unless indicated otherwise in the
table); diagrams corresponding to several models are denoted by $(D*)$.
In principle, these diagrams suffice to recover the lattices; details are
explained in the relevant parts of the proof.

Two of the quartics, \viz. \quartic{X64} and \quartic{X48}, carry faithful actions
of \emph{Mukai groups} (\ie, maximal finite groups of symplectic
automorphisms, see~\cite{Mukai} and \autoref{s.Mukai});
these are all quartics with this property
that contain lines (see \autoref{th.Mukai}).

Note that there are $97$ lattices~$T$, constituting $78$ genera,
satisfying the condition $\det T\le80$.
According to \autoref{th.others}, very few of the
corresponding $K3$-surfaces admit smooth spatial models.

The configurations of lines not found in~\cite{DIS} are \quartic{X48} (the
classical Fermat quartic), \quartic{Q52'''}, the alternative models
\quartic{Y48'}, \quartic{Y48''} of \quartic{Y52'}, \quartic{Y52''},
respectively, and the alternative models \quartic{Z48'}, \quartic{Z48''}
of \quartic{Z52}, \quartic{Z50} and \quartic{Q48} of \quartic{X52'},
\quartic{Q52''}.
(Following~\cite{DIS}, we denote by~$Z_*$ the configurations of lines
generating a sublattice of rank~$19$ in $\NS(X)$. In all other cases
in \autoref{tab.quartics},
the space
$\NS(X)\otimes\Q$ is generated by
the classes of
the lines contained in~$X$.)
All
quartics given by \autoref{th.others} contain many (at least $48$) lines
and,
conversely, all eight quartics with at least $56$ lines
(see~\cite{DIS}) do appear in the table. This observation may shed new light
on the line counting problem. For example, the following statement is an
alternative characterization of Schur's quartic \quartic{X64}, given by the
equation
\[*
z_0(z_0^3-z_1^3)=z_2(z_2^3-z_3^3).
\]
(Recall that, according to~\cite{DIS},
Schur's quartic
\quartic{X64} is the only smooth quartic that
contains the maximal
number~$64$ of lines).

\corollary\label{cor.Schur}
If a singular $K3$-surface $X(T)$ admits a smooth spatial model, then either
$T=[8,4,8]$ and the model is Schur's quartic \quartic{X64}, or $\det T\ge55$.
\done
\endcorollary

In other words,
$48$ is the minimal discriminant of a singular $K3$-surface
admitting a smooth spatial model, and Schur's quartic~\quartic{X64}
is the only one
minimizing this discriminant,
just as it is the only quartic maximizing the number of lines.
This statement, together with the next corollary
(which follows from the proof of \autoref{th.others}),
substantiates \autoref{conj.max} below. (Though, see also \autoref{rem.max}.)


\corollary[see \autoref{proof.P1xP1}]\label{cor.P1xP1}
Let
$T$ be as above, $\det T\le80$, and $T\ne[4,0,16]$. Then, $X(T)$ admits a
degree~$2$ map $X(T)\to Q:=\Cp1\times\Cp1$ with smooth ramification locus
$C\subset Q$
if and only if either
\roster*
\item
$T=[4,2,20]$, and the model is \quartic{PP48},
see \autopageref{PP48}, or
\item
$T=[8,4,12]$, and the model is \quartic{PP48''}, see~\eqref{PP48''}
on \autopageref{PP48''}.
\endroster
In each case, the model is unique up to
projective equivalence
and $C$ has
the maximal number $12$ of bitangents in each
of the two rulings of~$Q$
\rom(see also \autoref{rem.P1xP1}\rom).
\endcorollary

Certainly,
when speaking about bitangents of the ramification locus~$C$,
we admit the degeneration of a bitangent to a
generatrix intersecting~$C$ at a single point with multiplicity~$4$ (\cf. the
discussion of tritangents
right before \autoref{th.polarizations}
below). Note, though,
that in the extremal case as in \autoref{cor.P1xP1},
all twelve generatrices are true bitangents, as
follows immediately (together with the bound~$12$ itself)
from the Riemann--Hurwitz formula.

With few exceptions,
any pair $X',X''\subset\Cp3$ of smooth models (not necessarily distinct)
of the same
$K3$-surface $X(T)$ appearing in \autoref{tab.quartics} constitutes a
so-called \emph{Oguiso pair} (see \cite{Oguiso} and \autoref{s.Oguiso});
in particular, the models are Cremona equivalent.
The exceptions are:
\roster*
\item
the quadruple $Z_*$ of models of $X([4,0,20])$, which splits
into two pairs, \viz. $(\quartic{Z52},\quartic{Z48'})$ and
$(\quartic{Z50},\quartic{Z48''})$, each having the above property, and
\item
the Fermat quartic~\quartic{X48}:
there is no Oguiso pair $(\quartic{X48},\quartic{X48})$ (\cf.
\cite{Shimada:X56}).
\endroster
If $X'=X''$,
then
Oguiso's
construction~\cite{Oguiso}
gives us a Cremona self-equivalence $X'\to X''$ that is not
regular on the ambient space~$\Cp3$.
It follows also that each model is a \emph{Cayley $K3$-surface},
\ie, a smooth determinantal
quartic (see \cite{Cayley:quartics,Oguiso}).
These phenomena are specific to small discriminants, see
\autoref{th.infty} below.

\subsection{Other polarizations}\label{s.polarizations}
The approach applies as well to other polarizations of $K3$-surfaces, \ie,
projective models $\Gf\:X\into\Cp{n}$ defined by linear
systems~$\ls|h|$, $h\in\NS(X)$, $h^2=2n-2$. We will consider the following
commonly used models.
\roster
\item\label{i.planar}
$h^2=2$: a \emph{planar} model $\Gf\:X\to\Cp2$. The model is a degree~$2$ map
ramified at a sextic curve $C\subset\Cp2$; it is called \emph{smooth} if so
is~$C$ (\cf. \autoref{cor.P1xP1}).
\item\label{i.spatial}
$h^2=4$: a \emph{spatial} or \emph{quartic} model $\Gf\:X\to\Cp3$ considered
in \autoref{th.others}.
\item\label{i.sextic}
$h^2=6$: a \emph{sextic} model $\Gf\:X\to\Cp4$. The image of $\Gf$ is a
complete intersection (regular if $\Gf$ is smooth) of a quadric and a cubic.
\item\label{i.octic}
$h^2=8$: an \emph{octic} model $\Gf\:X\to\Cp5$.
Typically, the image of $\Gf$ is a
complete intersection (regular if $\Gf$ is smooth) of three quadrics
(\cf. \autoref{lem.octic}).
\endroster
In the first case, a \emph{line} in~$X$ is a smooth rational curve that
projects isomorphically to a line in~$\Cp2$.
The projection establishes a two-to-one correspondence between lines and
tritangents of~$C$. (Here,
we admit the possibility that a tritangent degenerates to a line
intersecting~$C$ at two points with multiplicities~$2$ and~$4$ or at a single
point with multiplicity~$6$. More generally, if $C$ is allowed to be
singular, the ``tritangents'' are the so-called \emph{splitting lines}, \ie,
lines whose local intersection index with~$C$ is even at each intersection
point.)

\theorem[see \autoref{proof.polarizations}]\label{th.polarizations}
Let $T$ be a positive definite even lattice of rank~$2$.
\roster
\item\label{i.p.planar}
If $\det T\le116$ and $X(T)$ admits a smooth planar model
$X(T)\to\Cp2$, then
$T=[12,6,12]$ and the only model is \quartic{2-144},
see \autopageref{2-144}.
\setcounter{enumi}2
\item\label{i.p.sextic}
If $\det T\le48$, then, up to projective equivalence,
any smooth sextic model
$X(T)\into\Cp4$ is one of those listed in \autoref{tab.sextics}.
\item\label{i.p.octic}
If $\det T\le40$, then, up to projective equivalence, any smooth octic model
$X(T)\into\Cp5$ is one of those listed in \autoref{tab.octics}.
\bgroup
\def\*{\rlap{$^*$}}
\def\d{\rlap{$^\diamond$}}
\def\+{\noalign{\kern2\lineskip}}
\def\-{\noalign{\kern-\lineskip}}
\def\balign{\vbox\bgroup\rm\halign\bgroup
 \strut\quad\hss$##$\hss\quad&\hss$##$\hss\quad&\hss$##$\hss\quad&
  \hss$##$\hss\quad&$##$\hss\quad\cr
 \noalign{\hrule\vspace{2pt}}
 \det&T&X&\text{Ranks}&\text{Pencils, remarks}\cr
 \noalign{\vspace{1pt}\hrule\vspace{2pt}}}
\table
\caption{Smooth sextic models (see \autoref{th.polarizations}\iref{i.p.sextic})}\label{tab.sextics}
\hbox to\hsize{\hss\balign
39&[2,1,20]&\quartic{6-42}&(19,19)&
 \PS{[ [ 9 ], 42 ]}\cr\-                                
  &[6,3,8]&\cr\+
48&[6,0,8]&\quartic{6-42'}&(19,19)&
 \PS{[ [ 9 ], 42 ]};\ \text{see \eqref{6-36'}}\cr\-     
  &&\quartic{6-38}&(19,19)&
 \PS{[ [ 11 ], 2 ], [ [ 9 ], 16 ], [ [ 7 ], 20 ]}\cr\-  
  &&\quartic{6-36'}&(19,19)&
 \PS{[ [ 8 ], 36 ]};\ \text{see \eqref{6-36'}}\cr\+     
48&[8,4,8]&\quartic{6-36''}&(18,18)&
 \PS{[ [ 8 ], 36 ]}\cr                                  
\ealign\hss}
\endtable
\table
\caption{Smooth octic models (see \autoref{th.polarizations}\iref{i.p.octic})}\label{tab.octics}
\hbox to\hsize{\hss\balign
32&[4,0,8]&\quartic{8-36}\*&(20,20)&
 \PS{[ [ 7 ], 16 ], [ [ 6 ], 16 ], [ [ 4 ], 4 ]}\cr\-   
  &&\quartic{8-32}&(17,17)&
 \PS{[ [ 6 ], 32 ]};\
 \text{Mukai group $\Mukai{F384}$}\cr\+                      
36&[6,0,6]&\quartic{8-36'}&(19,20)&
 \PS{[ [ 6 ], 36 ]};\ \text{see \eqref{8-33}}\cr\-      
  &&\quartic{8-33}\d&(18,18)&
 \PS{[ [ 8 ], 9 ], [ [ 5 ], 24 ]};
    \ \text{see \eqref{8-33}}\cr\-                      
  &&\quartic{8-32'}&(18,18)&
 \PS{[ [ 6 ], 32 ]}\cr\+                                
39&[2,1,20]&\quartic{8-30}&(18,19)&
 \PS{[ [ 6 ], 12 ], [ [ 5 ], 12 ], [ [ 4 ], 6 ]}\cr\-   
  &[6,3,8]&\cr
\ealign\hss}
\endtable
\egroup
\endroster
\endtheorem

In the tables, we use conventions similar to \autoref{tab.quartics},
referring to the diagrams found on
pages~\pageref{2-144}--\pageref{8-33-2}.
As an additional invariant, we list the ranks of the sublattice $F\subset\NS(X)$
generated by the classes of lines and its extension $F+\Z h$.
Instead of the pencil structure, we merely list the valencies of the vertices
of the dual adjacency graph of the lines.
The model $\quartic{8-33}$ marked with a diamond~${}^\diamond$ is the only one
(in \autoref{tab.octics}) whose
defining ideal is not generated by polynomials of degree~$2$
(see \autoref{lem.octic}).
With the exception of \quartic{6-42} and \quartic{6-42'},
all configurations found in the
tables are pairwise distinct,
as the invariants show.

The ramification locus of \quartic{2-144} admits a faithful action of the
Mukai group~$\Mukai{M9}$; hence, according to Sh\.~Mukai~\cite{Mukai}, its equation is
\[*
z_0^6+z_1^6+z_2^6=10(z_0^3z_1^3+z_1^3z_2^3+z_2^3z_0^3).
\]
In \autoref{s.Mukai}, we show that
no sextic model containing lines carries a faithful action of a Mukai group; the
two octics with this property are \quartic{8-32} and a model of
$X([4,0,12])$, with the same configuration of lines and an action
of~$\Mukai{H192}$ (see \autoref{th.Mukai}).

As a by-product, we obtain new lower bounds on the maximal number of lines in
a model. (According to S.~Rams, private communication, no interesting
examples of sextics are known, whereas the best known example of an octic has
$32$ lines; octics with $32$ lines have been studied in~\cite{Dolgachev:book}.)
Comparing \autoref{th.polarizations} and
Corollaries~\ref{cor.Schur} and~\ref{cor.P1xP1}, we conjecture that these
new bounds are sharp.

\conjecture\label{conj.max}
A smooth sextic curve $C\subset\Cp2$ has at most $72$ tritangents. A smooth
sextic \rom(octic\rom) model of a $K3$-surface has at most $42$
\rom(respectively, $36$\rom) lines.
\endconjecture

\remark\label{rem.max}
The
cases of sextic and octic models in \autoref{conj.max} are settled, in
the affirmative, in~\cite{degt:lines}, where we also give a sharp bound on
the maximal number of lines in a smooth $K3$-surface $X\into\Cp{D+1}$ for
all $2\le D\le15$. (For $D\ge16$, the bound is~$24$
and its sharpness depends on
the residue $D\bmod12$; for all $D\gg0$, large configurations of lines are
fiber components of elliptic pencils.)
Thus, the only case that remains open is that of
plane sextic curves.

Curiously, the motivating observation, \viz. the fact
that the number of lines is maximized
by the discriminant minimizing singular $K3$-surfaces, does not persist for
higher polarizations: for degree~$10$ surfaces in $\Cp6$,
the discriminant minimizing surface $X([2,0,16])$ has fewer (28) lines than
the maximum~$30$.
\endremark

\subsection{Further examples}\label{s.results}
Each $K3$-surface $X(T)$ found in \autoref{tab.quartics} has at most three (two up to
projective equivalence and complex conjugation)
distinct smooth spatial models. The number of models of any particular
$K3$-surface is always finite, but we show that this number is not bounded.
(The former statement, which is an immediate consequence from the finiteness of
each genus of
lattices, was first obtained by Sterk~\cite{Sterk}.)

\theorem[see \autoref{proof.infty}]\label{th.infty}
For each integer~$d$, the number of projective equivalence classes of
spatial models
\rom(smooth or not\rom) $X(T)\to\Cp3$ with $\det T\le d$ is finite.

However,
for each integer $M>0$, there exist a lattice $T$
and $M$ smooth spatial models $X(T)\into\Cp3$
that are pairwise not Cremona equivalent.
\endtheorem

An almost literate
analogue of this statement for the other three polarizations considered in
\autoref{s.polarizations} is discussed at the end of \autoref{proof.infty}.

Another misleading observation suggested by \autoref{tab.quartics} is the
fact that the number of lines in smooth spatial models of $X(T)$ tends to
decrease when $\det T$ increases.
This tendency persists only up to a certain limit, \viz. $52$ lines (which is
the maximal number of lines realized by an equilinear $1$-parameter family of
quartics, see~\cite{DIS}).
Furthermore, the number of such models of a given $K3$-surface is not bounded,
providing an alternative series of examples,
with a large number of lines but much lower growth rate (\cf. \autoref{rem.growth}),
for the statement of \autoref{th.infty}.

\proposition[see \autoref{proof.52}]\label{prop.52}
For each integer $n>1$, there is a smooth spatial model
$X([4n,0,24])\into\Cp3$ containing $52$ lines, namely, the configuration
\quartic{Z52}.
Denoting by $N(n)$ the number of such models, we have
$\limsup_{n\to\infty}N(n)=\infty$.
\endproposition

The proof of \autoref{prop.52} is based on the fact that the
lines constituting the configuration \quartic{Z52} span a lattice of
rank~$19$. Similar arguments would apply to all sextics in
\autoref{tab.sextics} and most octics in \autoref{tab.octics}, producing
infinitely many smooth models $X\into\Cp4$ or~$\Cp5$ with many
($42$ or $33$, respectively) lines. (In fact, as shown in~\cite{degt:lines},
there also is a $1$-parameter family of octics with $34$
lines.)

We conclude
with an example of a \emph{singular} quartic containing many
lines. At present, the best known bound on the number of lines in a quartic
surface with at worst simple singularities is~$64$, and the best known
example of a singular quartic has $40$ lines (and one
simple node).
Both statements are due to D.~Veniani~\cite{Veniani}.

\theorem[see \autoref{proof.singular}]\label{th.singular}
The $K3$-surface $X([4,0,12])$ has a spatial model with two simple nodes that
contains $52$ lines.
\endtheorem

Recently,
D.~Veniani (private communication) has found an explicit equation of this
quartic:
\[*
z_0z_1(z_0^2+z_1^2+z_2^2+z_3^2)+z_2z_3(z_0^2+z_1^2-z_2^2-z_3^2)=
 2z_2^2z_3^2-2z_0^2z_1^2-2z_0z_1z_2z_3.
\]
He also observed that, when reduced modulo~$5$, the quartic has four simple
nodes and $56$ lines: the best known example in characteristics other
than~$2$ and~$3$.

Remarkably, the discriminant $\det[4,0,12]=48$ equals that of Schur's
quartic.
A few other singular quartics with many lines are also discussed in
\autoref{proof.singular}.

\subsection{Other fields of definition}\label{s.fields}
Given a field $\kk\subset\C$,
one can define the maximal number~$M_\kk$ of lines defined over~$\kk$ in a smooth
quartic $X\subset\Cp3$ defined over~$\kk$.
Thus,
\[*
M_\C=64,\qquad M_\R=56,\qquad M_\Q\le52,
\]
see~\cite{DIS,rams.schuett}.
(Similar numbers can be defined for other polarizations as well, but very little is
known about them.)
The
precise value of~$M_\Q$ was left unsettled in~\cite{DIS}, as in all
interesting examples for which
defining equations are known the lines are only defined
over quadratic algebraic number fields (which can usually be shown by
computing the cross-ratios of appropriate quadruples of intersection points).
We discuss this problem in \autoref{s.Q} and
show that $M_\Q\ge46$.

\theorem[see \autoref{s.Q}]\label{th.Q}
Given a smooth quartic $X\subset\Cp3$ defined over~$\Q$,
denote by $\Fn_\Q X$ the set of lines in~$X$ defined over~$\Q$ and let
$\Fano_\Q(X)\subset\NS(X)$ be the lattice spanned by~$h$ and the
classes~$[\ell]$,
$\ell\in\Fn_\Q X$. Then one has\rom:
\roster
\item\label{i.Q.41}
if $\rank\Fano_\Q(X)=20$, then $\ls|\Fn_\Q X|\le41$, and this bound is
sharp\rom;
\item\label{i.Q.46}
there exists a model \quartic{Q46} of $X([2,1,82])$ with
$\ls|\Fn_\Q\quartic{Q46}|=\ls|\Fn\quartic{Q46}|=46$\rom;
one has
$\rank\Fano_\Q(\quartic{Q46})=19$ and
$\Ps(\quartic{Q46})=
\PS{[ [ 4, 3 ], 2 ], [ [ 3, 5 ], 16 ], [ [ 3, 4 ], 16 ], [ [ 2, 7 ], 12 ]}$.
\endroster
\endtheorem

In the course of the proof of \autoref{th.Q} we also observe that the
$K3$-surface $X([2,1,82])$ has more than three thousands of distinct smooth
spatial models. All models and lines therein are defined over~$\Q$
(see \autoref{s.Q}).

The bound~$M_\kk$ can be defined for any field~$\kk$,
including
the case $\fchar\kk>0$.
We have
\[*
M_{\bar\F_2}=60,\qquad
M_{\bar\F_3}=112,\qquad
M_{\bar\F_p}=64\quad \text{for $p\ge5$},
\]
see~\cite{degt:supersingular}, \cite{rams.schuett:char3},
and~\cite{rams.schuett}, respectively. (In fact, lines in the maximizing
examples
are also defined
over quadratic extensions of~$\F_p$.) Most other questions considered in this
paper also make sense over fields of positive characteristic.
If $\rank\NS(X)\le20$, the lattice $\NS(X)$ lifts to characteristic~$0$ (see,
\eg, \cite{Lieblich.Maulik}) and
appropriate versions of Theorems~\ref{th.others} and~\ref{th.polarizations}
still hold as \emph{upper bounds} on the number of models. We can no longer assert
the existence of each model over each field because of the lack of surjectivity
of the period map. Note also that, in the arithmetical settings, the
transcendental lattice~$T$ is not defined; however, we can still speak about
its genus (given by
$\Gs(T)=4-\rank T=\rank\NS(X)-18$
and $\discr T\cong-\discr\NS(X)$, see \autoref{s.lattice}
and \cite{Nikulin:forms} for details)
and, in particular, the discriminant $\det T=-\det\NS(X)$.
An analogue of
\autoref{cor.Schur} also holds if $\fchar\kk\ge5$.

In the case $\fchar\kk>0$, a more interesting phenomenon is that of
\emph{supersingular} $K3$-surfaces~$X$, \ie, such that $\rank\NS(X)=22$.
(Note, though, that the quartics over $\bar\F_p$, $p\ne3$, maximizing the
number of lines
are \emph{not} supersingular; for example, if $p=2$, the bound for
supersingular quartics is $40$ lines, see~\cite{degt:supersingular}.)
The N\'{e}ron--Severi lattice of
a supersingular $K3$-surface~$X$
is determined by the so-called \emph{Artin invariant}
\[*
\Artin(X):=\frac12\dim_{\F_p}\discr\NS(X)\in\{1,\ldots,10\},
\]
\cf. \autoref{s.lattice} below.
Unfortunately, the approach outlined in \autoref{s.idea}, \viz. embedding the
orthogonal complement $S:=h^\perp\subset\NS(X)$
to a Niemeier lattice, does not work
unless $\Artin(X)=1$. As an alternative, one can probably start from all, not
necessarily smooth, models with $\Artin(X)=1$ and study root-free finite index
sublattices of~$S$. We postpone this question until a future paper.
The intuition
(\cf. the discussion of long vectors in \autoref{s.idea} below)
suggests that there should be
a large number of smooth models whenever $\fchar\kk\gg0$ or
$h^2\gg0$; probably, the two most interesting cases would be
$\fchar\kk=2$ or~$3$
(\cf. \cite[Theorems 1.1, 1.2 and Remark 7.7]{degt:supersingular}).

\subsection{Idea of the proof}\label{s.idea}
Most proofs in the paper
reduce to a detailed study of the lattice $S:=h^\perp\subset\NS(X)$, where
$h$ is the polarization. We can easily control the genus of~$S$; however,
since $S$ is negative definite of rank~$19$, this genus typically consists of
a huge number of isomorphism classes. To list the classes, we represent
$S$
as the orthogonal complement of a certain fixed lattice~$\bbT$ in a
Niemeier lattice~$N$ (see \autoref{s.reduction}).
Certainly, this approach
is not new, \cf., \eg, Kond\=o~\cite{Kondo},
dealing with the Mukai groups, or
Nishiyama~\cite{Nishiyama},
where Jacobian
elliptic $K3$-surfaces are studied by means of the orthogonal complement
$\bU^\perp\subset\NS(X)$ of the sublattice generated by the distinguished
section and a generic fiber.

Typically, there are many isometric embeddings $\bbT\into N$, hence, many
models (\cf. \autoref{th.infty} and \autoref{prop.52}). However, if $\bbT$ is
sufficiently ``small'', all or most orthogonal complements
$S\cong\bbT^\perp\subset N$ contain roots and the corresponding models are
singular; this phenomenon is accountable for the fact that singular
$K3$-surfaces of small discriminant admit very few \emph{smooth} models.
A good quantitative restatement of this
intuitive
observation might be an interesting
lattice theoretical problem shedding more light to the models of
$K3$-surfaces. At present, I can only suggest the computation of the
so-called
\emph{minimal dense square}, covering partially the special case of a single
long vector (see the proof of \autoref{th.Q} in \autoref{s.Q}); most other cases
are handled by a routine \GAP-aided enumeration.

\subsection{Contents of the paper}\label{s.contents}
In \autoref{S.lattice},
we recall briefly a few notions and known results concerning integral
lattices and their extensions, which are the principal technical tools of the
paper. In \autoref{S.S}, Theorems~\ref{th.main} and~\ref{th.others}
are partially reduced to the classification of root-free
lattices in certain fixed genera;
this, in turn,
amounts to the study of appropriate sublattices in the Niemeier lattices.

A detailed proof of \autoref{th.main} is contained in \autoref{S.Niemeier}
(elimination of most Niemeier lattices)
and \autoref{proof.main} (a thorough study of the few lattices left).
A similar, but less detailed, proof of \autoref{th.others} and
\autoref{cor.P1xP1} is outlined in \autoref{S.others}; at the end of this
section, we also discuss Oguiso pairs and Cremona equivalence.
In \autoref{S.polarizations}, we extend the approach to the other
polarizations of $K3$-surfaces, prove \autoref{th.polarizations},
and discuss
smooth polarized $K3$-surfaces that carry a faithful action of
one of the Mukai
groups by symplectic projective automorphisms
(for short, \emph{models of Mukai groups}).
Finally, in \autoref{S.examples}, we consider a few sporadic examples, in
particular those constituting Theorems~\ref{th.infty}, \ref{th.singular},
\ref{th.Q} and \autoref{prop.52}.

\subsection{Acknowledgements}
I am grateful to
S{\l}awomir Rams, Matthias Sch\"{u}tt, Ichiro Shimada,
Tetsuji Shioda, and
Davide Veniani
for a number of fruitful and motivating discussions of the subject.
My special gratitude goes to Dmitrii Pasechnik, who patiently explained to me
the computational aspects of Mathieu groups and Golay codes.
I would also like to thank the anonymous referees of this paper for several
valuable suggestions.
This paper was mainly conceived during my short visit to the
\emph{Leibniz Universit\"{a}t Hannover};
I am grateful to the research unit of the
\emph{Institut f\"{u}r Algebraische Geometrie}
for their warm hospitality.

\section{Lattices}\label{S.lattice}

We recall briefly a few notions and known results concerning integral
lattices and their extensions. Principal references are~\cite{Nikulin:forms}
and~\cite{Conway.Sloane}.

\subsection{Integral lattices\pdfstr{}{
 {\rm(see~\cite{Nikulin:forms})}}}\label{s.lattice}
An \emph{\rom(integral\rom) lattice} is a finitely generated free abelian
group~$L$ equipped with a symmetric bilinear form
\[*
L\otimes L\to\Z,\quad x\otimes y\mapsto x\cdot y.
\]
We abbreviate $x\cdot x=x^2$. In this paper, all lattices are nondegenerate
and \emph{even}, \ie, $x^2=0\bmod2$ for each $x\in L$. The group of
autoisometries of a lattice~$L$ is denoted by $\OG(L)$.

We also consider $\Q$-valued symmetric bilinear forms, possibly
degenerate, on free abelian groups; to avoid confusion, they are referred to as
\emph{forms}. The kernel of a form~$Q$ is the subgroup
\[*
\ker Q=Q^\perp:=\bigl\{x\in Q\bigm|
 \text{$x\cdot y=0$ for each $y\in Q$}\bigr\}.
\]
The quotient $Q/\ker Q$ (often abbreviated to $Q/\ker$) is a nondegenerate form.

An example of a form is the \emph{dual group} $L\dual$ of a lattice~$L$,
\[*
L\dual:=\bigl\{x\in L\otimes\Q\bigm|
 \text{$x\cdot y\in\Z$ for each $y\in L$}\bigr\},
\]
with the bilinear form inherited from $L\otimes\Q$.
We have an obvious inclusion $L\subset L\dual$; the finite quotient group
$\discr L:=L\dual/L$ is called the \emph{discriminant group} of~$L$.
The lattice~$L$ is said to be \emph{unimodular} if $\discr L=0$, \ie,
if $L=L\dual$.
In general, one has $\ell(\discr L)\le\rank L$, where the
\emph{length} $\ell(A)$ of an abelian group~$A$ is
defined as
the minimal number of elements generating~$A$.

The discriminant
group inherits from $L\otimes\Q$ a
symmetric bilinear form
\[*
b\:\discr L\otimes\discr L\to\Q/\Z,\quad
 (x\bmod L)\otimes(y\bmod L)\mapsto (x\cdot y)\bmod\Z,
\]
and its quadratic extension
\[*
q:=q_L\:\discr L\to\Q/2\Z,\quad
 (x\bmod L)\mapsto x^2\bmod2\Z,
\]
called, respectively, the \emph{discriminant bilinear and quadratic forms}.
If there is no confusion, we use the abbreviation $b(\Ga,\Gb)=\Ga\cdot\Gb$ and
$q(\Ga)=\Ga^2$. The discriminant form is
\emph{nondegenerate} in the sense that the
homomorphism
\[*
\discr L\to\Hom(\discr L,\Q/\Z),\quad
\Ga\mapsto(\Gb\mapsto\Ga\cdot\Gb),
\]
is an isomorphism.
When
speaking about (auto-)morphisms of discriminant groups
$\discr L$, {\em the discriminant
forms are always taken into account}.

Note that \emph{not} any sublattice $S\subset L$ is an orthogonal direct
summand. However, if $S$ is nondegenerate, we have well-defined orthogonal
projections
\[
\pr_S\:L\to S\dual,\qquad \pr_S^\perp\:L\to(S^\perp)\dual
\label{eq.projections}
\]
to the dual forms; these projections extend to~$L\dual$.

Lattices are naturally grouped into \emph{genera}. Omitting the precise
definition, we merely use \cite[Corollary 1.9.4]{Nikulin:forms}
which states that
two nondegenerate even lattices $L',L''$
are in the same genus if and only if one has
$\rank L'=\rank L''$, $\Gs(L')=\Gs(L'')$
(where $\sigma(L)$ is the usual signature of
$L\otimes\R$),
and $\discr L'\cong\discr L''$.
Each genus consists
of finitely many isomorphism classes.

\subsection{The homomorphism \pdfstr{O(L) -> Aut discr L}
 {$\OG(L)\to\Aut\discr L$}}\label{s.O(L)}
When speaking about 
isometries of discriminant
groups, we always take~$q$ into account. The group of autoisometries of
$(\discr L,q)$ is denoted by $\Aut\discr L$. The action of $\OG(L)$ extends
to $L\otimes\Q$ by linearity, and the latter extension descends to
$\discr L$. Hence, there is a canonical homomorphism
$\OG(L)\to\Aut\discr L$. In general, this map is neither one-to-one nor onto;
nevertheless, we do not introduce a dedicated notation and freely apply
autoisometries $g\in\OG(L)$ to objects in $\discr L$.
The abbreviation $[\Aut\discr L:\OG(L)]$ stands for the index of
the image of the above canonical homomorphism.
Given an element $\Gg\in\discr L$, we denote by $\stab\Gg\subset\Aut\discr L$
and $\Stab\Gg\subset\OG(L)$ the stabilizer of~$\Gg$ and its pull-back
in $\OG(L)$, respectively.

Let $\CL:=\discr L$ and $A:=\Aut\CL$. Denote by $\<\Gg\>\subset\CL$ the
subgroup generated by an element $\Gg\in\CL$. The restriction $q|_{\<\Gg\>}$
is nondegenerate if and only if the order of $\Gg^2$ in $\Q/\Z$ equals that
of~$\Gg$ in~$\CL$. If this is the case, we have an orthogonal direct
sum decomposition $\CL=\<\Gg\>\oplus\Gg^\perp$ and $\Gg^\perp$ is also
nondegenerate. Fix $s\in\Q/2\Z$ and a nondegenerate quadratic
form~$\CV$ and consider the set
\[*
\CL_{s}(\CV):=\bigl\{\Gg\in\CL\bigm|
 \text{$q|_{\<\Gg\>}$ is nondegenerate, $\Gg^2=s$, $\Gg^\perp\cong\CV$}\}.
\]
It is immediate from the definitions that $\CL_{s}(\CV)$ consists of a single
$A$-orbit and that, for each $\Gg\in\CL_{s}(\CV)$, the stabilizer
$\stab\Gg$ is canonically identified with
the full automorphism group $\Aut\Gg^\perp$.
Hence, denoting by $\Gg\OG(L)$ the orbit, we have
\[
\ls|\CL_s(\CV)|\cdot[\Aut\Gg^\perp:\Stab\Gg]
 =\ls|\Gg\OG(L)|\cdot[\Aut\CL:\OG(L)].
\label{eq.orbits}
\]

\subsection{Extensions\pdfstr{}{
 {\rm(see~\cite{Nikulin:forms})}}}\label{s.extension}
An \emph{extension} of a lattice~$S$ is an even lattice $L\supset S$.
Two extensions $L',L''\supset S$ are \emph{isomorphic} if there is a
bijective
isometry $L'\to L''$ preserving~$S$ \emph{as a set}; they are
\emph{strictly isomorphic} if this isometry can be chosen identical on~$S$.
In the latter case, if $L$ is fixed, we will also speak about the
$\OG(L)$-orbits of isometries $S\into L$.

A \emph{finite index extension} of~$S$ is an even lattice
$L\supset S$ such that $\rank L=\rank S$, \ie, $L$ contains~$S$ as a subgroup
of finite index. An extension gives rise to an
isometry $L\into S\otimes\Q$ and, since
$L$ is also a lattice, we have $L\subset S\dual$.
Furthermore, the subgroup $\CK:=L/S\subset\discr S=S\dual\!/S$, called the
\emph{kernel} of the extension,
is \emph{isotropic}, \ie, $q|_\CK=0$. Conversely, if
$\CK\subset\discr S$ is isotropic, the group
\[*
L:=\bigl\{x\in S\otimes\Q\bigm|x\bmod S\in\CK\bigr\}
\]
is an integral lattice. This can be summarized in the following statement.

\proposition[see~\cite{Nikulin:forms}]\label{prop.extension}
Given a lattice~$S$, the correspondence
\[*
L\mapsto\CK:=L/S,\qquad
\CK\mapsto L:=\bigl\{x\in S\otimes\Q\bigm|x\bmod S\in\CK\bigr\}
\]
is a bijection between the set of
strict isomorphism classes of finite index extensions $L\supset S$
and that of isotropic subgroups $\CK\subset\discr S$.
Under this correspondence, one has $\discr L=\CK^\perp\!/\CK$.
Furthermore,
\roster
\item\label{i.prim.iso}
two extensions~$L'$, $L''$ are isomorphic if and only if
their kernels~$\CK'$, $\CK''$ are in the same
$\OG(S)$-orbit, \ie, there is $g\in\OG(S)$ such that $g(\CK')=\CK''$\rom;
\item\label{i.prim.aut}
an autoisometry $g\in\OG(S)$ extends to~$L$ if and only if $g(\CK)=\CK$.
\done
\endroster
\endproposition

A \emph{primitive extension} is an extension $L\supset S$ such that $S$ is
primitive in~$L$, \ie, $(S\otimes\Q)\cap L=S$.
In~\cite{Nikulin:forms}, such extensions are studied by fixing  (the
isomorphism class of) the orthogonal complement $T:=S^\perp\subset L$.
Then, $L$ is a finite index extension of $S\oplus T$ in which both~$S$
and~$T$ are primitive. According to \autoref{prop.extension}, this extension
is described by its kernel
\[*
\CK\subset\discr(S\oplus T)=\discr S\oplus\discr T,
\]
and the primitivity of~$S$ and~$T$ in~$L$ implies that
\[*
\CK\cap\discr S=\CK\cap\discr T=0.
\]
In other words, $\CK$ is the graph of a certain monomorphism
$\psi\:\Cal{D}\to\discr T$, where $\Cal{D}\subset\discr S$;
since $\CK$ is isotropic, $\psi$ is an anti-isometry.

To keep track of the sublattice $S\subset L$,
Statements~\iref{i.prim.iso} and~\iref{i.prim.aut} of
\autoref{prop.extension} should be restricted to the subgroup
$\OG(S)\times\OG(T)\subset\OG(S\oplus T)$.
Below, we use freely a number of other similar restrictions
taking into account additional structures.

An extension $L\supset S\oplus T$ is unimodular if and only if
$\CK^\perp\!/\CK=0$, \ie, $\CK^\perp=\CK$. Then,
$\ls|\CK|^2=\ls|\discr S|\cdot\ls|\discr T|$, which implies that
$\ls|\CK|=\ls|\discr S|=\ls|\discr T|$ and
$\psi$ above is a group isomorphism $\discr S\to\discr T$.

\corollary[see~\cite{Nikulin:forms}]\label{cor.extension}
Given a pair of lattices~$S$, $T$, there is a natural one-to-one
correspondence between the strict isomorphism classes of unimodular
finite index extensions $N\supset S\oplus T$
in which both~$S$ and~$T$ are primitive and bijective isometries
$\psi\:\discr S\to-\discr T$.
If an isometry~$\psi$
\rom(hence, an extension~$N$\rom) is fixed,
then
\roster
\item\label{i.ext.iso}
the isomorphism classes of extensions are in a one-to-one correspondence with
the double cosets $\OG(S)\backslash\Aut\discr S/\!\OG(T)$\rom;
\item\label{i.ext.aut}
a pair of autoisometries $g\in\OG(S)$, $h\in\OG(T)$ extends to~$N$ if and only
if one has $g=\psi\1h\psi$ in $\Aut\discr S$.
\done
\endroster
\endcorollary

\subsection{Lattices of rank~$2$}\label{s.2x2}
According to Gauss~\cite{Gauss:Disquisitiones}, any positive definite
even lattice $T$ of rank~$2$ has a basis in which the Gram matrix
of~$T$ is of the form
\[*
\bmatrix a&b\\b&c\endbmatrix,\qquad
0<a\le c,\quad 0\le2b\le a,\quad a=c=0\bmod2,
\]
and the ordered triple $(a,b,c)$ satisfying the above conditions is uniquely
determined by~$T$. We use the notation $[a,b,c]$ for the isomorphism class
of~$T$. This description implies that $\frac34c\le\det T\le c^2$, which
makes the enumeration of lattices within a given genus an easy task.

Let $\OG^+(T)\subset\OG(T)$ be the group of orientation preserving
autoisometries of~$T$.
Using the description above, one can see that, with few exceptions,
$\OG^+(T)=\{\pm\id\}$.
The exceptions are the lattices $T=[2n,0,2n]$ (with $\OG^+(T)\cong\Z/4$)
and $[2n,n,2n]$ (with $\OG^+(T)\cong\Z/6$), where $n$ is any positive integer.
Since the greatest common divisor of the entries of a Gram matrix is a genus
invariant, each of the exceptional lattices is unique in its genus.
It follows that, for a positive definite even lattice~$T$ of rank~$2$, the image of
$\OG^+(T)$ in $\Aut\discr T$ depends on the genus of~$T$
only.

\subsection{Root systems\pdfstr{}{
 {\rm(see~\cite[Chapter 4]{Conway.Sloane})}}}\label{s.root}
A \emph{root system}, or \emph{root lattice}, is a positive definite
lattice~$R$
generated by its \emph{roots}, \ie, vectors $a\in R$ of square~$2$. (Recall
that we consider even lattices only.) Any positive definite lattice~$L$
contains its maximal root lattice $\rt(L)$, which is generated by all roots
$a\in L$.

Each root lattice decomposes uniquely into an orthogonal direct sum of
irreducible ones, which are of type $\bA_n$, $n\ge1$, $\bD_n$, $n\ge4$, or
$\bE_n$, $n=6,7,8$. One has
\[*
\discr\bA_n\cong\Z/(n+1),\quad
\ls|\discr\bD_n|=4,\quad
\ls|\discr\bE_n|=9-n.
\]
For these groups, we use the numbering $\discr R=\{0=\Ga_0,\Ga_1,\ldots\}$
as in~\cite{Conway.Sloane}.

The lattices $\bA_{n-1}$ and $\bD_n$ are, respectively, the
orthogonal complement and $(\bmod\,2)$-orthogonal complement of the
characteristic vector $\bar\be:=\be_1+\ldots+\be_n$ in the odd unimodular lattice
$\bH_n:=\bigoplus_{i=1}^n\Z\be_n$, $\be_i^2=1$.
Then, $\bE_8\supset\bD_8$ and $\bE_7\supset\bA_7$ are the index~$2$
extensions by the vector $\frac12\bar\be-4\be_8$.
Alternatively, the lattices~$\bE_n$, $n=6,7$, can be described as
$\bE_n=\bA_{8-n}^\perp\subset\bE_8$.

If $n$ is large, ``short'' vectors in $\bH_n\otimes\Q$ tend
to have many equal coordinates, and we follow~\cite{Conway.Sloane} and use
the ``run-length encoding'' for the $\SG{n}$-orbits of such vectors: the
notation
\[*
(s_1)^{u_1}\ldots(s_t)^{u_t},\quad s_1<\ldots<s_t,\quad u_i>0,\quad
 u_1+\ldots+u_t=n
\]
designates the orbit whose representatives have $u_i$ coordinates~$s_i$,
$i=1,\ldots,t$.

In this notation, the shortest representatives of the nonzero elements of the
discriminant groups are as follows
(for $\bE_6$ and $\bE_7$, we only indicate the squares):
\begin{gather}
  \bA_{n}\:\quad\Ga_p=\frac1{n+1}\bigl((-p)^q(q)^p\bigr),\quad
    \Ga_p^2=\frac{pq}{n+1}\quad (q:=n+1-p);
   \label{eq.shortest.A}\\\allowbreak
  \bD_n\:\quad
   \gathered
      \Ga_{2k+1}=\frac12\bigl((-1)^p(1)^{n-p}\bigr),\quad
        \Ga_{2k+1}^2=\frac{n}4\quad (p=k\bmod2), \\
        \Ga_2=\bigl((0)^{n-1}(\pm1)^1\bigr),\quad \Ga_2^2=1;
   \endgathered
   \label{eq.shortest.D}\\\allowbreak
  \bE_6\:\quad \Ga_1^2=\Ga_2^2=\frac43;\qquad
  \bE_7\:\quad \Ga_1^2=\frac32.
   \label{eq.shortest.E}
\end{gather}

The groups $\OG(\bA_{n-1})$ and $\OG(\bD_n)$ are semi-direct products
${\GROUP{RG}}\rtimes\SG{n}$, where $\GROUP{RG}$ is generated by $-\id$ and,
for~$\bD_n$, by the reflections against all basis elements~$\be_i$.
If $n$ is large, \GAP's built-in orbit/stabilizer routines do not work very
well, and we use the run-length encoding to classify the pairs and triples of
vectors (\cf. \autoref{s.computation} below).
More precisely, given two $\SG{n}$-orbits encoded by
$(s_1^*)^{u_1^*}\ldots(s_{t^*}^*)^{u_{t^*}^*}$, $*=\prime$ or $\prime\prime$,
the $\SG{n}$-orbits of pairs of vectors can be encoded by the sequences of
the form
\[*
\ldots(s_i',s_j'')^{u_{ij}}\ldots,\quad 1\le i\le t',\quad 1\le j\le t''
\]
(after disregarding the entries with $u_{ij}=0$), where $u_{ij}\ge0$ are
integers such that
$\sum_ju_{ij}=u_i'$ for each~$i$ and $\sum_iu_{ij}=u_j''$ for each~$j$.
The classification of such balanced matrices $[u_{ij}]$ is straightforward.
The further passage from pairs to triples of vectors is done
in a similar way.

\subsection{The Niemeier lattices\pdfstr{}{
 {\rm(see~\cite[Chapter 16]{Conway.Sloane} and~\cite{Niemeier})}}}\label{s.Niemeier}
A \emph{Niemeier lattice} is
a positive definite unimodular even lattice of
rank~$24$. Up to isomorphism, there are $24$ Niemeier lattices.
One of them,
the so-called \emph{Leech lattice}, has no roots, and each of the other
$23$ lattices~$N$ is a finite index extension of $\rt(N)$.
Furthermore, the isomorphism class of~$N$ is determined by that of $\rt(N)$,
see \autoref{tab.Niemeier}
(where we also refer to the relevant parts of the proof of
\autoref{th.main});
for this reason, the Niemeier lattice with
$\rt(N)=R$ is often denoted by $N(R)$.
\table
\caption{The $24$ Niemeier lattices}\label{tab.Niemeier}
\let\+\oplus
\def\balign{\vtop\bgroup\halign\bgroup\strut\quad$##$\hss\quad&\quad##\hss\quad\cr
\noalign{\hrule\vspace{2pt}}%
\text{Roots}&Details\cr
\noalign{\vspace{1pt}\hrule\vspace{2pt}}}
\hbox to\hsize\bgroup\hss\balign
\bD_{24}&see \autoref{ss.D24}\cr
\bD_{16}\+\bE_8&see \autoref{ss.E8}\cr
\bE_8^3&see \autoref{ss.E8}\cr
\bA_{24}&see \autoref{ss.D24}\cr
\bD_{12}^2&see \autoref{ss.D12}\cr
\bA_{17}\+\bE_7&see \autoref{ss.E7}\cr
\bD_{10}\+\bE_7^2&see \autoref{ss.E7}\cr
\bA_{15}\+\bD_9&see \autoref{ss.A15}\cr
\bD_8^3&see \autoref{ss.n<=6}\cr
\bA_{12}^2&see \autoref{ss.D12}\cr
\bA_{11}\+\bD_7\+\bE_6&see \autoref{ss.E6}\cr
\bE_6^4&see \autoref{ss.E6}\cr
\ealign\qquad\balign
\bA_9^2\+\bD_6&see \autoref{ss.n<=6}\cr
\bD_6^4&see \autoref{ss.n<=6}\cr
\bA_8^3&see \autoref{ss.n<=6}\cr
\bA_7^2\+\bD_5^2&see \autoref{ss.n<=6}\cr
\bA_6^4&see \autoref{ss.n<=6}\cr
\bA_5^4\+\bD_4&see \autoref{ss.n<=6}\cr
\bD_4^6&see \autoref{ss.n<=6}\cr
\bA_4^6&see \autoref{ss.n<=6}\cr
\bA_3^8&see \autoref{s.8A3}\cr
\bA_2^{12}&see \autoref{s.12A2}\cr
\bA_1^{24}&see \autoref{s.24A1}\cr
\text{Leech}&no roots\cr
\ealign\hss\egroup
\endtable

Consider two positive definite lattices $S$, $V$ and assume that
$\rank S+\rank V=24$ and $\discr S\cong-\discr V$, so that $V\oplus S$
admits a finite index extension to a Niemeier lattice, see
\autoref{cor.extension}.

\lemma\label{lem.index}
Given~$S,V$ as above, the index $[\Aut\discr S:\OG(S)]$ equals the
number of strict isomorphism classes of extensions $N\supset V$, to
\emph{all} Niemeier lattices $N$, with the
property that $V^\perp\cong S$.
\endlemma

\proof
Let $A:=\Aut\discr S$, and let $H,G\subset A$ be the images of $\OG(S)$ and
$\OG(V)$, respectively. (For the latter, we fix an anti-isometry
$\psi\:\discr S\to\discr V$.)
The statement of the lemma follows from the \emph{double coset
formula}
\[*
\ls|A/G\times A/H|=\sum_{g\in G\backslash A/H}\ls|A/(G^g\cap H)|,
\]
where we let $G^g:=g\1Gg$. Indeed, the formula
implies
$[A:H]=\sum[G^g:G^g\cap H]$,
the summation running over all
double cosets $g\in G\backslash A/H$, \ie, over all isomorphism classes of
extensions $N\supset V\oplus S$, see \autoref{cor.extension}.
By the same corollary, $G^g\cap H$ is the group of autoisometries of~$V$
extending to~$N$; hence, the
contribution of a class $g$ in the above sum is the number of
strict isomorphism classes of extensions $N\supset V$ that are contained in~$g$.
\endproof

\section{The reduction}\label{S.S}

Theorems~\ref{th.main} and~\ref{th.others}
can partially be reduced to the classification and study of root-free
lattices in certain fixed genera, see \autoref{th.K3}.
This, in turn,
amounts to the study of appropriate sublattices in the Niemeier lattices, see
\autoref{lem.genus}.

\subsection{Quartic $K3$-surfaces}\label{s.K3}
Consider a $K3$-surface~$X$, and let $L:=H_2(X)$ be its homology group,
equipped with the intersection paring. It is a unimodular even lattice of
signature $(3,19)$; such a lattice is unique up to isomorphism.

For the sake of simplicity, we confine ourselves to the \emph{singular}
$K3$-surfaces, \ie, we assume that
$\rank\NS(X)=20$.
A singular $K3$-surface~$X$ is
characterized by the \emph{oriented} isomorphism type of its
\emph{transcendental lattice}
\[*
T:=\NS(X)^\perp\subset L=H_2(X),
\]
which is a positive definite
even lattice of rank~$2$,
see~\cite{Shioda.Inose}.
(The vector space $T\otimes\R$ is spanned by the real and imaginary parts of
a holomorphic
$2$-form on~$X$, and this basis defines a distinguished orientation.)
This correspondence is emphasized by the notation $X:=X(T)$.
The $K3$-surface $X(\bar T)$ corresponding to the lattice~$T$ with the
opposite orientation is the complex conjugate surface
\smash{$\,\overline{\!X(T)\!}\,$}. The surface $X(T)$ is \emph{real} (defined
over~$\R$) if and only if $T$ has an orientation reversing automorphism
(which, in rank~$2$, is always involutive).

A \emph{spatial model} of
a singular $K3$-surface~$X$ is a map $\Gf\:X\to\Cp3$ defined by a
fixed point free ample
linear system of degree~$4$. Two models~$\Gf_1$, $\Gf_2$ are
\emph{projectively equivalent} if there exist automorphisms $a\:\Cp3\to\Cp3$
and $a_X\:X\to X$ such that $\Gf_2\circ a_X=a\circ\Gf_1$.
Denote by $h:=h_\Gf\in\NS(X)$, $h^2=4$,
the class of
a hyperplane section and let
$S_\Gf:=h^\perp\subset\NS(X)$; it is an even negative definite lattice of
rank~$19$. We can represent~$S_\Gf$ as the orthogonal complement
$(T\oplus\Z h)^\perp\subset L$. The new sublattice $T\oplus\Z h\subset L$ is not
necessarily primitive, and its primitive hull is the finite index extension
determined by a certain isotropic subgroup
\[
\tker\subset\discr T\oplus\discr\Z h,\qquad\tker\cap\discr T=0.
\label{eq.C-kernel}
\]
(The last identity is due to the fact that $T$ \emph{is} primitive in~$L$.)
This subgroup~$\tker$ is cyclic of order~$1$, $2$, or~$4$, and we define
the \emph{depth} of~$\Gf$ as $\depth\Gf:=\ls|\tker|$.
Then, by \autoref{cor.extension},
\[
\discr S_\Gf\cong-\tker^\perp\!/\tker,\qquad
\ls|\discr S_\Gf|=4\ls|\discr T|/(\depth\Gf)^2.
\label{eq.depth}
\]

Fix a lattice $S:=S_\Gf$ and let $d=d(S):=\depth\Gf$; this number is recovered
from~$S$ and~$T$ \via~\eqref{eq.depth}.
Consider the set
\[*
\CS_{dh}:=\bigl\{\Gg\in\discr S\bigm|
 (4/d)\Gg=0,\ \Gg^2=-d^2\!/4\bmod2\Z\bigr\}.
\]
Each element $\Gg\in\CS_{dh}$ gives rise to the isotropic subgroup
$\CK_\Gg\subset\discr\Z h\oplus\discr S$ generated by $(d/4)h\oplus\Gg$,
and we define
\[*
\CS_{dh}^+:=\bigl\{\Gg\in\CS_{dh}\bigm|
 \CK_\Gg^\perp/\CK_\Gg\cong-\discr T\bigr\}.
\]
Clearly,
$\CK_\Gg$ is the kernel of the extension $\NS(X)\supset\Z h\oplus S$;
one has $\ls|\CK_\Gg|=4/d$.
Recall that $\stab\Gg\subset\Aut\discr S$
and
$\Stab\Gg\subset\OG(S)$
are the stabilizers of an element $\Gg\in\CS_{dh}^+$.
Fixing an isometry
$\CK_\Gg^\perp/\CK_\Gg\cong-\discr T$,
we can regard both groups acting
on the discriminant $\discr T$.

\remark\label{rem.S+}
We always have $\CS_{4h}^+=\CS_{4h}=\{0\}$ and $\CK_0^\perp\!/\CK_0=\discr S$.
If $d=1$ or~$2$, the inclusion $\CS_{dh}^+\subset\CS_{dh}$ may be proper. If
$d=1$, then
$\CK_\Gg^\perp\!/\CK_\Gg=\Gg^\perp\subset\CS$,
where $\CS:=\discr S$;
hence, in this case,
the set $\CS_h=\CS_{-1/4}(-\discr T)$
is a single orbit of $\Aut\discr S$ and we have
$\stab\Gg=\Aut\Gg^\perp\cong\Aut\discr T$,
see \autoref{s.O(L)}.
\endremark

\theorem\label{th.K3}
The projective equivalence classes of spatial models
$\Gf\:X(T)\to\Cp3$
are in a one-to-one correspondence with the triples
consisting of
\roster*
\item
a negative definite lattice~$S$ of rank~$19$ and
$\discr S\cong-\tker^\perp\!/\tker$ as in~\eqref{eq.C-kernel},
\item
an $\OG(S)$-orbit $[\Gg]\subset\CS_{dh}^+$
\rom(where $d=d(S)$ is the depth, see~\eqref{eq.depth}\rom), and
\item
a double coset $c\in\OG^+(T)\backslash\Aut\discr T/\Stab\Gg$
\endroster
and such that
\roster
\item\label{i.model.moved}
$d\ge2$ or $d=1$ and $\Gg$ is not represented by a vector
$a\in S\dual$, $a^2=-\frac14$.
\endroster
Under this correspondence, the following statements hold\rom:
\roster[\lastitem]
\item\label{i.model.birational}
a model~$\Gf$ is birational onto its image if and only if either $d=4$ or $d\le2$
and the class $(2/d)\Gg$ is not represented by a $(-1)$-vector
in~$S\dual$\rom;
\item\label{i.model.nonsingular}
a birational model~$\Gf$ is smooth if and only if
$S$ is root free\rom;
\item\label{i.model.lines}
the straight lines contained in the image of a smooth model~$\Gf$ are in a
one-to-one correspondence with the vectors $a\in S\dual$, $a^2=-\frac94$,
representing~$\Gg$.
\endroster
\endtheorem

\remark\label{rem.2=>1}
Note that the hypothesis of \autoref{th.K3}\iref{i.model.birational},
which is equivalent to the
absence of a vector
$e\in\NS(X)$ such that $e^2=0$ and $e\cdot h=2$, implies
Condition~\iref{i.model.moved} in the theorem.
Since we are mainly interested in birational spatial models, we will only
check the hypothesis of~\iref{i.model.birational}.
\endremark

\remark\label{rem.genus.constant}
The number of projective equivalence classes of models and their properties
(birational/smooth/number of lines and their adjacency graph) depend on the
genus of the transcendental lattice~$T$ only.
The last statement follows directly from
\autoref{th.K3}\iref{i.model.moved}--\iref{i.model.lines},
and for the number of models one uses, in addition,
the description of the group $\OG^+(T)$ in
\autoref{s.2x2} to conclude that, with $(S,\Gg)$ fixed,
the quotient set $\OG^+(T)\backslash\Aut\discr T/\Stab\Gg$
is independent of~$T$.

This phenomenon has a simple geometric explanation: if $T'$ and $T''$
are in the same genus, the
corresponding $K3$-surfaces $X(T')$ and $X(T'')$
are Galois conjugate over some algebraic number field
(see, \eg, \cite{Schutt.singular.K3}),
and so
are their spatial models.
\endremark

\proof[Proof of \autoref{th.K3}]
Let $X:=X(T)$.
As explained above in this section, a spatial model gives rise to a lattice
$S:=h^\perp\subset\NS(X)$, a class $\Gg$ which defines the extension
$\NS(X)\supset\Z h\oplus S$, and a double coset~$c$ defining the
extension
$L\supset T\oplus\NS(X)$, see \autoref{s.extension}. A projective equivalence
induces an autoisometry of $H_2(X)$ preserving the pair $h\in\NS(X)$
and the orientation of $T$;
hence, these data are defined up to the
actions listed in the statement.
Conversely, a set of data as in the statement determines the extensions
$H_2(X)\cong L\supset\NS\ni h$ uniquely up to automorphism of~$L$ preserving
the orientation of~$\NS^\perp$.
Multiplying, if necessary, by~$(-1)$ and applying reflections, we can assume
that the \emph{marking} $L\cong H_2(X)$ is chosen so that $\NS$ is taken to
$\NS(X)$ and $h$ is taken to the closure of the K\"{a}hler cone, so that $h$ is
nef. Condition~\iref{i.model.moved}, stating that there is no vector
$e\in\NS(X)$ such that $e^2=0$ and $e\cdot h=1$,
is equivalent to the requirement that the linear system $\ls|h|$ has no fixed
components, see~\cite{Nikulin:Weil}.
Then, by
\cite[Corollary~3.2]{Saint-Donat},
the system $\ls|h|$ is fixed point free
and $\dim\ls|h|=3$; hence,
$h$ does define a spatial model $\Gf\:X\to\Cp3$.

Statement~\iref{i.model.birational}
follows from \cite[Theorem 5.2]{Saint-Donat}.

Statement~\iref{i.model.nonsingular} is well known. By the Riemann--Roch
theorem, if $a\in\NS(X)$ and $a^2=-2$, then either~$a$ or $-a$ is effective.
If also $a\cdot h=0$, then $a$ represents a curve of projective degree~$0$;
this curve is contracted by~$\Gf$. Conversely, by the adjunction formula, the
genus of a curve~$C$ in a $K3$-surface is given by $\frac12C^2+1$.
It follows that all exceptional divisors are rational $(-2)$-curves (in
particular, all singularities are simple);
clearly, the classes of these curves are orthogonal to~$h$.


Statement~\iref{i.model.lines} is also known, see, \eg,~\cite{DIS}. If
$\Gf(X)$ is smooth, then the classes of lines are the vectors
$l\in\NS(X)$ such that $l^2=-2$ and $l\cdot h=1$. Any such vector is of the
form $\frac14h\oplus a$, where $a\in S\dual$ is as in the statement.
\endproof

\subsection{The pencil structure}\label{s.pencils}
Most quartics considered in this paper contain many lines. A convenient way
to identify/distinguish such quartics is the so-called \emph{pencil structure}, which
is an easily computable projective invariant.

Fix a smooth quartic $X\subset\Cp3$. Given a line~$\ell\subset X$, we can
consider the pencil $\{\pi_t,t\in\Cp1\}$ of planes containing~$\ell$. Each
intersection $\pi_t\cap X$ is a planar quartic curve which splits into~$\ell$
itself and the \emph{residual cubic} $C_t\subset\pi_t$.
All but finitely many residual cubics are irreducible; a certain number~$p$
of them split into three lines, and a certain number~$q$ split into a line
and an irreducible conic. The pair $(p,q)$ is called the \emph{type} of the
original line~$\ell$, and the \emph{pencil structure} of~$X$ is the multiset
of the types of all lines contained in~$X$. Following~\cite{DIS}, we use the
partition notation
\[*
\Ps(X)=(p_1,q_1)^{u_1}\ldots(p_r,q_r)^{u_r},
\]
meaning that $X$ contains $u_i$ lines of type $(p_i,q_i)$, $i=1,\ldots r$.
(The total number of lines in~$X$ equals $u_1+\ldots+u_r$; this
number is used as the subscript in the notation for quartics and
configurations of lines.)

The pencil structure is easily computed in terms of the adjacency matrix of
lines. If lines are represented as in \autoref{th.K3}\iref{i.model.lines}, by
vectors $a_i\in S\dual$, $a_i^2=-\frac94$, then the adjacency matrix is the
Gram matrix of $\{a_i\}$ with each entry increased by~$\frac14$.

If $X$ is singular, the types of the fibers~$\pi_t$ are much more diverse
(see \cite{Veniani}) and we do not use the notion of pencil structure.

\subsection{Reduction to Niemeier lattices}\label{s.reduction}
\autoref{th.K3} reduces the classification of (smooth) spatial models of a
fixed $K3$-surface $X(T)$ to that of (root-free)
definite even lattices~$S$ within a certain
collection of genera (which is determined,
\via~\eqref{eq.C-kernel} and~\eqref{eq.depth}, by the genus of~$T$).
For the convenience of the further exposition, we will switch to the positive
definite lattice $-S$, so that $\discr(-S)=\tker^\perp\!/\tker$.

By the construction, $\ell(\tker^\perp\!/\tker)\le\rank T+1=3$. Hence, using
\cite[Theorem 1.10.1]{Nikulin:forms}, one can find a positive definite
lattice~$V$ of rank~$5$ such that $\discr V\cong-\tker^\perp\!/\tker$. (The
other conditions of the theorem hold trivially due to the additivity of the
signature.)
Fix one such lattice~$V$ and call it the \emph{test lattice}.
Then, according to \autoref{cor.extension}, the direct sum
$V\oplus{-S}$ has a unimodular finite index extension in which both $V$ and
$-S$ are primitive.
Since $\rank V+\rank S=24$, this extension is a Niemeier lattice.
Thus, we have the following statement.

\lemma\label{lem.genus}
Any lattice~$S$ in \autoref{th.K3} is of the form $-V^\perp$, where $V$ is
any fixed test lattice and $V\into N$ is a primitive embedding to
a Niemeier lattice.
\done
\endlemma

The finite index extension $N\supset(V\oplus{-S})$ depends on
certain extra data (see \autoref{cor.extension}).
Hence,
\latin{a priori},
distinct embeddings $V\into N'$, $V\into N''$ may result in
isomorphic orthogonal complements $S=-V^\perp$, \cf. \autoref{lem.index}.

\remark\label{rem.genus}
In \autoref{lem.genus},
one can choose a ``universal'' test lattice~$V$ satisfying
\[*
\discr V\cong{-\discr T}\oplus{-\discr\Z h}
\]
and depending on the genus of~$T$ only. Then, lifting the primitivity
requirement, one should check an analogue of~\eqref{eq.C-kernel} for each embedding
$V\into N$.
Note also that one can always assume that $V$ has at least one root
and, thus, exclude the Leech lattice in \autoref{lem.genus}
(\cf. Kond\=o~\cite{Kondo});
in general, one cannot assert that $\rank\rt(V)\ge2$
(\cf. \autoref{conv.V} below;
this is the reason for excluding the lattice
$T=[4,0,16]$ from the statement of \autoref{th.others}).
\endremark

\subsection{The Fermat quartic: the genus $\genus$}\label{s.genus}
The abstract Fermat quartic~$\Fermat$ is characterized by the transcendental
lattice $\TL:=[8,0,8]$, see, \eg,~\cite{Shioda.Inose}.
Hence, it is immediate that any spatial model of~$\Fermat$ has depth~$1$ and
one has
\[*
\discr S\cong{-\discr\TL}\oplus{-\discr\Z h}
\]
in \autoref{th.K3}. Passing, as above, to $-S$, we define
the \emph{genus~$\genus$} as the set of isomorphism classes of positive
definite even lattices of rank~$19$ and discriminant
$\discr\TL\oplus\discr\Z h$, $h^2=4$.

To use \autoref{lem.genus},
consider the test lattice~$\bbT$ given,
in a certain distinguished basis
$\bba_2^1$, $\bba_2^2$, $\bba_4$, $\bbc_4$, $\bba_8$,
by the Gram matrix
\[*
\bbT:=\begin{bmatrix}
        2 & 0 & 0 & 1 & 0 \\
        0 & 2 & 0 & 1 & 0 \\
        0 & 0 & 4 & 2 & 0 \\
        1 & 1 & 2 & 4 & 0 \\
        0 & 0 & 0 & 0 & 8
      \end{bmatrix}.
\]
One has $\rt(\bbT)=\Z\bba_2^1\oplus\Z\bba_2^2\cong\bA_1^2$, and
$\rt(\bbT)^\perp$ is the diagonal lattice
\[*
\bT:=\Z\ba_8\oplus\Z\bc_8\oplus\Z\ba_4,
 \quad (\ba_8)^2=(\bc_8)^2=8,\quad (\ba_4)^2=4,
\]
where
\[
\ba_4=\bba_4,\quad\ba_8=\bba_8,\quad\bc_8=2\bbc_4-\bba_4-\bba_2^1-\bba_2^2.
\label{eq.basis}
\]
The sublattice $\bT\oplus\Z\bba_2^1\oplus\Z\bba_2^2\subset\bbT$ has
index~$2$, since
$\tfrac12(\bc_8+\ba_4+\bba_2^1+\bba_2^2)\in\bbT$.
Therefore, \autoref{lem.genus} can be restated in the following form.

\lemma\label{lem.genus.S.prime}
One has $\genus=\bigl\{(\bT\oplus\Z r_1\oplus\Z r_2)^\perp\bigr\}$,
where
$\bT\into N$ runs through all primitive embedding to Niemeier
lattices~$N$ and $r_1,r_2\in\rt(\bT^\perp)$ run through pairs of orthogonal roots such
that
$\tfrac12(\bc_8+\ba_4+r_1+r_2)\in N$.
\done
\endlemma

In practice, we usually simplify \autoref{lem.genus.S.prime} even further and
merely analyze a triple $\ba_8,\bc_8,\ba_4\in N$ of pairwise orthogonal
primitive vectors of squares $8$, $8$, and~$4$, respectively.
Note that $\ba_4$ is always primitive, whereas $\ba_8$ or $\bc_8$ is
imprimitive if and only the vector equals $2r$ for a root $r\in N$.

\section{Eliminating Niemeier lattices}\label{S.Niemeier}

The principal result of this section is \autoref{th.2xS} in
\autoref{s.classification}. It is proved by eliminating the Niemeier lattices
one-by-one, mainly using \autoref{lem.genus.S.prime}.

\subsection{Projections to the components}\label{s.proj}
Fix a Niemeier lattice~$N$ and let
\[*
\rt(N)=\bigoplus R_k,\quad k\in\IS:=\{1,\ldots,m\},
\]
be the decomposition of $\rt(N)$ into irreducible components.
For a subset $I\subset\IS$ of the index set, we define
$R_I:=\bigoplus_{k\in I}R_k$.
We use repeatedly the following statement, which is, essentially, the
definition of $\rt(N)$.

\lemma\label{lem.roots}
If $r\in N$ is a root, then $r\in R_k$ for some index $k\in\IS$.
\done
\endlemma

Consider the orthogonal projections
\[*
\pr_k\:N\to R_k\dual,\quad k\in\IS,\quad\text{and}\quad
\pr_I\:N\to R_I\dual,\quad I\subset\IS,
\]
see~\eqref{eq.projections}.
In most cases, we will analyze a sublattice
$\bT\subset N\subset\bigoplus_{k}R_k\dual$ by means of its images
$\pr_k(\bT)$. More precisely, we will speak about isometries
$\T_k\to R_k\dual$, not necessarily
injective,
of $\Q$-valued $3$-forms \emph{with distinguished
bases} $\ba_8,\bc_8,\ba_4$.
By means of these bases, we can identify a form and its Gram matrix and
consider sums and differences of forms; thus,
$\pr_I(\bT)=\sum_{k\in I}V_k$ for a subset $I\subset\IS$,
and we must have
$\bT=\sum_{k\in\IS}\T_k$.
The \emph{orthogonal complement} $\T^\perp\subset R\dual$ of
an isometry $\T\to R\dual$ is the orthogonal complement of its image.

\definition\label{def.dense}
Let $R$ be a direct summand of $\rt(N)$, \ie,
$R=R_I$ for some $I\subset\IS$.
A $3$-form $\T\to R\dual$ is said to be
\emph{bounded}, $\T\le\bT$, if the difference $\bT-\T$ is positive
semi-definite. The form $\T\to R\dual$ is \emph{$n$-dense}, $n=0,1,2$, if there are $n$
pairwise orthogonal roots $r_i\in \T^\perp\cap R$ such that
\roster
\item\label{i.dense.roots}
the orthogonal complement $\bigl(\T\oplus\bigoplus_i\Z r_i\bigr)^\perp$
has no roots in~$R$, and
\item\label{i.dense.impr}
the vector $\bc_8+\ba_4+\sum_ir_i$ is divisible by~$2$ in $R\dual$.
\endroster
We will denote by $\dense_n(R)$, $n=0,1,2$, the
set of the Gram matrices of bounded $n$-dense $3$-forms
$\T\to R\dual$ satisfying the additional primitivity condition
\roster[\lastitem]
\item\label{i.dense.8}
if $\ba_8$ or $\bc_8$ projects to a square~$8$ vector in a single summand
$R_k\dual$,
then the image is \emph{not} of the form $2r$, $r\in R_k$.
\endroster
Given $\T\in\dense_nR$, the \emph{reduced complement} $\red_n\T$ is the
\emph{abstract} isomorphism class of the $\Q$-valued quadratic form obtained
as follows: start with $(\bT-\T)/\ker$ and, if $n=2$, pass to its extension
\via\ $\frac12(\bc_8+\ba_4)$.

As usual, we let $\dense_*(R):=\bigcup_n\dense_n(R)$, $n=0,1,2$.
\enddefinition

Clearly, $n$ roots~$r_i$ as in \autoref{def.dense}\iref{i.dense.roots}
exist if and only if the root lattice $R':=\rt(\T^\perp\cap R)$ is isomorphic to
\roster*
\item
$0$, if $n=0$,
\item
$\bA_1$ or $\bA_2$, if $n=1$, or
\item
$\bA_1^2$, $\bA_2\oplus\bA_1$, $\bA_2^2$, $\bA_3$, or $\bA_4$, if $n=2$;
\endroster
in particular, $\rank R'\le2n$.
Then, such a collection $\{r_i\}$ is unique up to the action of the Weyl
group of~$R'$;
hence, Condition~\iref{i.dense.impr} in the
definition does not depend on the choice of $\{r_i\}$.
Furthermore, in view of \autoref{lem.roots},
both conditions are ``local'' and can be checked independently on irreducible
components: if $R=R'\oplus R''$, then
the set $\dense_n(R)$ is
the union over $n'+n''=n$ of the subsets
\[
\bigl\{\T'+\T''\bigm|(\T',\T'')\in\dense_{n'}(R')\times\dense_{n''}(R''),\
  \T'+\T''\le\bT\bigr\}.
\label{eq.pair}
\]

It follows from~\eqref{eq.basis} that, if $\bbT\subset N$ is embedded so that
$\bbT^\perp$ is root free, then the sublattice $\bT\subset\bbT$, regarded as
a $3$-form $\bT\to\rt(N)\dual$, is $2$-dense. An immediate consequence
of this observation and~\eqref{eq.pair} is the
following simple lemma.

\lemma\label{lem.root-free}
Assume that any one of the following conditions holds\rom:
\roster
\item\label{i.rf.m-1}
$\dense_nR_I=\varnothing$ whenever $n=0,1$ and $\ls|I|=m-1$\rom;
\item\label{i.rf.m-2}
$\dense_0R_I=\varnothing$ whenever $\ls|I|=m-2$\rom;
\item\label{i.rf.*}
$\dense_*R_I=\varnothing$ for a subset $I\subset\IS$\rom;
\item\label{i.rf.red}
there is an index $k\in\IS$ such that, for each
$\T\in\dense_n(R_{\IS\sminus k})$, $n=0,1,2$, and each isometry
$C:=\red_n\T\to R_k\dual$, one has $\rank\rt(C^\perp\cap R_k)>4-2n$.
\endroster
Then, a root-free lattice in the genus $\genus$ does not admit an embedding
to~$N$ with the orthogonal complement isomorphic to~$\bbT$.
\done
\endlemma

\subsection{The computation}\label{s.computation}
The sets $\dense_*R_k$ can be computed by \GAP~\cite{GAP4}, using the following
straightforward algorithm.
\roster
\item\label{i.step.A}
Consider the sets $\fA_s:=\{a\in R_k\dual\,|\,a^2\le s\}$, $s=4,8$.
\item\label{i.step.a8}
For~$\ba_8$, pick a representative of each $\OG(R_k)$-orbit of $\fA_8$.
\item\label{i.step.c8}
For~$\bc_8$, pick a representative of each $\stab(\ba_8)$-orbit of $\fA_8$.
\item\label{i.step.a4}
For~$\ba_4$, pick a representative of each $\stab(\ba_8,\bc_8)$-orbit of $\fA_4$.
\item\label{i.step.roots}
Compute $\rt(\T^\perp\cap R_k)$ and check if $\T\to R_k\dual$ is dense.
\endroster
At step~\iref{i.step.A}, we disregard square~$8$ vectors violating
Condition~\iref{i.dense.8} in \autoref{def.dense},
and at steps~\iref{i.step.c8} and~\iref{i.step.a4}, we check that the form
obtained
is bounded. A similar algorithm (with appropriate
verification at each step) can be used to enumerate all isometries
$\T\to R_k\dual$ of a particular $\Q$-valued form~$\T$.

If $\rank R_k>9$, the built-in group action algorithms are slow and we
replace them with
the run-length encoding as explained in \autoref{s.root}.

In
this section, we do not make use of any information about the group
$N/\!\rt(N)$ defining the extension $N\supset\rt(N)$ and merely compute
the sets
$\dense_*R_I$, $I\subset\IS$, inductively, using~\eqref{eq.pair}.
Below, we outline a few details.

\subsubsection{The lattices $\bD_{16}^+\oplus\bE_8$ and $\bE_8^3$}\label{ss.E8}
(Here, $\bD_{16}^+$, also denoted
$\bE_{16}$, is the so-called \emph{Barnes--Wall lattice}; it can be defined
as the only, up to isomorphism, even index~$2$ extension of~$\bD_{16}$.)
We have $\dense_n(\bE_8)=\varnothing$ unless $n=2$,
and \autoref{lem.root-free}\iref{i.rf.m-2} eliminates the
lattice $\bE_8^3$. Furthermore, if $\T\in\dense_2(\bE_8)$, then
$\red_2\T=0$ or $\Z a$, $a^2=4$.
On the other hand, for any vector $a\in\bD_{16}^+$,
$a^2=4$, one
has $\rt(a^\perp)\ne0$ (\cf. \autoref{s.root}), and we apply
\autoref{lem.root-free}\iref{i.rf.red}.

\subsubsection{The root systems $\bA_{17}\oplus\bE_7$ and $\bD_{10}\oplus\bE_7^2$}\label{ss.E7}
We have $\dense_0(\bE_7)=\varnothing$.
Then, $\dense_*(\bE_7^2)=\varnothing$ and
\autoref{lem.root-free}\iref{i.rf.*} eliminates $\bD_{10}\oplus\bE_7^2$.
There are several dozens of forms
$C:=\red_n\T$, $\T\in\dense_n(\bE_7)$, $n=1,2$.
Listing the isometries $C\to\bA_{17}\dual$, we conclude that
$\rank\rt(C^\perp\cap\bA_{17})\ge10$ and apply
\autoref{lem.root-free}\iref{i.rf.red}.

\subsubsection{The root systems $\bA_{11}\oplus\bD_7\oplus\bE_6$ and $\bE_6^4$}\label{ss.E6}
All sets $\dense_*(\bE_6)$, $\dense_*(\bD_7)$ are
nonempty; however, $\dense_0(\bE_6^2)=\varnothing$
and \autoref{lem.root-free}\iref{i.rf.m-2} eliminates the root system $\bE_6^4$.
Besides, $\dense_n(\bD_7\oplus\bE_6)=\varnothing$ unless $n=2$ and the forms
$C:=\red_2\T$, $\T\in\dense_2(\bD_7\oplus\bE_6)$, either
contain a vector of square~$\frac14$ or are $\Z a$, $a^2=\frac53$.
In the former case, $C$ admits no isometry to $\bA_{11}\dual$,
whereas in the latter
case, one has $\rt(a^\perp\cap\bA_{11})=\bA_9\oplus\bA_1$ and
\autoref{lem.root-free}\iref{i.rf.*} applies.
Both assertions on vectors in $\bA_{11}\dual$ follow from~\eqref{eq.shortest.A}.

\subsubsection{The root systems $\bD_{12}^2$ and $\bA_{12}^2$}\label{ss.D12}
The sets $\dense_*(\bD_{12})$ and $\dense_n(\bA_{12})$ for $n=0,1$ are empty,
and the lattices are eliminated by \autoref{lem.root-free}\iref{i.rf.m-1}.

\subsubsection{The root system $\bA_{15}\oplus\bD_9$}\label{ss.A15}
We have $\dense_n(\bD_9)=\varnothing$ unless $n=2$. Then, for each
class $C:=\red_2\T$, $\T\in\dense_2(\bD_9)$ (there are but two dozens of forms
not representing~$1$), we
enumerate the isometries
$C\to\bA_{15}\dual$, obtaining $\rank\rt(C^\perp)\ge10$.
Hence, \autoref{lem.root-free}\iref{i.rf.red} eliminates this lattice.

\subsubsection{Other root systems with $m\le6$ components}\label{ss.n<=6}
For the lattices listed below, all sets $\dense_*(R_k)$ are nonempty, and
we use~\eqref{eq.pair} to compute $\dense_*(R_I)$, $I\subset\IS$.
Then, we eliminate the lattice by applying
\roster*
\item
\autoref{lem.root-free}\iref{i.rf.m-1}, if
$\rt(N)=\bD_8^3$, $\bA_9^2\oplus\bD_6$, $\bA_8^3$, $\bA_7^2\oplus\bD_5^2$,
or $\bA_6^4$, or
\item
\autoref{lem.root-free}\iref{i.rf.m-2}, if
$\rt(N)=\bD_6^4$, $\bA_5^4\oplus\bD_4$, $\bD_4^6$, or $\bA_4^6$.
\endroster

\subsection{The lattices $\bD\sb{24}\sp+$ and $\bA\sb{24}\sp+$}\label{ss.D24}
By \cite[Chaprer 16]{Conway.Sloane},
the Niemeier lattice~$N$ is the extension of $\rt(N)=\bD_{24}$ or $\bA_{24}$
by the vector $\Ga_1\in\bD_{24}\dual$ or $\Ga_5\in\bA_{24}\dual$,
respectively, see~\eqref{eq.shortest.D} and~\eqref{eq.shortest.A}.
Thus, we apply the algorithm at the
beginning of \autoref{s.computation}, restricted to the subsets $\fA_s\cap N$,
and list all embeddings $\bT\into N$, arriving at
$\rank\rt(\bT^\perp)\ge15$ or $11$ for $N=\bD_{24}^+$ or $\bA_{24}^+$,
respectively. It follows that $\bT$ cannot be $2$-dense.

\remark\label{rem.singular}
In fact, it is not difficult to list all primitive embeddings
$\bbT\into N$, obtaining more than a thousand of lattices $\bbT^\perp$ in
the genus $\genus$ that contain roots.
Thus, already these two Niemeier lattices
give rise to a large number of singular spatial
models of the Fermat quartic.
\endremark

\subsection{Root systems with many components}\label{s.many.components}
Till the end
of this section, we consider a root system
$R:=\rt(N)=\bigoplus_{k=1}^mR_k$, where $m=8,12,24$ and
each $R_k$ is a copy of the
same irreducible root system $\bA_{n}$, $n:=24/m$.
Let
$\Ga_{i,k}\in R_k\dual$, $k\in\IS$, $0\le i\le n$, be
a distinguished shortest representative, given by~\eqref{eq.shortest.A}, of
the $i$-th element in the cyclic group $\discr R_k\cong\Z/(n+1)$;
sometimes, we will use the shortcut $\Ga_{i,I}:=\sum_{k\in I}\Ga_{i,k}$
for a subset $I\subset\IS$.

In all three cases, the kernel $\CK:=N/R\subset\discr R$
of the extension is generated by the $(m-1)$ elements of the form
\[*
{\textstyle\sum_{k=1}^{m}\Ga_{p_k,k}}\bmod R\in\discr R,
\]
where the sequences $(p_k)$ are obtained from the one given below
for each lattice by all
cyclic permutations of the subset $\{2,\ldots,m\}\subset\IS$
(see~\cite[Chaprer 16]{Conway.Sloane}).


An element $\Gb\in\discr R$ has the form
$\Gb_1+\ldots+\Gb_m$, where $\Gb_k\in\discr R_k$ is the projection, $k\in\IS$.
Similarly, an element $b\in R\dual$ has the form $b_1+\ldots+b_m$, where
$b_k:=\pr_kb\in R_k\dual$. For such an element, we define the \emph{support}
\[*
\supp\Gb:=\bigl\{k\in\IS\bigm|\Gb_k\ne0\bmod R_k\bigl\},\quad
\supp b:=\bigl\{k\in\IS\bigm|b_k\ne0\bigl\}
\]
and \emph{Hamming norm} $\norm|\,\cdot\,|:=\ls|\supp(\,\cdot\,)|$.
Clearly, one has $\supp b\supset\supp(b\bmod R)$.
The following statement is a consequence of \autoref{lem.roots}.

\lemma\label{lem.support}
If $b\in R\dual$ and $r\in R$ is a root, then $b\cdot r=0$
unless $r\in R_k$ for some index $k\in\supp b$.
In particular, for any isometry
$\bbT\into N$,
the two sets $\supp\bba_2^i$, $i=1,2$, are singletons contained in
$\supp\bbc_4$.
\done
\endlemma

We use a version of run-length encoding:
an element $b=b_1+\ldots+b_m\in R\dual$ is said to be
\emph{of the form}
\[*
\rle(b):=(s_1)^{u_1}\ldots(s_t)^{u_t},\quad 0<s_1<\ldots<s_t,\quad u_i>0,
\]
if, among the
projections $0\ne b_k\in R_k\dual$, there are exactly
$u_i$ vectors of square~$s_i$ for each $i=1,\ldots,t$, and there are no other
nonzero projections.

Unlike the previous two sections, below we consider an embedding
$\bbT\into N$ of the original lattice~$\bbT$ of rank~$5$. Most computations
are done up to the group $\OG(R)$; in fact, we study isometries
$\bbT\into R\dual$
using some limited information about the sublattice $N\subset R\dual$ which
must contain the image.
At the end, when classifying the root-free
lattices found, we switch to the finer group
\[
\OG(N)=\bigl\{g\in\OG(R)\bigm|g(\CK)=\CK\bigr\}
\label{eq.G(N)}
\]
of autoisometries of~$N$.
The \emph{combinatorial type} of an isometry $\bbT\into R\dual$ is defined as its
$\OG(R)$-orbit. (Recall that a basis for~$\bbT$ is assumed fixed; hence, we
can merely speak about $\OG(R)$-orbits of
ordered quintuples of vectors in $R\dual$.)

\subsection{The root system $\bA\sb3\sp{8}$}\label{s.8A3}
The kernel~$\CK$ is described
as in \autoref{s.many.components} by
\[*
(p_k)=(3,2,0,0,1,0,1,1).
\]
Since $\discr\bA_3\cong\Z/4$, we can refine the Hamming norm of
$\Gb\in\discr R$ to the \emph{type} $\type\Gb:=(\norm|\Gb|,\norm|2\Gb|)$.
We have
\[
\type\Gb\in\{(0,0), (4,0), (5,4), (7,4), (8,0), (8,8)\}\quad
 \text{for each $\Gb\in\CK$}.
\label{eq.A3.type}
\]
From this and~\eqref{eq.shortest.A},
one can see that the square~$4$
vectors $b\in N$ are of the form
\[
\bigl(\tfrac34\bigr)^4(1)^1,\ (1)^4,\ (2)^2,\ (4)^1.
\label{eq.roots.8A3}
\]
Another simple observation is the fact
that, for a vector $b\in\bA_3\dual$,
$b^2\le8$, one has $\rt(b^\perp)\ne0$ unless $b^2=5$.
Using~\eqref{eq.A3.type} again, we conclude that,
\[
\text{if $b\in N$, $b^2=8$, $\rle(b)\not\ni(5)$,
 then $\rt(b^\perp\cap R_k)\ne\varnothing$ for all $k\in\IS$}.
\label{eq.rle.5.1}
\]
In the exceptional cases
$\rle(b)=\bigl(\tfrac34\bigr)^4(5)^1$ or $(1)^3(5)^1$,
there is exactly one trivial intersection $\rt(b^\perp\cap R_k)$.

Consider a sublattice $\bbT\subset N$. By~\eqref{eq.rle.5.1}, a
necessary condition for $\rt(\bbT^\perp)=0$ is the bound
$\ls|\supp\bba_4\cup\supp\bbc_4|\ge7$. Since $\bba_4\cdot\bbc_4=2$,
using~\eqref{eq.roots.8A3}, we find three
combinatorial types of
pairs, with
$\bba_4=\Ga_{2,1}+\Ga_{1,\{2,\ldots,5\}}$ and $\bbc_4$ one of
\[*
\Ga_{2,1}+\Ga_{3,2}+\Ga_{1,\{3,6,7\}},\quad
\Ga_{2,7}+\Ga_{1,\{1,2,3,6\}},\quad
\Ga_{2,2}+\Ga_{1,\{3,4,6,7\}}.
\]
In each case,
the total support has length~$7$ and
$\rt\bigl((\Z\bba_4+\Z\bbc_4)^\perp\cap R_k\bigr)\ne0$ for each
index $k\in\IS$.
Then,
by~\eqref{eq.rle.5.1} again, $\bba_8$ must have a component of
length~$5$ in $R_8\dual$, and then it has
at most four other nonzero components. At most two components
contain~$\bba_2^1$ and~$\bba_2^2$; hence, there is at least one index
$k\in\IS$ left for which
\[*
\rt(\bbT^\perp\cap R_k)=
\rt\bigl((\Z\bba_4+\Z\bbc_4)^\perp\cap R_k\bigr)\ne0.
\]

\subsection{The root system $\bA\sb2\sp{12}$}\label{s.12A2}
The kernel~$\CK$ is described
as in \autoref{s.many.components} by
the sequence
\[*
(p_k)=(2,1,1,2,1,1,1,2,2,2,1,2);
\]
it
is the ternary Golay code~$\Golay_{12}$.
We have
\[
\norm|\Gb|\in\{0,6,9,12\}\quad\text{for each $\Gb\in\CK$}.
\label{eq.A2.type}
\]
In view of~\eqref{eq.shortest.A}, it follows that
$\rle(b)=\bigl(\tfrac23\bigr)^6$ or $(2)^2$ if $b\in N$ and $b^2=4$,
whereas the square~$8$ vectors $b\in N$ are of the form
\[*
\bigl(\tfrac23\bigr)^{12},\
\bigl(\tfrac23\bigr)^9(2),\
\bigl(\tfrac23\bigr)^8\bigl(\tfrac83\bigr),\
\bigl(\tfrac23\bigr)^6(2)^2,\
\bigl(\tfrac23\bigr)^5(2)\bigl(\tfrac83\bigr),\
\bigl(\tfrac23\bigr)^5\bigl(\tfrac{14}3\bigr),\
\bigl(\tfrac23\bigr)^4\bigl(\tfrac83\bigr)^2,\
(2)^4,\
(2)(6).
\]
If $b\in\bA_2\dual$, $b^2<8$, then $\rt(b^\perp\cap\bA_2)\ne0$
unless $b^2=2$ or $\frac{14}3$.
Hence, a
necessary condition for $\rt(\bbT^\perp)=0$ is that
$\ls|\supp\bba_4\cup\supp\bbc_4|\ge8$. Since
$\bba_4\cdot\bbc_4=2$,
up to the action of $\OG(\bA_2^{12})$ we have
\[*
\bba_4=\Ga_{1,\{1,\ldots,6\}}\quad\text{and}\quad
\bbc_4=\Ga_{1,\{4,\ldots,9\}}\quad\text{or}\quad
\Ga_{2,\{1,2\}}+\Ga_{1,\{5,\ldots,8\}}.
\]
In order to eliminate the roots in $R_{10}$ through~$R_{12}$, the remaining
vector $\bba_8$ must be of the form $(2)^4$; then, at least four pairs of
roots survive to $\bbT^\perp$.

\subsection{The root system $\bA\sb1\sp{24}$}\label{s.24A1}
This is the only Niemeier lattice containing root-free sublattices in the
genus $\genus$.
The kernel~$\CK$ is described
as in \autoref{s.many.components} by
\[*
(p_k)=(1,0,0,0,0,0,1,0,1,0,0,1,1,0,0,1,1,0,1,0,1,1,1,1);
\]
it is the binary Golay code~$\Golay_{24}$.
We have
\[
\norm|\Gb|\in\{0,8,12,16,24\}\quad\text{for each $\Gb\in\CK$}.
\label{eq.A1.type}
\]
The vectors of Hamming norm~$8$ are called
\emph{octads}; their supports are complementary to those of norm~$16$.
The vectors $b\in N$ of interest are of the form
\[*
\bigl(\tfrac12\bigr)^8,\
 (2)^2\
 \text{if $b^2=4$};\quad
\bigl(\tfrac12\bigr)^{16},\
 \bigl(\tfrac12\bigr)^{12}(2),\
 \bigl(\tfrac12\bigr)^8(2)^2,\
 (2)^4\
 \text{if $b^2=8$}.
\]
If $\bbT\subset N$, then
\roster*
\item
$\bbc_4$ is of the form $\bigl(\tfrac12\bigr)^8$ and
$\supp\bbc_4\supset\supp\bba_2^i$, $i=1,2$,
\cf. \autoref{lem.support}, and
\item
$\rt(\bbT^\perp)=0$ if and only if
$\supp\bba_4\cup\supp\bbc_4\cup\supp\bba_8=\IS$.
\endroster
Taking into account the fact that $\bba_4\cdot\bbc_4=2$, we arrive at
three combinatorial types, which can be encoded by the following diagrams:
\begin{gather}
\diagram[X_{48}]{
\n\n\.\.\\
\-\-\-\-\-\-\-\-\\
\ \ \ \ \ \ \ \ \-\-\-\-\-\-\-\-\-\-\-\-\-\-\-\-
}\punct,\label{X48}\\\noalign{\allowbreak}
\diagram[X_{56}]{
\n\n\ \ \-\-\-\-\-\-\-\-\\
\-\-\-\-\-\-\-\-\\
\ \ \ \ \ \ \ \ \=\=\-\-\-\-\-\-\-\-\-\-\-\-\-\-
}\punct,\label{X56}\\\noalign{\allowbreak}
\diagram[X_{56}\bis]{
\n\n\ \ \-\-\-\-\-\-\-\-\\
\-\-\-\-\-\-\-\-\\
\ \ \ \ \=\=\-\-\ \ \ \ \-\-\-\-\-\-\-\-\-\-\-\-
}\punct.\label{X56-2}
\end{gather}
In the diagrams, we list, line by line, the images of the basis vectors
of~$\bbT$ other than roots, indicating their nonzero coordinates in a
standard basis for $\bA_1^{24}$ with the following symbols
(some of which are used later in the paper):
\[
\symb-\tfrac12,\quad
\symb=-\tfrac12,\quad
\symb+\tfrac32,\quad
\symb\circ1.
\label{eq.diagram}
\]
The images of the basis elements of square~$2$ are shown
by~$\bullet$, which can appear in any line; we assume that the only
nonvanishing coordinate of each root equals~$1$.
Similar diagrams will be used throughout the paper.

\subsection{The classification}\label{s.classification}
The principal result of this section can be summarized in the following
statement.

\theorem\label{th.2xS}
There are two isomorphism classes of extensions $N\supset S$, where $N$ is a
Niemeier lattice, $S$ is a root-free lattice in the genus $\genus$, and
$S^\perp\cong\bbT$. In both cases, $N=N(\bA_1^{24})$\rom;
the isometry $\bbT=S^\perp\into N(\bA_1^{24})$ has
combinatorial type as in diagram~\eqref{X48} or~\eqref{X56} in
\autoref{s.24A1}.
\endtheorem

\proof
The Leech lattice has no roots and, hence, cannot contain $\bbT$ as a
sublattice.
All other Niemeier lattices except $N(\bA_1^{24})$ have been eliminated in
\autoref{s.computation}--\autoref{s.12A2},
and the embeddings $\bbT\into N:=N(\bA_1^{24})$ have been reduced to
three combinatorial types
in \autoref{s.24A1}. The type~\eqref{X56-2} differs from~\eqref{X56}
by the basis change
$\bba_4\mapsto-\bba_4$,
$\bbc_4\mapsto\bbc_4-\bba_4$; hence, it results in an isomorphic pair $(N,S)$.

It is straightforward that any isometry $\bbT\into R\dual$ as in the statement
is $\OG(R)$-equivalent to an isometry $\bbT\into N$, and there only remains
to show that the latter is unique up to the action of~$\OG(N)$.
Recall that the group of
automorphisms of the Golay code~$\Golay_{24}$
is the Mathieu group $M_{24}\subset\SG{24}$,
and the action of $M_{24}$ on~$\IS$ has the following properties
(see \cite[Chapter 10]{Conway.Sloane}):
\roster
\item\label{i.M.transitive}
the action is transitive on the $759$ octads in $\Golay_{24}$;
\item\label{i.M.A8}
the stabilizer of an octad $b\in\Golay_{24}$ factors to
$\AG8\subset\Sym(\supp b)\cong\SG8$.
\endroster
By~\iref{i.M.transitive}, we can fix the octad~$\bbc_4$.
Then, the uniqueness of all further choices in the case~\eqref{X48} follows
immediately from Statement~\iref{i.M.A8}. For the other case~\eqref{X56},
we need an
additional observation, which is easily confirmed by \GAP~\cite{GAP4}:
\roster[\lastitem]
\item\label{i.M.pairs}
the action of $M_{24}$ is transitive on the set of ordered pairs
$(b_1,b_2)\in\Golay_{24}\times\Golay_{24}$ of octads such that
the set $C:=\supp b_1\cap\supp b_2$ is of size~$4$;
\item\label{i.M.A4xA4}
the stabilizer of an ordered pair $(b_1,b_2)\in\Golay_{24}\times\Golay_{24}$ as
above factors to an index~$2$ subgroup of $\Sym(C_1)\times\Sym(C_2)$,
where $C_i:=\supp b_i\sminus C$.
\endroster
Thus, by~\iref{i.M.pairs}, we can fix the pair $(\bba_4,\bbc_4)$; then,
by~\iref{i.M.A4xA4}, there is a unique choice for the two singletons
$\supp\bba_2^{1,2}\subset\supp\bbc_4\sminus C$ and the subset
$V_-\subset\supp\bba_4\sminus C$ where the coordinates of~$\bba_8$ are equal
to $-\frac12$ (the two ``$=$'' in the diagram).
\endproof

\corollary[of the proof]\label{cor.aut.T}
Let $\bbT\into N$ be an isometry as in \autoref{th.2xS}. Then, an
autoisometry $g\in\OG(\bbT)$ extends to
to an autoisometry of~$N$ if and only if $g$ preserves the
combinatorial type of $\bbT\into N$.
\done
\endcorollary

\section{Proof of \autoref{th.main}}\label{proof.main}

In this section, we complete the proof of \autoref{th.main} by
analyzing
extensions $H_2(\Fermat)\supset\NS(\Fermat)\supset{\bS_n}$ of the
two lattices $\bS_n$, $n=48,56$, constructed in \autoref{s.24A1}.

\subsection{Automorphisms of $\bbT$}\label{s.aut.T}
The discriminant $\CS:=\discr\bbT\cong(\Z/4)\oplus(\Z/8)^2$ is generated by three pairwise orthogonal
elements $\Ga,\Gb_1,\Gb_2$ of orders $4,8,8$ and squares
$-\frac14,-\frac18,-\frac18\bmod2\Z$, respectively.
In the notation of \autoref{s.K3}, we have
\[*
\CS_h=\CS_h^+
 =\bigl\{\pm\Ga+c_1\Gb_1+c_2\Gb_2\bigm|c_1,c_2=0,4\bigr\}.
\]
The group $A:=\Aut\CS$ has order~$128$, and its image~$A_h$ in $\Sym(\CS_h)$
is isomorphic to $(\Z/2)\times\DG8$. In an appropriate ordering of~$\CS_h$,
this image is generated by
\[*
a:=(1,2)(3,4)(5,6)(7,8),\quad
u:=(2,3)(6,7),\quad
-\id=(1,5)(2,6)(3,7)(4,8).
\]
In particular, $A_h$ is transitive on~$\CS_h$, \cf. \autoref{rem.S+}.

The following few statements are straightforward.
\roster
\item\label{i.O(\T)}
The group $\OG(\bbT)$ has order $64$.
\item\label{i.G}
The image $G\subset A$ of $\OG(\bbT)$ is isomorphic to $(\Z/2)^4$.
\item\label{i.Gh}
The image $G_h\subset A_h$ of $G$ is the order~$8$ subgroup generated by
$a$, $u\1au$, and~$-\id$; it acts on~$\CS_h$ simply transitively.
\endroster

Below, we fix an isometry $\iota\:\bbT\into N:=N(\bA_1^{24})$ as in
\autoref{th.2xS}, consider the lattice $S:=-\bbT^\perp$, and identify
$\discr S$ and~$\CS$ \via\ the isometry~$\psi$ corresponding to the
extension $N\supset(-S)\oplus\bbT$, see \autoref{cor.extension}.
By means of this identification, we can speak about the images
$H\subset A$ and $H_h\subset A_h$ of the
group $\OG(S)$. Consider also the subgroup
$\OG^\iota(\bbT)\subset\OG(\bbT)$
consisting of the automorphisms preserving the
combinatorial type of~$\iota$, see \autoref{s.many.components}, and its images
$G^\iota\subset G$ and $G_h^\iota\subset G_h$. We have
\roster[\lastitem]
\item\label{i.G.cap.H}
$G^\iota=G\cap H$ and, hence, $G^\iota_h\subset G_h\cap H_h$,
see Corollaries~\ref{cor.extension}\iref{i.ext.aut} and~\ref{cor.aut.T}.
\endroster
On the other hand, the two lattices~$S$ corresponding to the two isometries
given by \autoref{th.2xS} are not isomorphic, see \autoref{lem.S56} below; hence,
each extension $N\supset(-S)\oplus\bbT$ is unique up to isomorphism and,
by \autoref{cor.extension}\iref{i.ext.iso}, the group~$A$ is a single double
coset, \ie, $A=\{gh\,|\,g\in G,\ h\in H\}$.
Then, by~\iref{i.G.cap.H} and~\iref{i.Gh} above,
\begin{gather}
[A:H]\le[G:G^\iota],
 \label{eq.index}\\
H_h\supset G^\iota_h\cdot\<\bar g\>\quad
 \text{for some $\bar g\in A_h\sminus G_h$}.
 \label{eq.H.h}
\end{gather}

For each $\Gg\in\CS_h$,
we have $\CK_\Gg^\perp/\CK_\Gg=\Gg^\perp\cong-\discr\TL$,
see \autoref{rem.S+}, and
the restriction establishes an isomorphism
\[
\stab\Gg=\Aut\Gg^\perp\cong(\Z/2)\times\DG8.
\label{eq.stab.gamma}
\]
The intersection $G\cap\stab\Gg$ is a subgroup of order~$2$.

\subsection{The lattice $\bS\sb{48}$}\label{s.S48}
Let $\iota\:\bbT\into N$ be the isometry with the combinatorial type as
in~\eqref{X48}
in \autoref{s.24A1}, and denote $\bS_{48}:=-\bbT^\perp$.
The following lemma is proved by a straightforward computation.

\lemma\label{lem.S48}
The form $\bS_{48}\dual$ contains a unique pair $\pm a$ of vectors of
square~$(-1)$\rom; they represent the element
$2\Ga+4(\Gb_1+\Gb_2)\in\discr\bS_{48}$. Furthermore, each element
$\Gg\in\CS_h$ is represented by exactly $48$ vectors $a\in\bS_{48}\dual$,
$a^2=-\frac94$.
\done
\endlemma

\lemma\label{lem.S48.aut}
The canonical homomorphism $\OG(\bS_{48})\to\Aut\CS$ is surjective.
Hence,
the action of $\OG(\bS_{48})$ is transitive on~$\CS_h$ and, for each
$\Gg\in\CS_h$, the stabilizer $\Stab\Gg$ projects to the
full automorphism group $\Aut\Gg^\perp$.
\endlemma

\proof
Any autoisometry of~$\bbT$ preserves the combinatorial type of~$\iota$.
Hence, we have $G^\iota=G$ and, by~\eqref{eq.index}, $H=A$
(\cf. also \autoref{lem.index}).
The other statements follow from the properties of the group $A=\Aut\CS$
discussed in \autoref{s.aut.T}.
\endproof

\subsection{The lattice $\bS\sb{56}$}\label{s.S56}
Let $\iota\:\bbT\into N$ be the isometry with the combinatorial type as
in~\eqref{X56}
in \autoref{s.24A1}, and denote $\bS_{56}:=-\bbT^\perp$.
The following lemma is proved by a straightforward computation.

\lemma\label{lem.S56}
The form $\bS_{56}\dual$ has no vectors of square~$(-1)$.
Each element
$\Gg\in\CS_h$ is represented by exactly $56$ vectors $a\in\bS_{56}\dual$,
$a^2=-\frac94$.
As a consequence, we have $\bS_{56}\not\cong\bS_{48}$, \cf.
\autoref{lem.S48}.
\done
\endlemma

\lemma\label{lem.S56.aut}
The image $H_{56}$ of the canonical homomorphism $\OG(\bS_{56})\to\Aut\CS$ is a
subgroup of index~$2$.
This subgroup has the following properties\rom:
\roster
\item\label{i.H56.Sh}
the subgroup~$H_{56}$ is transitive on~$\CS_h$\rom;
\item\label{i.H56.8}
the subgroup~$H_{56}$ is transitive on the set
$\CS_{56}:=\{\Gg\in\CS\,|\,\Gg^2=-\frac18\bmod2\Z\}$\rom;
\item\label{i.H56.refl}
the subgroup $H_{56}$ contains the index~$4$ subgroup~$H'$ generated by reflections
defined by the elements of
$\CS_r:=\{\Ga\in\CS\,|\,\text{$\Ga^2=\frac38$ or $\frac34\bmod2\Z$,
$4\Ga\ne0$}\}$.
\endroster
\endlemma

\proof
We have $[G:G^\iota]=2$, as the combinatorial type depends on the choice of a
basis, see \eqref{X56-2} \vs.~\eqref{X56} in \autoref{s.24A1}.
Hence, $[A:H]=2$ by~\eqref{eq.index}.
(Alternatively, the same conclusion follows from \autoref{lem.index}.)
Precise computation shows that
$G_h^\iota\subset G_h$ is the index~$2$ subgroup generated by~$a$
and~$-\id$ (see \autoref{s.aut.T} for the notation),
and one can check that $G_h^\iota\cdot\<g\>=A_h$
for \emph{each} element $\bar g\in A_h\sminus G_h$. Hence, \eqref{eq.H.h}
gives us
$H_h=A_h$, and
Statement~\iref{i.H56.Sh}
follows from Statement~\iref{i.Gh}
in \autoref{s.aut.T}.

Fix $\Gg\in\CS_h$. The group $\Stab\Gg$ has been studied in
\cite[Lemma~6.19]{DIS} as $\OG_h(\bX_{56})$; its image in $\Aut\Gg^\perp$
is the index~$2$ subgroup generated by reflections
defined by the elements $\Ga\in\Gg^\perp\cap\CS_r$.
Any element of~$\CS_r$
is orthogonal to some $\Gg\in\CS_h$; hence,
$H'\subset H_{56}$. There remains to observe that $H'$ acts transitively on
$\CS_{56}$.
\endproof


\subsection{End of the proof}\label{s.end.proof}
Since we are interested in a smooth model, the lattice $S\in{\genus}$ in
\autoref{th.K3} must be root free, \ie, one of the two lattices
$\bS_n$, $n=48$ or~$56$, introduced above.
By Lemmas~\ref{lem.S48.aut} and~\ref{lem.S56.aut},
the action of $\OG(\bS_n)$ on $\CS_h=\CS_h^+$ is transitive and,
for each lattice, it suffices to
consider one representative $\Gg\in\CS_h$.
The models obtained are birational (and then smooth) due to
\autoref{th.K3}\iref{i.model.birational}, \iref{i.model.nonsingular}
and Lemmas~\ref{lem.S48}
and~\ref{lem.S56}; by the same lemmas
and \autoref{th.K3}\iref{i.model.lines}, the quartic obtained contains $n$
lines.

To complete the proof, we need to analyze the set
\[*
\OG^+(\TL)\backslash\Aut\discr\TL/\Stab\Gg=\OG^+(\TL)\backslash\Aut\Gg^\perp/\Stab\Gg.
\]
If $n=48$, \autoref{lem.S48.aut} gives us a unique double coset,
hence one model.
If $n=56$, \autoref{lem.S56.aut} implies that the image
of $\Stab\Gg$ in $\Aut\Gg^\perp$ is the index~$2$ subgroup
$H_{56}\cap\stab\Gg=H'\cap\stab\Gg$ which contains the image
of $\OG^+(\TL)$.
Hence, there are two double cosets resulting in two
quartics.
One can easily check (or merely refer to~\cite{DIS}) that the two double
cosets are interchanged by the full group $\OG(\TL)$; hence, the two quartics
are complex conjugate (see the beginning of \autoref{s.K3}).

Finally, the three projective quartics obtained are identified with
the classical Fermat
quartic $X_{48}$ or the pair $X_{56}$, $\bar X_{56}$ constructed in
\cite{DIS} and studied further in \cite{Shimada:X56}
according to the number of lines contained in the
surface.
\qed

\section{Proof of \autoref{th.others}}\label{S.others}

The approach used in this section is similar, but not identical to that of
\autoref{S.Niemeier}. For some mysterious reason, it does not work very well
for the Fermat quartic. Hence, below we refer to \autoref{th.main} and assume
that $T\ne[8,0,8]$.

\subsection{Niemeier lattices with many roots}\label{s.many.roots}
Throughout the proof, $\TL$ is a fixed positive definite even lattice of
rank~$2$, $\det\TL\le80$, and $\TL\ne[8,0,0]$ or $[4,0,16]$.

Any test lattice~$\bbT$ (see
\autoref{s.reduction}) can be decomposed into a sum (direct, but
not orthogonal) $\bbT=\rt(\bbT)+\bT$. We assume that the basis for $\bbT$ is
obtained from one for $\rt(\bbT)$ inductively, by adding at each step a
shortest vector possible, and then $\bT$ is spanned by the basis vectors
$\bbc_1,\ldots,\bbc_p$, $p=5-\rank\rt(\bbT)$, that are not in $\rt(\bbT)$.

\convention\label{conv.V}
We also assume that
\roster*
\item
$\bbc_i^2\le16$ for each $i=1,\ldots,p$, and
\item
$\rank\rt(\bbT)\ge2$ (this is the reason for excluding $\TL=[4,0,16]$);
\endroster
among the lattices satisfying these conditions, we choose $\bbT$ by
minimizing the rank $p=5-\rank\rt(\bbT)$, then the maximal square~$\bbc_p^2$,
and then the trace $\bbc_1^2+\ldots+\bbc_p^2$ (see
\url{http://www.fen.bilkent.edu.tr/~degt/papers/test_matrix.zip}).

With a pair $\bT\subset\bbT$ fixed, we call a $\Q$-valued $p$-form~$V$ with a
distinguished basis $\bbc_1,\ldots,\bbc_p$ \emph{bounded}, $V\le\bT$, if
$\bT-V$ is positive semi-definite.
\endconvention

Fix $\bT\subset\bbT$ as in \autoref{conv.V}. Any isometry
$\bbT\into N$ to a Niemeier lattice~$N$ restricts to an isometry
$\rt(\bbT)\into\rt(N)$, and we can consider the root lattice
\[*
R:=\rt\bigl(\rt(\bbT)^\perp\bigr)
 =\bigoplus R_k,\quad k\in\IS:=\{1,\ldots,m\},
\]
where $R_k$ are the irreducible components.
As in \autoref{s.proj}, we let $R_I:=\bigoplus_{k\in I}R_k$ for a subset $I\subset\IS$.
Note that this lattice~$R$ differs
from the maximal root lattice $\rt(N)$ considered in \autoref{s.proj},
and \autoref{lem.roots} takes
the following form.

\lemma\label{lem.roots.2}
If $r\in\rt(\bbT)^\perp$ is a root, then $r\in R_k$ for some index $k\in\IS$.
\done
\endlemma

\warning
The
lattices
$\bT\subset\bbT$ and $R$ introduced here differ from
$\bT=\rt(\bbT)^\perp$
in \autoref{s.genus} and $R=\rt(N)$ in \autoref{s.proj}
(and so are the irreducible components~$R_k$ and root lattices $R_I$).
This difference is due to a slight change in the approach: instead of considering
all isometries $\rt(\bbT)^\perp\into N$ and selecting those with at most two
roots in the orthogonal component, we \emph{start} with an embedding
$\rt(\bbT)\into\rt(N)$ and try to find dense (see below) bounded sublattices
in $\rt(\bbT)^\perp\subset N$.
\endwarning

At this stage, we do not insist that the isometry $\bbT\into N$ or its
restriction to $\rt(\bbT)$ should be primitive, \cf.
\autoref{rem.genus}.
The possible isometries $\rt(\bbT)\into\rt(N)$ can easily be classified
step by step, by embedding an irreducible component $Q\subset\rt(\bbT)$ to an
irreducible component $P\subset\rt(N)$ and replacing~$P$ with $\rt(Q^\perp)$.
Isometries of irreducible root systems are well known; for the reader's
convenience, we list them in \autoref{tab.roots}, where we stretch
the notation and let
$\bA_{n-1}=\bD_n:=0$ for $n\le1$ and
\[*
\bD_2:=\bA_1^2,\quad
\bD_3:=\bA_3,\quad
\bE_3:=\bA_1\oplus\bA_2,\quad
\bE_4:=\bA_4,\quad
\bE_5:=\bD_5.
\]
(These conventions are based on the structure of the discriminant group.)
In a few cases, an isometry $Q\into P$ is not unique; most notably, $\bA_3$
can be embedded to~$\bD_n$ as $\bA_3$ or $\bD_3$.
\table
\caption{Embeddings of irreducible root lattices}\label{tab.roots}
\def\+{\noalign{\kern2\lineskip}}
\def\-{\noalign{\kern-\lineskip}}
\def\balign{\vtop\bgroup\halign\bgroup
 \strut\quad$##$\hss\quad&&$##$\hss\quad\cr
 P&Q\into P&\rt(Q^\perp)\cr
\noalign{\vspace{1pt}\hrule\vspace{2pt}}}
\def\ealign{\crcr
 \egroup\egroup}
\hbox to\hsize{\hss\valign{\hrule\kern2pt#\vss\kern2pt\hrule\cr\balign
\bA_n&\bA_m,\ 1\le m\le n&\bA_{n-m-1}\cr\+
\bD_n&\bA_1&\bA_1\oplus\bD_{n-2}\cr\-
     &\bA_m,\ 2\le m<   n&\bD_{n-m-1}\cr\-
     &\bD_m,\ 3\le m\le n&\bD_{n-m}\cr\+
\bE_6&\bA_1&\bA_5\cr\-
     &\bA_2&\bA_2^2\cr\-
     &\bA_3&\bA_1^2\cr\-
     &\bA_4, \bA_5&\bA_1\cr\-
     &\bD_4, \bD_5, \bE_6&0\cr
\ealign\cr\noalign{\quad}
\balign
\bE_7&\bA_1&\bD_6\cr\-
     &\bA_2&\bA_5\cr\-
     &\bA_3&\bA_3\oplus\bA_1\cr\-
     &\bA_4, \bA_5&\bA_2\cr\-
     &\bA_5, \bD_5, \bD_6&\bA_1\cr\-
     &\bA_6, \bA_7, \bE_6, \bE_7&0\cr\-
     &\bD_4&\bA_1^3\cr\+
\bE_8&\bA_m,\ 1\le m\le 5&\bE_{8-m}\cr\-
     &\bA_6, \bA_7, \bE_7&\bA_1\cr\-
     &\bA_7, \bA_8, \bE_8&0\cr\-
     &\bD_m,\ 4\le m\le 8&\bD_{8-m}\cr\-
     &\bE_6&\bA_2\cr
\ealign\crcr}\hss}
\endtable

Fix an isometry $\rt(\bbT)\into\rt(N)$ and let $R$ be as above.
Assume that there is an extension
$\bbT\into N$ and consider the projections of $\bT$ to the groups
$R_I\dual$; as in \autoref{s.proj}, these projections are regarded
as isometries
\[*
V_k\to R_k\dual,\quad k\in\IS,\quad\text{and}\quad
V_I\to R_I\dual,\quad I\subset\IS,
\]
and referred to as $p$-forms.
Clearly, all forms $V_I$ are bounded.
A $p$-form $V\to R_I\dual$ is called
\emph{dense} if $\rt(V^\perp\cap R_I)=0$; we denote by $\dense(R_I)$ the
set of the Gram matrices of bounded dense $p$-forms $V\to R_I$.
If $I=I'\cup I''$ and $I'\cap I''=\varnothing$, then, by
\autoref{lem.roots.2}, we have
\[*
\dense(R_I)=
\bigl\{\T'+\T''\bigm|(\T',\T'')\in\dense(R_{I'})\times\dense(R_{I''}),\
  \T'+\T''\le\bT\bigr\},
\]
\cf.~\eqref{eq.pair}, \ie, the sets $\dense(R_I)$
for all subsets $I\subset\IS$ can be computed inductively
starting from $\dense(R_k)$, $k\in\IS$.
Note that this computation is reusable, as the set $\dense(R_I)$ depends
on~$R_I$ and $\bT$ only, and the sublattice $\bT$ chosen according to
\autoref{conv.V} is often shared by many transcendental lattices~$\TL$.

The following statement has been obtained using \GAP~\cite{GAP4} and the
algorithm outlined in \autoref{s.computation} to enumerate
the bounded isometries $V\to R_k$.

\lemma\label{lem.dense}
Assume that the test lattices $\bT\subset\bbT$ are chosen as in
\autoref{conv.V} and that the root system $\rt(N)$ has at most six
irreducible components. Then, for any isometry $\rt(\bbT)\into\rt(N)$, one
has $\dense(R)=\dense(R_\IS)=\varnothing$.
\endlemma

As an immediate consequence, we conclude that, for $N$ as in
\autoref{lem.dense} and any isometry $\bbT\into N$, the orthogonal complement
$\bbT^\perp$ is not root free.

\subsection{Root systems with many components}\label{s.many.2}
Fix $\bT\subset\bbT$ as above and consider one of the remaining three
Niemeier lattices (see \autoref{s.many.components}),
assuming that $\rt(\bbT)$ admits an isometry to
the root lattice $R:=\rt(N)$ and
using the definitions
and notation introduced in \autoref{s.many.components}--\autoref{s.24A1}.
This time, we start with building an isometry $\bT\into R\dual$, considering
the latter up to the action of $\OG(R)$ and using the list of types/Hamming
norms of the elements of~$\CK$ only, see~\eqref{eq.A3.type}, \eqref{eq.A2.type},
\eqref{eq.A1.type}; this restriction is applied to each element of the group
$\bT\bmod R$.
It is not difficult to enumerate pairs of vectors $\bbc_1,\bbc_2\in R\dual$;
if a third vector $\bbc_3$ is to be added, one can usually limit the choices
similar to \autoref{s.8A3}--\autoref{s.24A1}, by analyzing the position of
the roots in~$\bbT$ with respect to the supports $\supp(\bbc_i)$.

There remains to extend the isometries found to $\rt(\bbT)\into R$ and select
those for which the lattice $\bbT^\perp$ is root free.
We arrive at the following
fifteen
combinatorial types
of isometries
$\bbT\into R\dual$ for $R=\bA_1^{24}$:
\begin{gather}
\diagram[X_{64}]{
\-\-\-\-\-\-\-\-\\
\ \ \ \ \ \ \ \ \-\-\-\-\-\-\-\-\\
\n\n\ \ \ \ \ \ \ \ \ \ \ \ \ \ \-\-\-\-\-\-\-\-
}\punct,\label{X64}\\\noalign{\allowbreak}
\diagram[X_{60}'']{
\n\n\ \ \ \ \ \ \ \ \ \ \n\\
\-\-\-\-\-\-\-\-\-\-\-\-\\
\ \ \ \ \ \ \ \ \ \ \ \ \-\-\-\-\-\-\-\-\-\-\-\-
}\punct,\label{X60''}\\\noalign{\allowbreak}
\diagram[X_{60}']{
\-\-\-\-\-\-\-\-\\
\-\-\ \ \ \ \ \ \-\-\-\-\-\-\\
\n\n\ \ \ \ \ \ \=\-\-\-\-\-\-\-\-\-\-\-\-\-\-\-
}\punct,\label{X60'}\\\noalign{\allowbreak}
\diagram[Q_{56}]{
\-\-\-\-\ \ \ \ \-\-\-\-\\
\-\-\-\=\-\-\-\-\\
\n\n\ \ \ \ \ \ \-\-\-\-\-\-\-\-\-\-\-\-\-\-\-\-
}\punct,\label{Q56}\\\noalign{\allowbreak}
\diagram[\tilde Y_{48}']{
\ \ \ \ \n\ \ \ \ \ \ \ \ \ \ \ \ \ \ \ \ \ \n\n\\
\-\-\-\-\-\-\-\-\\
\-\-\-\-\ \ \ \ \.\.\-\-\-\-\-\-\-\-\-\-\-\-
}\punct,\label{Y48'}\\\noalign{\allowbreak}
\diagram[\quadric_{48}']{
\ \ \n\n\n\ \ \ \ \ \ \ \ \ \ \ \ \ \ \ \ \ \ \ \\
\.\.\\
\-\-\-\-\-\-\-\-\-\-\-\-\-\-\-\-\-\-\-\-\-\-\-\-
}\punct,\label{PP48}\\\noalign{\allowbreak}
\diagram[*]{
\ \ \ \ \ \ \ \ \n\n\n\ \ \ \ \ \ \ \ \ \ \ \ \ \\
\-\-\-\-\-\-\-\-\\
\.\.\ \ \ \ \ \ \-\-\-\-\-\-\-\-\-\-\-\-\-\-\-\-
}\punct,
\label{X52''}\\\noalign{\allowbreak}
\diagram[\theA\bis]{
\ \ \ \ \ \ \ \ \n\n\n\ \ \ \ \ \ \ \ \ \ \ \ \ \\
\-\-\-\-\-\-\-\-\\
\=\=\-\-\-\-\-\-\-\-\-\-\-\-\-\-\-\-\-\-\-\-\-\-
}\punct,
\label{X52''2}\\\noalign{\allowbreak}
\diagram[Y_{52}'']{
\ \n\n\ \ \ \ \ \ \n\ \ \ \ \ \ \ \ \ \ \ \ \ \ \\
\-\-\-\-\-\-\-\-\\
\.\ \ \ \ \ \ \ \+\-\-\-\-\-\-\-\-\-\-\-\-\-\-\-
}\punct,\label{Y52''}\\\noalign{\allowbreak}
\diagram[\tilde Y_{48}'']{
\ \ \ \n\n\ \ \ \n\ \ \ \ \ \ \ \ \ \ \ \ \ \ \ \\
\=\-\-\-\-\-\-\-\\
\.\.\.\ \ \ \ \ \-\-\-\-\-\-\-\-\-\-\-\-\-\-\-\-
}\punct,\label{Y48''}\\\noalign{\allowbreak}
\diagram[*]{
\-\-\-\-\-\-\-\-\ \ \ \ \ \ \ \ \ \ \ \ \ \ \n\n\\
\ \ \ \ \ \ \ \ \-\-\-\-\-\-\-\-\\
\-\-\-\-\ \ \ \ \-\-\-\-\ \ \ \ \-\-\-\-\-\-\-\-
}\punct,\label{Z48'}\\\noalign{\allowbreak}
\diagram[*]{
\-\-\-\-\-\-\-\-\\
\ \ \ \ \ \ \ \ \-\-\-\-\-\-\-\-\\
\.\.\ \ \n\n\ \ \ \ \ \ \ \ \ \ \-\-\-\-\-\-\-\-
}\punct,\label{PP48''}\\\noalign{\allowbreak}
\diagram[*]{
\-\-\-\-\-\-\-\-\\
\ \ \ \ \ \ \ \ \-\-\-\-\-\-\-\-\\
\-\-\-\-\n\n\ \ \=\=\-\-\ \ \ \ \-\-\-\-\-\-\-\-
}\punct,\label{Q48}\\\noalign{\allowbreak}
\diagram[Q_{52}'']{
\-\-\-\-\-\-\-\-\\
\=\=\-\-\ \ \ \ \-\-\-\-\\
\-\-\-\-\n\n\ \ \ \ \ \ \-\-\-\-\-\-\-\-\-\-\-\-
}\punct,\label{Q52''}\\\noalign{\allowbreak}
\diagram[X_{52}']{
\-\-\-\-\-\-\-\-\\
\=\-\ \ \ \ \ \ \-\-\-\-\-\-\\
\-\-\.\ \n\n\ \ \ \ \ \ \ \ \-\-\-\-\-\-\-\-\-\-
}\punct{}\label{X52'}
\end{gather}
\qlabel{Z48'}{Z_{48}'}{Z48'}%
\qlabel{Q48}{Q_{48}}{Q48}%
\qlabel{Q54}{Q_{54}}{X52''}%
\qlabel{X52''}{X_{52}''}{X52''}%
\qlabel{PP48''}{\quadric_{48}''}{PP48''}%
(see~\eqref{eq.diagram} for the notation), eight
combinatorial types
for $R=\bA_2^{12}$:
\begin{gather}
\cusptrue\diagram[Y_{56}]{
\rr{\.\-\*}\ \ \ \ \ \ \ \ \ \ \c\\
\>\>\>\>\>\>\R\R\R\\
\-\-\-\-\-\-\=\=\=\c\c
}\punct,\label{Y56}\\\noalign{\allowbreak}
\cusptrue\diagram[Q_{52}''']{
\c\ \ \ \ \ \ \ \ \ \ \rr{\*\-\*}\\
\>\>\>\>\>\>\R\R\R\c\\
\<\<\<\<\<\<\L\L\L\ \c
}\punct,\label{Q52'''}\\\noalign{\allowbreak}
\cusptrue\diagram[Y_{52}']{
\l\ \ \ \ \ \ \ \ \ \c\c\\
\>\>\>\>\>\>\\
\<\<\<\<\<\<\c\c\c\c
}\punct,\label{Y52'}\\\noalign{\allowbreak}
\cusptrue\diagram[Z_{50}]{
\r\ \ \ \ \ \r\\
\>\>\>\>\>\>\ \ \ \ \ \ \\
\ \ \ \ \ \ \>\>\>\>\>\>\\
\<\<\<\<\<\<\<\<\<\<\<\<
}\punct{$(\times2)$,}\label{Z50}\\\noalign{\allowbreak}
\cusptrue\diagram[*]{
\r\r\\
\>\>\>\>\>\>\ \ \ \ \ \ \\
\ \ \ \ \ \ \>\>\>\>\>\>\\
\<\<\<\<\<\<\<\<\<\<\<\<
}\punct{$(\times2)$,}\label{Z52}\\\noalign{\allowbreak}
\cusptrue\diagram[X_{52}'\bis]{
\l\l\\
\>\>\>\>\>\>\ \ \ \ \ \ \\
\ \ \ \ \ \ \>\>\>\>\>\>\\
\ \ \l\<\<\<\<\<\<\L\L\L
}\punct,\label{X52'2}
\end{gather}
\qlabel{Z48''}{Z_{48}''}{Z52}%
\qlabel{Z52}{Z_{52}}{Z52}%
(each of~\eqref{Z50} and~\eqref{Z52} represents two
combinatorial types, which differ
by the transposition $\bbc_1\leftrightarrow\bbc_2$), and one
combinatorial type
for $R=\bA_3^{8}$:
\begin{gather}
\tacnodetrue\diagram[Q_{52}''\bis]{
\C\C\\
\+\+\+\+\\
\ \ \ \ \+\+\+\+\\
\ \ \C\C\-\-\-\-
}\punct.\label{Q52''-2}
\end{gather}
In the diagrams for $\bA_2^{12}$,
unlike $\bA_1^{24}$,
we list the images of the basis vectors of $\rt(\bbT)$
\emph{in the first line}, and those of $\bbc_1,\bbc_2,\ldots$ in the other lines.
Representing each copy of~$\bA_2$ as
$\Z a_1+\Z a_2$, $a_1^2=a_2^2=2$, $a_1\cdot a_2=-1$, we use the following
notation:
\[
\begin{alignedat}3
&\bbl\bullet a_1,\quad
&&\bbr\bullet a_2,\quad
&&\bbcc\circ -a_1-a_2,\\
&\bbcc{\relbar\joinrel\joinrel\to}\tfrac13(2a_1+a_2),\quad
&&\bbcc{\leftarrow\joinrel\joinrel\relbar}\tfrac13(a_1+2a_2),\quad
&&\bbcc{\relbar\joinrel\relbar}\tfrac13(a_2-a_1),
\end{alignedat}
\label{eq.diagram.A2}
\]
a double line indicating the negation.
Certainly, within each column of a diagram,
the components should be regarded
up to simultaneous action of $\OG(\bA_2)$. A similar convention applies
to~$\bA_3^8$, where, in addition to roots,
we only use the following two elements:
\[
\bbccc{\relbar\joinrel\relbar\joinrel\relbar}\tfrac12(a_1+a_3),\quad
\bbccc{\relbar\joinrel\mathrel+\joinrel\relbar}\tfrac12(a_1+2a_2+a_3).
\label{eq.diagram.A3}
\]

\subsection{Proof of \autoref{th.others}}\label{proof.others}
Each of the $24$ combinatorial types $\bbT\into R\dual$ found in
\autoref{s.many.2} is represented by an isometry $\bbT\into N$.
Thus, we only need to classify these isometries up to the finer group
$\OG(N)$
and analyze the quartics obtained.

The group $\OG(N)$ is given by~\eqref{eq.G(N)}, and the stabilizers of the
corresponding codes $N/R\subset R\dual\!/R$ are well known and found, \eg,
in~\cite{Conway.Sloane}. Omitting the details (\cf.
the proof of \autoref{th.2xS}), we merely state that the first five
diagrams
\eqref{X52''}--\eqref{Q48}
in \autoref{tab.ambiguous} give rise to two orbits each, whereas
all other diagrams are represented by a single orbit each.

\table
\caption{Ambiguous combinatorial types}\label{tab.ambiguous}
\hbox to\hsize\bgroup\hss\vbox\bgroup\def\2{\rlap{$^2$}}
\halign\bgroup
 &\strut\ \hss$#$\hss\ \cr
\eqref{X52''}&\eqref{X52''2}&\eqref{Z48'}&
 \eqref{PP48''}&\eqref{Q48}&\eqref{Q52''}&\eqref{X52'}&
 \eqref{Z50}\2&\eqref{Z52}\2&\eqref{X52'2}&\eqref{Q52''-2}
 \cr\noalign{\kern2pt\hrule\kern1pt}
\quartic{Q54}&\quartic{Q54}&*&
 \quartic{Q52''}&\quartic{Q52''}&\quartic{Q52''}&\quartic{X52'}&
 \quartic{Z50}&*&\quartic{X52'}&\quartic{Q52''}
 \cr
\quartic{X52''}&\quartic{X52''}&\quartic{Z48'}&
 \quartic{PP48''}&\quartic{Q48}
 \cr\noalign{\kern2pt\hrule\kern1pt}
\multispan{11}\ \ \strut\vtop{\small\hbox{$*$ stands for the pair of configurations
 \quartic{Z52} and \quartic{Z48''}}
 \hbox{$\eqref{Z50}^2$ and $\eqref{Z52}^2$ represent
 two combinatorial types each}}\hss\cr
\crcr\egroup\egroup\hss\egroup
\endtable

For each orbit, we compute the lattice $S:=-\bbT^\perp$ and verify
Condition~\iref{i.model.birational} of \autoref{th.K3}.
Two lattices, \viz. \quartic{PP48} and \quartic{PP48''} (one of the orbits of
\eqref{PP48''}, see \autoref{tab.ambiguous}) fail; these lattices are
discussed in \autoref{proof.P1xP1} below.
For each of the remaining lattices, we compute the set~$\CS_h$ (see
\autoref{s.K3}) and, for each $\Gg\in\CS_h$, use
\autoref{th.K3}\iref{i.model.lines} to compute the configuration of lines in
the corresponding quartic. These configurations are used to distinguish or
identify (if the configuration is found in \cite{DIS}) lattices
obtained from distinct orbits. With the few exceptions listed in
\autoref{tab.ambiguous}, each lattice~$S$ is obtained from a unique orbit;
then, by \autoref{lem.index}, we have $[\Aut\discr S:\OG(S)]=1$
and \autoref{th.K3} gives us a unique spatial model.

Analyzing \autoref{tab.ambiguous}, we conclude that the index
$m:=[\Aut\discr S:\OG(S)]$ given by \autoref{lem.index} is greater than~$1$
in the following cases:
\roster*
\item
\quartic{Q54}, \quartic{X52'}, \quartic{X52''} with $m=2$ and
\quartic{Q52''} with $m=4$: these configurations of lines
and corresponding quartics are classified
in~\cite{DIS};
\item
\quartic{Z50} with $m=2$ and the pair $(\quartic{Z52},\quartic{Z48''})$
(represented by a $*$ in \autoref{tab.ambiguous}) with $m=3$:
these configurations are treated separately below.
\endroster
Note that,
for the $Z_*$ series, even those known in~\cite{DIS}, the results
of~\cite{DIS} do not apply directly: in this case, lines do not generate
$\TL^\perp\otimes\Q$ and, hence, \cite{DIS} establishes the connectedness of
a
$1$-parameter family of quartics of typical Picard rank~$19$ rather than the
uniqueness of any particular singular quartic (\cf. \autoref{prop.52}).

\subsubsection{The pair $(\quartici{Z52},\quartici{Z48''})$}
Let $\iota\:\bbT\into N$ be one of the \emph{three} isometries to which
\autoref{tab.ambiguous} assigns the pair of configurations
$(\quartic{Z52},\quartic{Z48''})$ ($*$ in the table), and consider the
lattice $S:=-\bbT^\perp$. In the notation of \autoref{s.K3}, we have a
splitting $\CS_h=\CS_{52}\cup\CS_{48}$, so that the quartic corresponding to
an element $\Gg\in\CS_n$ has $n$ lines. (If $n=52$, the configuration of
lines is identified with \quartic{Z52} in \cite{DIS}.)
Since $\ls|\CS_{52}|=4$ and $\ls|\CS_h|=12$, from \eqref{eq.orbits} we have
$[\Aut\discr S:\OG(S)]\ge3$; then, using \autoref{lem.index}, we conclude
that $[\Aut\discr S:\OG(S)]=3$ and all three lattices are isomorphic.
By \eqref{eq.orbits} again, $\CS_{52}$ is a single $\OG(S)$-orbit and
$[\Aut\Gg^\perp:\Stab\Gg]=1$ for each $\Gg\in\CS_{52}$; hence, there is a
unique quartic with this configuration of lines.

For the other configuration~\quartic{Z48''}, consider the only isometry
$\bbT\into N:=N(\bA_1^{24})$, see~\eqref{Z48'}, identify $\discr S$ and
$\discr\bbT$ accordingly, and let $G,H\subset\Aut\discr S$ be the images of
$\OG(\bbT)$ and $\OG(S)$, respectively, \cf. \autoref{s.aut.T}. Any
autoisometry of~$\bbT$ extends to~$N$
(due to the uniqueness of this extension); hence, $G\subset H$.
On the other hand, $\CS_{52}$ and $\CS_{48}$ are precisely the $G$-orbits
on~$\CS_h$ and the $G$-stabilizer of each element $\Gg\in\CS_{48}$ contains
$\pm\id\in\Aut\Gg^\perp$. It follows that $\CS_{48}$ is a single
$\OG(S)$-orbit and the image of $\Stab\Gg$ in $\Aut\Gg^\perp$
contains $\{\pm\id\}$,
which is the image of $\OG^+(T)$. Since
$[\Aut\Gg^\perp:\Stab\Gg]=2$ by~\eqref{eq.orbits}, we conclude that
the configuration
\quartic{Z48''} is realized by two complex conjugate quartics.

\subsubsection{The configuration \quartici{Z50}}
Consider the two isometries $\iota\:\bbT\into N:=N(\bA_1^{24})$ given
by~\eqref{Z50}; they differ by an autoisometry of~$\bbT$ and, hence, have the
same orthogonal complement $S:=-\bbT^\perp$. Let $G,H\subset\Aut\discr S$ be
as above, and let $G^\iota\subset G$ be the index~$2$ subgroup preserving the
combinatorial type, so that we have $G^\iota\subset H$.
The set $\CS_h$ splits into three $G^\iota$-orbits, each of length~$4$, and
for two of these orbits, the $G^\iota$-stabilizers of elements contain
$\pm\id$. On the other hand, we have $[\Aut\discr S:H]=2$ by
\autoref{lem.index} and the orbits of~$H$ are unions of
those of~$G^\iota$. Hence, in view of~\eqref{eq.orbits}, $\CS_h$ is a single
orbit and $[\Aut\Gg^\perp:\Stab\Gg]=2$ for each $\Gg\in\CS_h$, the image of
$\Stab\Gg$ containing $\pm\id$.
By \autoref{th.K3}, the lattice~$S$ gives rise to two complex conjugate quartics.
\qed

\subsection{Proof of \autoref{cor.P1xP1}}\label{proof.P1xP1}
We can represent $Q:=\Cp1\times\Cp1$ as a smooth quadric in~$\Cp3$.
Hence, any map $X(T)\to Q$ can be regarded as a spatial model
$\Gf\:X\to\Cp3$. According to
\cite[Theorem~5.2]{Saint-Donat} (see also \cite[Proposition~5.7]{Saint-Donat}),
such a map~$\Gf$
factors through
a quadric~$Q$ if and only if the corresponding set of data $(S,[\Gg])$
in \autoref{th.K3} violates the hypotheses of
Statement~\iref{i.model.birational} of the theorem.
If this is the case and $(2/d)\Gg$ is represented by a vector $e\in S\dual$,
$e^2=-1$, then the linear systems $\ls|\frac12h\pm e|$ are the pull-backs of
the two rulings of~$Q$.
Furthermore, Statement~\iref{i.model.nonsingular} can be restated as
\roster[3]
\item
the ramification locus $C\subset Q$ is smooth if and only if $S$ is root
free,
\endroster
and Statement~\iref{i.model.lines} holds literally if ``lines'' are
understood as smooth rational curves mapped isomorphically to generatrices
of~$Q$. (The image of each line is a bitangent of~$Q$, and
the pull-back of each bitangent consists of two lines. The two rulings differ
by the intersection $a\cdot e=\pm1$.)

The
complete list of root-free lattices satisfying the imposed conditions
on the discriminant
is found in \autoref{s.many.2};
for only two of these lattices,
\viz. \quartic{PP48} and \quartic{PP48''} (one of the orbits of
\eqref{PP48''}, see \autoref{tab.ambiguous}),
there is a class $\Gg\in\CS_h$ that violates the hypothesis of
\autoref{th.K3}\iref{i.model.birational}.
By \autoref{th.K3} combined with \autoref{lem.index} and~\eqref{eq.orbits},
each of these two lattices gives rise to a
single model $\Gf\:X(T)\to Q$ and this model has $48$ lines, \ie, $C$ has $24$
bitangents. On the other hand, a curve of bidegree $(4,4)$
in~$Q$ may have at most
$12$ bitangents in each ruling, as an elliptic pencil on a $K3$-surface may
have at most $12$ singular fibers of Kodaira type $\mathrm{I}_2$.
\qed

\remark\label{rem.P1xP1}
Using an appropriate version of \cite[Lemma 4.6]{DIS}, it is not very
difficult to show that all maximal (\ie, containing $48$ lines)
configurations of lines in
smooth hyperelliptic
models $X\to Q:=\Cp1\times\Cp1$ of $K3$-surfaces are isomorphic to each
other. The lifts to~$X$ of the $24$ bitangents to the ramification locus can
be chosen and ordered so as to intersect according to the pattern shown in
\autoref{fig.V}.
\figure
\hbox to\hsize{\hss$
\makeatletter
\setbox0\hbox{$\m@th\joinrel\relbar\joinrel\relbar\joinrel$}%
\def\zbox#1{\kern-\wd0\hbox to\wd0{\hss$#1$\hss}}%
\def\1{\zbox\bullet}\def\.{}\let\\\cr
\vcenter{\offinterlineskip\halign{&\copy0\zbox|#\cr
\1&\1&\1&\1&\1&\1&\1&\1&\1&\1&\1&\1\\
\1&\1&\1&\1&\1&\1&\1&\1&\.&\.&\.&\.\\
\1&\1&\1&\1&\.&\.&\.&\.&\1&\1&\.&\.\\
\1&\1&\1&\1&\.&\.&\.&\.&\.&\.&\1&\1\\
\1&\1&\.&\.&\1&\1&\.&\.&\1&\1&\.&\.\\
\1&\1&\.&\.&\1&\1&\.&\.&\.&\.&\1&\1\\
\1&\1&\.&\.&\.&\.&\1&\1&\1&\1&\1&\1\\
\1&\1&\.&\.&\.&\.&\1&\1&\.&\.&\.&\.\\
\1&\.&\1&\.&\1&\.&\1&\.&\1&\.&\1&\.\\
\1&\.&\1&\.&\1&\.&\1&\.&\.&\1&\.&\1\\
\1&\.&\.&\1&\.&\1&\1&\.&\1&\.&\1&\.\\
\1&\.&\.&\1&\.&\1&\1&\.&\.&\1&\.&\1\\}}
$\hss}

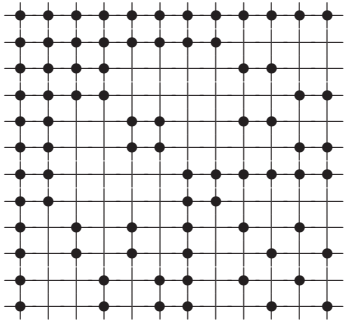
\caption{The maximal configuration of bitangents (see \autoref{rem.P1xP1})}\label{fig.V}
\endfigure
The lines span a rank~$17$ sublattice in $\NS(X)$;
hence,
there is a $3$-parameter equilinear family of models, and,
similar to \autoref{prop.52}, this family contains infinitely many singular
$K3$-surfaces. (Two of these singular $K3$-surfaces constitute
\autoref{cor.P1xP1};
another one is $X([8,0,12])$ discussed in \autoref{s.Mukai}.)
Note also that this maximal configuration admits a faithful transitive action
of the Mukai group~$\Mukai{H192}$, which is induced from the action on
$X([8,0,12])$
(see \autoref{th.Mukai} and remarks thereafter).
\endremark

\subsection{Oguiso pairs\pdfstr{}{ {\rm(see~\cite{Oguiso})}}}\label{s.Oguiso}
An \emph{Oguiso pair} is a pair of smooth spatial models $\Gf_i\:X\to\Cp3$,
$i=1,2$, of the same $K3$-surface~$X$ such that
one has
\[*
h_1\cdot h_2=6\quad\text{for}\quad
 h_i:=\Gf_i^*[\Cp2\cap\Gf_i(X)]\in\NS(X),\ i=1,2.
\]
According to \cite[Theorem 1.5]{Oguiso}, the two smooth
quartics $X_i:=\Gf_i(X)$ constituting an Oguiso pair are Cremona
equivalent and each quartic is a Cayley $K3$-surface.

\lemma\label{lem.Oguiso}
Consider a smooth quartic $X\in\Cp3$ and let $h\in\NS(X)$, $h^2=4$, be
its polarization. Then, another class $h'\in\NS(X)$ satisfying
$(h')^2=4$ and $h'\cdot h=6$ is very ample if and only if
\roster
\item\label{i.Oguiso.0}
there is no class $e\in\NS(X)$ such that $e^2=0$ and $e\cdot h'=2$\rom;
\item\label{i.Oguiso.2}
there is no class $e\in\NS(X)$ such that $e^2=-2$ and $e\cdot h'=0$\rom;
\item\label{i.Oguiso.nef}
one has $l\cdot h'>0$ for the class $l=[\ell]$ of each line $\ell\subset X$.
\endroster
\endlemma

\proof
We only need to show that $h'$ is nef; then, Conditions~\iref{i.Oguiso.0}
and~\iref{i.Oguiso.2} would imply that $h'$ is very ample (\cf.
\autoref{th.K3} and its proof). A class $h'\in\NS(X)$ is nef if and only if
$h'\cdot h>0$ (which is given) and $h'$ and~$h$ belong to the same
fundamental polyhedron of the subgroup of $\OG(\NS(X))$
generated by reflections, \ie, there is no root
$r\in\NS(X)$ such that $r\cdot h>0$ and $r\cdot h'<0$.
(Both $h$ and $h'$ are in the interior of a
fundamental polyhedron, due to the
fact that $h$ is very ample and Condition~\iref{i.Oguiso.2}, respectively).

Assume that there is a root~$r$ such that $\Ga:=r\cdot h\ge1$ and
$\Gb:=-r\cdot h'\ge1$ and consider the sublattice spanned by $h$, $h'$,
and~$r$. The determinant of the Gram matrix, which equals
$4(10-\Ga^2-3\Ga\Gb-\Gb^2)$, must be nonnegative, as the lattice is hyperbolic,
and we conclude that $\Ga=\Gb=1$, \ie, $r$ is the class of a line in~$X$.
\endproof

Let $(\NS,h)$ be one of the polarized N\'{e}ron--Severi lattices
listed in \autoref{tab.quartics}. (The transcendental lattice~$T$ is not
important in Oguiso's construction.)
In each case, it is straightforward to list all vectors $h'\in\NS$ satisfying
$(h')^2=4$ and $h'\cdot h=6$, use \autoref{lem.Oguiso} to select those that
are very ample, and compute the number of lines with respect to the new
polarization $h'$. (According to the table, this number identifies the
polarization within each lattice~$\NS$.)
With the two exceptions stated at the end of \autoref{s.det<80},
starting from any polarization~$h$,
we obtain all numbers of lines that are possible for the given lattice;
thus, with the same exceptions, any two smooth
models
of a $K3$-surface~$X$ as in \autoref{th.others} constitute an Oguiso pair.

\section{Other polarizations}\label{S.polarizations}

In this section, we consider other
(than quartics in $\Cp3$) polarizations of singular $K3$-surfaces
and prove \autoref{th.polarizations}.
Then, in \autoref{s.Mukai}, we discuss
projective models of the eleven Mukai groups.

\subsection{The set-up}\label{s.set-up}
Fix an even integer $2D>0$ and define a \emph{model} of a
singular $K3$-surface $X:=X(T)$ as
a map $\Gf\:X\to\Cp{D+1}$ defined by a fixed point free ample linear system
$\ls|h|$
of degree $h^2=2D$,
where $h:=h_\Gf\in\NS(X)$ stands for the class of a hyperplane section.
Two models~$\Gf_1$, $\Gf_2$ are
\emph{projectively equivalent} if there exists a pair of automorphisms
$a\:\Cp{D+1}\to\Cp{D+1}$
and $a_X\:X\to X$ such that $\Gf_2\circ a_X=a\circ\Gf_1$.
A model is \emph{smooth} if it does not contract a curve in~$X$, \ie, if the
pull-back of each point is finite. A \emph{line} in a model~$\Gf$ is a smooth
rational curve $C\subset X$ such that
the restriction $\Gf|_C$ is an isomorphism onto a line
in $\Cp{D+1}$.

Following \autoref{s.K3}, consider the lattice $S_\Gf:=h^\perp\subset\NS(X)$.
We have
\[
\discr S_\Gf\cong-\tker^\perp\!/\tker,\qquad
\ls|\discr S_\Gf|=2D\ls|\discr T|/(\depth\Gf)^2,
\label{eq.depth.D}
\]
where
\[
\tker\subset\discr T\oplus\discr\Z h,\qquad\tker\cap\discr T=0,
\label{eq.C-kernel.D}
\]
is a cyclic group; its order
$\depth\Gf:=\ls|\tker|$, called the \emph{depth} of~$\Gf$, divides~$2D$.
Fixing a lattice $S:=S_\Gf$ and, hence, the depth $d=d(S):=\depth\Gf$,
consider the sets
\begin{align*}
\CS_{dh}&:=\bigl\{\Gg\in\discr S\bigm|
 (2D/d)\Gg=0,\ \Gg^2=-d^2\!/2D\bmod2\Z\bigr\},\\
\CS_{dh}^+&:=\bigl\{\Gg\in\CS_{dh}\bigm|
 \CK_\Gg^\perp/\CK_\Gg\cong-\discr T\bigr\},
\end{align*}
where, for $\Gg\in\CS_{dh}$, the isotropic subgroup
$\CK_\Gg\subset\discr\Z h\oplus\discr S$ is generated by
$(d/2D)h\oplus\Gg$.
As usual, fixing an isometry $\CK_\Gg\cong-\discr T$, we regard both
stabilizers $\stab\Gg\subset\Aut\discr S$ and
$\Stab\Gg\subset\OG(S)$ acting
on the discriminant $\discr T$.

\autoref{th.K3} and its proof translate almost literally to
the general case.

\theorem\label{th.K3.D}
The projective equivalence classes of models
$\Gf\:X(T)\to\Cp{D+1}$
are in a one-to-one correspondence with the triples
consisting of
\roster*
\item
a negative definite lattice~$S$ of rank~$19$ and
$\discr S\cong-\tker^\perp\!/\tker$ as in~\eqref{eq.C-kernel.D},
\item
an $\OG(S)$-orbit $[\Gg]\subset\CS_{dh}^+$
\rom(where $d=d(S)$ is the depth, see~\eqref{eq.depth.D}\rom), and
\item
a double coset $c\in\OG^+(T)\backslash\Aut\discr T/\Stab\Gg$
\endroster
and such that
\roster
\item\label{i.D.moved}
$d>1$ or $d=1$ and $\Gg$ is not represented by a vector
$a\in S\dual$, $a^2=-1/2D$.
\endroster
Under this correspondence, the following statements hold\rom:
\roster[\lastitem]
\item\label{i.D.birational}
a model~$\Gf$ is birational onto its image if and only if
\roster*
\item
$D\ne1$ and either $D\ne4$ or $\tker\cap\discr\Z h=0$, and
\item
either $d>2$ or $d\le2$
and the class $(2/d)\Gg$ is not represented by a vector
$a\in S\dual$, $a^2=-2/D$\rom;
\endroster
\item\label{i.D.nonsingular}
a model~$\Gf$ is smooth if and only if
$S$ is root free\rom;
\item\label{i.D.lines}
the lines in a smooth model~$\Gf$ are in a
one-to-one correspondence with the vectors $a\in S\dual$,
$a^2=-(4D+1)/2D$,
representing~$\Gg$.
\done
\endroster
\endtheorem

The condition $\tker\cap\discr\Z h=0$ in
\autoref{th.K3.D}\iref{i.D.birational} means that, if $h^2=8$, the vector
$h\in\NS(X)$
must be primitive, see \cite[Theorem 5.2]{Saint-Donat}.
Unlike \autoref{th.K3},
in \autoref{th.K3.D}\iref{i.D.nonsingular},~\iref{i.D.lines}
we no longer require that the model should be birational. As
in \autoref{s.polarizations} (\cf. also \autoref{cor.P1xP1}
and its proof in \autoref{proof.P1xP1}),
the smoothness of a hyperelliptic model $\Gf\:\X\to\Cp{D+1}$ is
understood as the smoothness of its ramification locus~$C$, and \emph{lines} are
defined as rational curves $L\subset X$ that project isomorphically onto lines
in~$\Cp{D+1}$. (Typically, the images $\Gf(L)$ are
generatrices of the scroll $\Gf(X)$ that have even
intersection index with
the ramification locus $C\subset\Gf(X)$
at each point of intersection, \cf. bitangents in
\autoref{cor.P1xP1} and tritangents in \autoref{s.polarizations}.
Each generatrix with this property splits into two lines in~$X$.)

\lemma\label{lem.octic}
In the notation of \autoref{th.K3.D}, assume that $D=4$ and that the model
$\Gf\:X\to\Cp5$ given by a triple $(S,[\Gg],c)$ is birational.
Then the image $\Gf(X)$ is an intersection of quadrics if and only if
$3\Gg$ is not represented by a vector $u\in S\dual$ such that
$u^2=-\frac98$ and $u\cdot a\ge-\frac38$ for each
vector $a\in S\dual$ as in
\autoref{th.K3.D}\iref{i.D.lines}.
\endlemma

\proof
According to \cite[Theorem~7.2]{Saint-Donat}, the defining ideal of
$\Gf(X)$ is generated by its elements of degree~$2$ if and only if there is
no nef class $e\in\NS(X)$ such that $e^2=0$ and $e\cdot h=3$.
Arguing as in the proof of \autoref{lem.Oguiso}, one can easily show that a
class $e\in\NS(X)$ satisfying $e^2=0$ and $e\cdot h=3$ is nef if and only if
one has $e\cdot\ell\ge0$ for each curve $\ell\subset X$ that
is either a line
($\ell\cdot h=1$)
or an exceptional divisor ($\ell\cdot h=0$).
For the latter, one can
map~$e$
to the distinguished Weyl chamber of $-S$
by an appropriate element of the Weyl group
(\cf. \autoref{proof.singular} below).
Then, the requirement $e\cdot\ell\ge0$ can be extended to all, not necessarily
irreducible, $(-2)$-curves $\ell\subset X$ satisfying $\ell\cdot h=1$.
The statement of the lemma is a translation of this
latter condition in terms of $S$
and~$\Gg$.
\endproof

\subsection{Proof of \autoref{th.polarizations}}\label{proof.polarizations}
Fix an integer $h^2=2D=2$, $6$, or~$8$ and a positive definite
even lattice~$\TL$ of rank~$2$. As in \autoref{s.reduction}, we can find a
\emph{test lattice}~$\bbT$, which is a positive definite even lattice of
rank~$5$ satisfying the identity
\[*
\discr\bbT\cong{-\discr\TL}\oplus{-\discr\Z h}.
\]
Then, an analog of \autoref{lem.genus} holds: any lattice~$S$ as in
\autoref{th.K3.D} is of the form $S=-\bbT^\perp$ for an appropriate
isometry $\bbT\into N$ to a Niemeier lattice~$N$.
(\latin{A priori}, this isometry does not need to be primitive, \cf.
\autoref{rem.genus}.) Since we are only interested in smooth models,
an additional
requirement is that there should be no root $r\in\rt(N)$ orthogonal to~$\bbT$.

Assuming that $\det\TL$ is bounded as in \autoref{th.polarizations}, the test
lattice~$\bbT$ (depending on~$h^2$)
can be chosen according to \autoref{conv.V}.
Then, arguing as in \autoref{s.many.roots}, we can easily use \GAP~\cite{GAP4}
and eliminate all Niemeier lattices~$N$ with $\rt(N)$ consisting of six or
less irreducible components. The remaining three lattices are treated as in
\autoref{s.many.2}, by using~\eqref{eq.A3.type}, \eqref{eq.A2.type},
\eqref{eq.A1.type} and enumerating the combinatorial types of isometries
first. Then, the isometries are classified up to the finer group
$\OG(N)$ given by~\eqref{eq.G(N)} and the resulting models are studied using
\autoref{th.K3.D}.


\subsubsection{Planar models, $h^2=2$}\label{ss.curves}
If $\det\TL\le116$, there is one combinatorial type, \viz. the isometry
$\bbT\into(\bA_1^{24})\dual$ given by the diagram
(see~\eqref{eq.diagram} for the notation)
\[*
\diagram[\2_{144}]{
\n\n\n\\
\-\-\-\-\-\-\-\-\-\-\-\-\ \ \ \ \ \ \ \ \ \ \ \ \\
\ \ \ \ \ \ \ \ \ \ \ \ \-\-\-\-\-\-\-\-\-\-\-\-
}\punct.\label{2-144}
\]
It gives rise to a unique isometry $\bbT\into N$; hence, due to \autoref{th.K3.D}
combined with \autoref{lem.index} and \eqref{eq.orbits}, there is a unique
model.

\subsubsection{Sextic models in $\Cp4$, $h^2=6$}\label{ss.sextics}
If $\det\TL\le48$,
there is one combinatorial type of isometries $\bbT\into(\bA_1^{24})\dual$
(see~\eqref{eq.diagram} for the notation):
\[*
\diagram[*]{
\ \ \ \ \ \ \ \ \ \ \ \ \n\ \ \ \ \ \ \ \ \ \n\n\\
\-\-\-\-\-\-\-\-\-\-\-\-\ \ \ \ \ \ \ \ \ \ \ \ \\
\-\-\-\-\-\-\ \ \ \ \ \ \-\-\-\-\-\-\-\-\-\-\ \
}\punct{}\label{6-36'}
\]
and four combinatorial types of isometries $\bbT\into(\bA_2^{12})\dual$
(see~\eqref{eq.diagram.A2}):
\begin{gather}
\cusptrue\diagram[\6_{42}]{
\ \ \ \ \ \ \ \ \ \ \c\rr{\.\-\*}\\
\>\>\>\>\>\>\R\R\R\\
\<\<\<\<\<\<\L\L\L\c
}\punct,\label{6-42}\\\noalign{\allowbreak}
\cusptrue\diagram[\6_{38}]{
\r\ \ \ \ \ \ \ \ \ \c\c\\
\>\>\>\>\>\>\R\R\R\\
\<\<\<\<\<\<\L\L\L\c
}\punct,\label{6-38}\\\noalign{\allowbreak}
\cusptrue\diagram[\theA\bis]{
\ \ \ \ \ \ \ \ \ \r\c\c\\
\>\>\>\>\>\>\R\R\R\\
\<\<\<\<\<\<\L\L\L\c
}\punct,\label{6-36'-2}\\\noalign{\allowbreak}
\cusptrue\diagram[\6_{36}'']{
\ \ \ \ \ \ \ \ \ \ \c\rr{\.\-\*}\\
\>\>\>\>\>\>\R\R\R\r\\
\<\<\<\<\<\<\L\L\L\c
}\punct.\label{6-36''}
\end{gather}
\qlabel{6-42'}{6_{42}'}{6-36'}%
\qlabel{6-36'}{6_{36}'}{6-36'}%
Each combinatorial type gives rise to a unique $\OG(N)$-orbit of isometries,
and, with the exception of~\eqref{6-36'}, \eqref{6-36'-2}, all orthogonal
complements are pairwise distinct: this fact can be established by computing
the configurations of lines. (In~\eqref{6-36''}, the subset
$\CS_h^+\subset\CS_h$ is proper.) As above, each of these lattices gives us
a single model.

In each of the exceptional cases~\eqref{6-36'}, \eqref{6-36'-2}, the set
$\CS_h$ splits into two subsets, with the configurations of
lines~\quartic{6-42'} or~\quartic{6-36'}.
Using \autoref{lem.index} and~\eqref{eq.orbits},
we conclude that the two orthogonal
complement are isomorphic to each other (as otherwise
the group
$\OG(S)$ would
have to
be transitive
on~$\CS_h^+$) and each configuration is realized by a single model.

\subsubsection{Octic models in $\Cp5$, $h^2=8$}\label{ss.octics}
If $\det\TL\le40$,
there are five combinatorial types of isometries $\bbT\into(\bA_1^{24})\dual$
(see~\eqref{eq.diagram} for the notation):
\begin{gather}
\diagram[\8_{32}]{
\n\n\.\.\\
\-\-\-\-\-\-\-\-\ \ \ \ \ \ \ \ \ \ \ \ \ \ \ \ \\
\ \ \ \ \ \ \ \ \-\-\-\-\-\-\-\-\-\-\-\-\-\-\-\-
}\punct,\label{8-32}\\\noalign{\allowbreak}
\diagram[\8_{36}]{
\n\n\ \ \-\-\-\-\-\-\-\-\ \ \ \ \ \ \ \ \ \ \ \ \\
\-\-\-\-\-\-\-\-\\
\ \ \ \ \ \ \ \ \=\=\-\-\-\-\-\-\-\-\-\-\-\-\-\-
}\punct,\label{8-36}\\\noalign{\allowbreak}
\diagram[\8_{36}\bis]{
\n\n\ \ \-\-\-\-\-\-\-\-\ \ \ \ \ \ \ \ \ \ \ \ \\
\-\-\-\-\-\-\-\-\ \ \ \ \\
\ \ \ \ \=\=\-\-\ \ \ \ \-\-\-\-\-\-\-\-\-\-\-\-
}\punct,\label{8-36-2}\\\noalign{\allowbreak}
\diagram[*]{
\ \ \ \ \ \ \ \ \ \ \ \ \n\ \ \ \ \ \ \ \ \ \n\n\\
\-\-\-\-\-\-\-\-\-\-\-\-\ \ \ \ \ \ \ \ \ \ \ \ \\
\-\-\-\-\-\-\ \ \ \ \ \ \-\-\-\-\-\-\-\-\-\-\ \
}\punct,\label{8-33}\\\noalign{\allowbreak}
\diagram[\8_{30}]{
\-\-\-\-\-\-\-\-\-\-\-\-\n\ \ \ \ \n\\
\ \ \ \ \ \ \ \ \ \ \-\-\-\-\-\-\-\-\ \ \ \ \ \ \\
\ \ \ \ \ \ \-\-\-\-\ \ \ \ \ \ \=\-\-\-\-\-\-\-
}\punct{}\label{8-30}
\end{gather}
\qlabel{8-36'}{8_{36}'}{8-33}%
\qlabel{8-33}{8_{33}}{8-33}%
and two combinatorial types of isometries $\bbT\into(\bA_2^{12})\dual$
(see~\eqref{eq.diagram.A2}):
\begin{gather}
\cusptrue\diagram[\8_{32}']{
\r\ \ \ \ \ \ \ \ \ \c\c\\
\>\>\>\>\>\>\R\R\R\\
\<\<\<\<\<\<\L\L\L\c
}\punct,\label{8-32'}\\\noalign{\allowbreak}
\cusptrue\diagram[\theA\bis]{
\ \ \ \ \ \ \ \ \ \r\c\c\\
\>\>\>\>\>\>\R\R\R\\
\<\<\<\<\<\<\L\L\L\c
}\punct.\label{8-33-2}
\end{gather}
Note that diagrams~\eqref{8-33}, \eqref{8-32'}, and~\eqref{8-33-2} are
identical to~\eqref{6-36'}, \eqref{6-38}, and~\eqref{6-36'-2}, respectively.
This coincidence
is due to the fact that the lattice $[6,0,6]\oplus\Z h_8$, $h_8^2=8$, is
obviously isomorphic to $[6,0,8]\oplus\Z h_6$, $h_6^2=6$; hence, they have
the same sets of root-free orthogonal complements.
For a similar reason, diagrams~\eqref{8-32}, \eqref{8-36}, and~\eqref{8-36-2}
are identical to~\eqref{X48}, \eqref{X56}, and~\eqref{X56-2}, respectively.

Each combinatorial type gives rise to a unique $\OG(N)$-orbit of isometries.
For~\eqref{8-32}, \eqref{8-30}, and~\eqref{8-32'}, this fact immediately
implies the uniqueness of the corresponding model.
The pair~\eqref{8-33}, \eqref{8-33-2} is treated similar
to~\eqref{6-36'}, \eqref{6-36'-2} in \autoref{ss.sextics}: we obtain two
configurations of lines, \quartic{8-36'} and~\quartic{8-33}, each realized by
a single model. Finally, the lattice~$S$ corresponding to~\eqref{8-36}
and~\eqref{8-36-2} is $\bS_{56}$ in \autoref{s.S56}. By
\autoref{lem.S56.aut}\iref{i.H56.8}, the set $\CS_h=\CS_{56}$ is a single
$\OG(S)$-orbit and, by \autoref{lem.S56.aut}\iref{i.H56.refl},
$\Stab\Gg\supset\{\pm\id\}$ for each $\Gg\in\CS_h$. Hence, there are two
complex conjugate models.
\qed

\subsection{Mukai groups}\label{s.Mukai}
A \emph{Mukai group} is a maximal (with respect to inclusion) finite group
admitting a faithful symplectic action on a $K3$-surface.
There are $11$ Mukai groups (see \cite{Mukai,Kondo}); they are listed in
\autoref{tab.Mukai}.
\table
\caption{Models of Mukai groups (see \autoref{th.Mukai})}\label{tab.Mukai}
\def\hh#1{\vcenter{\offinterlineskip\def\\##1##2{\hbox{\strut##1\
 ##2}}#1}}
\def\nolines#1{\\{$h^2=#1$}{}}
\def\+{\noalign{\kern2\lineskip}}
\hbox to\hsize\bgroup\hss\vbox\bgroup
\halign\bgroup\strut\quad\hss$#$\hss\quad&\hss$#$\hss\quad&$#$\hss\quad\cr
\noalign{\hrule\vspace{2pt}}%
G&\text{Combinatorial type of $N^G\into N$}&\text{Models, remarks}\cr
\noalign{\vspace{1pt}\hrule\vspace{3pt}}
\Mukaigroup{L2-7}{L_2(7)}&\diagram{
\-\-\-\-\-\-\-\-\\
\ \ \ \ \ \ \ \ \-\-\-\-\-\-\-\-\\
\n\ \ \ \ \ \ \ \n\ \ \ \ \ \ \ \-\-\-\-\-\-\-\-
}&h^2=2,4,8\cr\+
\Mukaigroup{A6}{\AG6}&\diagram{
\n\n\ \ \ \ \ \ \n\\
\-\-\-\-\-\-\-\-\\
\ \ \ \ \ \ \ \ \-\-\-\-\-\-\-\-\-\-\-\-\-\-\-\-
}&h^2=2,6,8\cr\+
\Mukaigroup{S5}{\SG5}&\diagram{
\n\.\.\ \ \ \ \ \n\\
\-\-\-\-\-\-\-\-\\
\ \ \ \ \ \ \ \ \-\-\-\-\-\-\-\-\-\-\-\-\-\-\-\-
}&h^2=4,6\cr\+
\Mukaigroup{M20}{M_{20}}&\diagram{
\n\n\n\\
\-\-\-\-\-\-\-\-\\
\ \ \ \ \ \ \ \ \-\-\-\-\-\-\-\-\-\-\-\-\-\-\-\-
}&h^2=4,8\cr
\Mukaigroup{F384}{F_{384}}&\diagram{
\n\n\.\.\\
\-\-\-\-\-\-\-\-\\
\ \ \ \ \ \ \ \ \-\-\-\-\-\-\-\-\-\-\-\-\-\-\-\-
}&\hh{\\{$h^2=4$:}{\quartic{X48}}
 \\{$h^2=8$:}{\quartic{8-32}}}\cr
\Mukaigroup{A4-4}{\AG{4,4}}&\diagram{
\n\n\.\.\.\\
\-\-\-\-\-\-\-\-\\
\ \ \ \ \ \ \ \ \-\-\-\-\-\-\-\-\-\-\-\-\-\-\-\-
}&h^2=8\cr
\Mukaigroup{T192}{T_{192}}&\diagram{
\-\-\-\-\-\-\-\-\\
\ \ \ \ \ \ \ \ \-\-\-\-\-\-\-\-\\
\n\n\ \ \ \ \ \ \ \ \ \ \ \ \ \ \-\-\-\-\-\-\-\-
}&\hh{\\{$h^2=4$:}{\quartic{X64}}
 \nolines8}\cr
\Mukaigroup{H192}{H_{192}}&\diagram{
\n\.\.\.\\
\-\-\-\-\-\-\-\-\\
\ \ \ \ \ \ \ \ \-\-\-\-\-\-\-\-\\
\ \ \ \ \ \ \ \ \ \ \ \ \ \ \ \ \-\-\-\-\-\-\-\-\\
}&\hh{\\{$h^2=4$:}{\autoref{fig.V}}
 \\{$h^2=8$:}{$(6)^{32}$}}\cr
\Mukaigroup{N72}{N_{72}}&\diagram{
\n\.\.\\
\ \-\-\-\-\-\-\-\-\\
\ \-\-\ \ \ \ \ \ \-\-\-\-\-\-\\
\-\-\-\ \ \ \ \ \ \ \ \ \ \ \ \-\-\-\-\-\-\-\-\-\\
}&h^2=6\cr
\Mukaigroup{M9}{M_9}&\diagram{
\n\n\n\\
\-\-\-\-\-\-\-\-\-\-\-\-\\
\ \ \ \ \ \ \ \ \ \ \ \ \-\-\-\-\-\-\-\-\-\-\-\-
}&\hh{\\{$h^2=2$:}{\quartic{2-144}}
 \nolines8}\cr
\Mukaigroup{T48}{T_{48}}&\diagram{
\n\.\.\.\\
\-\-\-\-\-\-\-\-\\
\ \ \ \ \-\-\-\-\-\-\-\-\-\-\-\-\\
\ \ \ \ \-\-\-\-\ \ \ \ \ \ \ \ \-\-\-\-\-\-\-\-\\
}&\hh{\\{$h^2=2$:}{$108$ lines}
 \nolines8}
\ealign\hss\egroup
\endtable


Each Mukai group~$G$ acts on the lattice $L:=H_2(X)$. The invariant
sublattice~$L^G$ is positive definite of rank~$3$ ($G$ preserves the
holomorphic form by the definition, and, since $G$ is finite, one can find an
invariant K\"{a}hler form); hence, the \emph{coinvariant lattice}
$L_G:=(L^G)^\perp$ is negative definite of rank~$19$.
The lattice~$L_G$ is root free and its isomorphism class
determines and is determined by the group~$G$.
The number of isomorphism classes of
isometries $L_G\into L$ and, hence, actions of~$G$ on~$L$ (see
\autoref{cor.extension}; by definition, $G$ acts identically
on $\discr L_G$ and, hence, extends to any overlattice),
is two for the first three groups and one for the others (see, \eg,
\cite{Hashimoto}; an isometry is determined by
the orthogonal complement $L^G$).
It follows that we have one or two $G$-equivariant
$1$-parameter families of $K3$-surfaces; generic members of these families
are not algebraic, and the algebraic ones are singular.

A \emph{projective model} of a Mukai group~$G$ is a model $\Gf\:X\to\Cp{n}$ on
which $G$ acts by projective transformations. Since the lattice
$L_G=h^\perp\subset\NS(X)$ is root free, any such model is smooth, although
it may be hyperelliptic.
According to~\cite{Mukai,Kondo}, each coinvariant lattice $-S:=-L_G$ admits an
equivariant isometry $-S\into N:=N(\bA_1^{24})$, embedding~$G$ to~$M_{24}$.
The combinatorial types of the invariant lattices $N^G=(-S)^\perp$ can easily
be described using the orbit structure of $G\subset M_{24}$ given
in~\cite{Kondo}; they are listed in \autoref{tab.Mukai}
(see~\eqref{eq.diagram} for the notation).
This construction gives us a convenient
description of~$L_G$ and, together with \autoref{th.K3.D},
leads to the following statement.

\theorem\label{th.Mukai}
Each Mukai group~$G$ admits a projective model $\Gf\:X\to\Cp{n}$ of degree
$h^2=2n-2\le8$, see \autoref{tab.Mukai}.
Either one has $\depth(\Gf)>1$ and $\Gf(X)$ contains no lines, or
\rom(in seven cases\rom)
$G$ acts faithfully on the set of lines in~$\Gf(X)$.
\done
\endtheorem

The degree~$4$ model of~$\Mukai{H192}$ is hyperelliptic,
$X([8,0,12])\to\Cp1\times\Cp1$,
with the maximal
configuration of lines (see \autoref{rem.P1xP1}
and \autoref{fig.V}). The configurations of lines
in the octic models of~$\Mukai{F384}$ and~$\Mukai{H192}$ are isomorphic to each other;
hence, the same configuration admits a faithful action of both groups.
The action of~$G$ on the set of lines
is transitive with the exception of the following
three cases:
\roster*
\item
on \quartic{X64}, there are two orbits distinguished by the type of the
lines, see \autoref{s.pencils};
\item
on \quartic{2-144}, there are two orbits interchanged by the hyperelliptic
involution;
\item
on the degree~$2$ model of~$\Mukai{T48}$, there are three orbits:
two (of length~$48$) are interchanged by the hyperelliptic involution~$\tau$,
and one is $\tau$-invariant.
\endroster

\section{Examples}\label{S.examples}

In this section, we construct explicit examples to prove
Theorems~\ref{th.infty}, \ref{th.Q}, \ref{th.singular}
and \autoref{prop.52} and discuss a few other interesting examples of
quartics.

\subsection{Lattices without short vectors}\label{s.no.short}
Consider the test lattice with the Gram matrix
\[*
V_m:=\bmatrix
    2&  0&  0&  1&  0\\
    0&  2&  0&  1&  0\\
    0&  0&  2&  0&  0\\
    1&  1&  0&  4&  0\\
    0&  0&  0&  0&  m
\endbmatrix,
\]
where the even integer $m>0$ is to be specified later.
The group
\[*
\CS_m:=-\discr V_m\cong\Z/2\oplus\Z/4\oplus\Z/3\oplus\Z/m
\]
is generated by pairwise orthogonal elements of squares
$\frac12,\frac14,\frac23,-\frac1m\bmod2\Z$,
and
\[
\ls|\Aut\CS_m|\le4\cdot2^{\Omega(m)}\le4m,
\label{eq.size.Aut}
\]
where $\Omega(m)\le\log_2m$ stands for the number of prime divisors of~$m$, counted
with multiplicity.
(Indeed, if
$p\mathrel|m$ and $p>3$, then the group $\CS_m\otimes\Z_p$ is cyclic and one
has
$\Aut(\CS_m\otimes\Z_p)=\{\pm\id\}$. The remaining two
groups $\Aut(\CS_m\otimes\Z_2)$
and $\Aut(\CS_m\otimes\Z_3)$ are easily bounded.)

According to \cite[Theorem 1.10.1]{Nikulin:forms},
if $m\ne0\bmod3$, then there exists an even positive
definite lattice~$T_m$ of rank~$2$ such that
\[
\CS_m\cong{\discr T_m}\oplus{\discr\Z h},\quad h^2=4.
\label{eq.Tm}
\]
(If $m\ne0\bmod3$, all conditions in Nikulin's theorem hold trivially. Note
that a lattice~$T_m$ does \emph{not} exist if $m=0\bmod3$.)

Define the \emph{genus $\genus_m$} as the set of isomorphism classes of
even positive definite
lattices~$S$ of rank~$19$ such that $\discr S\cong\CS_m$.
Our proof of \autoref{th.infty} is based on the following statement.

\lemma\label{lem.infty}
For any pair of integers $M,Q>0$, there exists a positive even integer
$m\ne0\bmod3$ and $M$ pairwise
distinct representatives $S\in{\genus_m}$
with the property that $a^2>Q$ for each vector $a\in S\sminus0$.
\endlemma

\proof
We will construct sufficiently many primitive embeddings
$V_m\into N$, where $N:=N(\bA_1^{24})$, and let $S:=V_m^\perp$.
Then $S\in{\genus_m}$ by \autoref{cor.extension}.

Choose an orthogonal basis $\be_1,\ldots,\be_{24}$, $\be_i^2=2$, for
$\rt(N)$ and denote
\[*
N_k:=\bigl((\Z\be_k\oplus\ldots\oplus\Z\be_{24})\otimes\Q\bigr)\cap N
\]
for $k=1,\ldots,24$. Let, further,
\[*
\fQ_k:=\bigl\{a\in N_k\sminus0\bigm|a^2\le Q\bigr\};
\]
these are finite sets independent of~$M$.

The first four basis elements of $V_m$ are mapped according to the
diagram
\[*
\diagram{
\n\n\n\ \ \ \ \ \ \ \ \ \ \ \ \ \ \ \ \ \ \ \ \ \\
\ \-\-\ \ \-\-\-\-\-\-
}
\]
(see~\eqref{eq.diagram} for the notation); certainly, we assume that the
basis $\{\be_k\}$ is ordered so that the
image~$\bba$ of the square~$4$ vector is an octad in the Golay code
$N/\!\rt(N)$.
Let $\barN_k:=\bba^\perp\cap N_k$ and $\bar\fQ_k:=\bba^\perp\cap\fQ_k$;
obviously, $V_m^\perp\subset\barN_4$. The fifth basis vector is to be mapped
to $\bbc:=\sum_{k=4}^{24}c_k\be_k$, where the coefficients $c_4,c_5$ will be
fixed and the others will vary. (We will use the
``free'' coefficients $c_4,c_5$ to
adjust the number theoretical properties of~$\bbc$.)

Let $H\subset\barN_6\otimes\R$ be the finite union of hyperplanes $a^\perp$,
$a\in\bar\fQ_6$.
The number of integral points in the ball $B_r\subset\barN_r\otimes\R$
of radius~$r$ grows
as $O(\QOPNAME{vol}B_r)=O(r^{18})$,
whereas the number of points in $B_r\cap H$ grows as
$O(r^{17})$.
Subtracting and passing to spheres, we find that there is a sequence of
integers $s_n:=r_n^2\to\infty$ such that
\[*
\ls|\fC_n|\ge C_1r_n^{17},\quad\text{where}\quad
\fC_n:=\bigl\{u\in\barN_6\sminus H\bigm|u^2=s_n\bigr\}.
\]
(Here and below, $C_i:=C_i(Q)$ are positive constants independent of~$M$ and~$n$.)

Each coordinate of each vector $u\in\fC_n$ is bounded by~$r_n$; hence,
$\ls|u\cdot a|\le C_2r_n$ for all $u\in\fC_n$ and
$a\in\bar\fQ_5\sminus\bar\fQ_6$.
Since each vector $a\in\bar\fQ_5\sminus\bar\fQ_6$ has a nonzero coordinate
at~$\be_5$,
we have $(c_5\be_5+u)\cdot a\ne0$ for all $u\in\fC_n$ and
$a\in\bar\fQ_5$ whenever $c_5>C_2r_n$.
Fix an integer~$c_5$ with this property; we can assume that
$\ls|c_5|\le C_3r_n$.
In a similar way, we can find a positive integer
$c_4\le C_4r_n$ such that
$(c_4\be_4+c_5\be_5+u)\cdot a\ne0$ for all $u\in\fC_n$ and $a\in\bar\fQ_4$.
By slightly stretching the bounds
(say, replacing $C_3$ and~$C_4$ with $C_3+3$ and $C_4+C_3+6$, respectively),
we can also assume that
$c_4,c_5$ are coprime and of opposite parity and that the common square
\[*
m:=(c_4\be_4+c_5\be_5+u)^2=2c_4^2+2c_5^2+s_n\ne0\bmod3.
\]

Now,
taking any vector $c_4\be_4+c_5\be_5+u$, $u\in\fC_n$,
for the image~$\bbc$ of the fifth
generator, we obtain a primitive embedding $\bbT\into N$ such that
$\bbT^\perp\cap\fQ_1=\varnothing$. Hence, in view of~\eqref{eq.size.Aut},
combined with the bound $m\le(2C_3^2+2C_4^2+1)r_n^2$,
and
\autoref{cor.extension}, it suffices to choose $s_n=r_n^2$ so that
\[*
C_1r_n^{15}\ge4M(2C_3^2+2C_4^2+1)\,\ls|\OG(N)|
\]
to obtain at least $M$ distinct isomorphism classes
of orthogonal complements.
\endproof

\subsection{Proof of \autoref{th.infty}}\label{proof.infty}
The first statement is immediate:
we have finitely
many genera for the lattice $S:=h^\perp\subset\NS(X)$
(since $\ls|\det S|\le4\det T$ is bounded), each genus contains
finitely many isomorphism classes, and, for each
class, the extensions
$H_2(X)\supset\NS(X)\supset S$ are determined by a finite set of data.

For the second statement, we take for~$T$ the lattice~$T_m$
as in~\eqref{eq.Tm}, where
$m$ is given by \autoref{lem.infty} with $Q=932$.
Let $S_1,\ldots,S_M$ be $M$ distinct lattices given by the lemma.
By \autoref{th.K3}, each lattice~$-S_i$ gives rise to at least one spatial
model $\Gf_i\:X:=X(T)\to\Cp3$, and this model is smooth. (To ensure that
$\Gf_i$ satisfy the
hypotheses of Statements~\iref{i.model.birational}
and~\iref{i.model.nonsingular} of the theorem, it
would suffice to let $Q=4$
and~$2$, respectively, in \autoref{lem.infty}.)
Let $h_i\in\NS(X)$ be the hyperplane section class corresponding to~$\Gf_i$,
so that $h_i^\perp\cong{-S_i}$.
Any projective equivalence between~$\Gf_i$ and~$\Gf_j$ induces an
automorphism of $\NS(X)$ taking~$h_i^\perp$ to~$h_j^\perp$; hence,
$i=j$. According to \cite[Proposition 1.7]{Oguiso},
if the smooth quartic $\Gf_i(X)$ is taken to another smooth quartic by a Cremona
transformation $\Cp3\dashrightarrow\Cp3$ \emph{that is not regular on~$\Cp3$},
then $X$ contains a reduced irreducible curve~$C$ such that
$h_i\cdot[C]<16$ and the classes $h_i,[C]\in\NS(X)$ are linearly independent.
Letting $d:=h_i\cdot[C]\le15$, we would obtain a class
\[*
0\ne a:=4[C]-dh_i\in h_i^\perp\cong{-S_i}\quad\text{with}\quad
-a^2=4d^2-16C^2\le932
\]
(since $C^2\ge-2$),
which would contradict our choice of the lattices~$S_i$.
\qed

\remark\label{rem.growth}
The proof of \autoref{lem.infty} shows that, in fact, we have constructed
a sequence of lattices $T_m$, $\det T_m\to\infty$, such that
the number of distinct spatial models of the $K3$-surfaces
$X(T_m)$ grows at least as fast as $O\bigl((\det T_m)^{15/2}\bigr)$.
\endremark

A statement similar to \autoref{th.infty}, \ie, the fact that the
maximal number of
projective equivalence classes of smooth models of a fixed singular
$K3$-surface is not bounded, holds for the other polarizations, and the proof
is similar to that of \autoref{th.infty}.
(An assertion on the Cremona equivalence would need an
analogue of Proposition 1.7 in~\cite{Oguiso}.)
For example, for the three
polarizations $h^2=2D=2,6,8$ considered in \autoref{th.polarizations}, one
can start with the lattice
\[*
V_{m,D}:=V_D\oplus\Z a,\quad a^2=m,
\]
where $V_D\subset N(\bA_1^{24})$ is the rank~$4$ sublattice described by one
of the following diagrams (see~\eqref{eq.diagram} for the notation):
\begin{gather*}
\diagram[D=1]{
\n\n\ \ \ \ \ \ \ \ \ \ \ \ \ \ \ \ \ \ \ \ \ \ \\
\-\-\ \ \-\-\-\-\-\-\\
\ \ \ \ \ \ \ \ \-\-\-\-\-\-\-\-
}\punct,\\\noalign{\allowbreak}
\diagram[D=3]{
\n\n\n\ \ \ \ \ \ \ \ \ \ \ \ \ \ \ \ \ \ \ \ \ \\
\ \ \-\ \ \-\-\-\-\-\-\-\-\-\-\-\-\-\-\-
}\punct,\\\noalign{\allowbreak}
\diagram[D=4]{
\n\n\ \ \ \ \ \ \ \ \ \ \ \ \ \ \ \ \ \ \ \ \ \ \\
\ \-\ \ \-\-\-\-\-\-\-\\
\ \ \ \ \ \ \ \ \ \-\-\-\-\-\-\-\-\-\-\-\-
}\punct.
\end{gather*}
Then, an even positive definite lattice $T_{m,D}$ of rank~$2$ such that
\[*
-\discr V_{m,D}\cong{\discr T_{m,D}}\oplus{\discr\Z h},\quad h^2=2D,
\]
exists whenever $m\ne0\bmod11$ (if $D=1$) or $m\ne0\bmod5$ (in the two other
cases).
Denoting by $\genus_{m,D}$ the set of isomorphism classes of
even positive definite
lattices of rank~$19$ and discriminant $-\discr V_{m,D}$, we have an
obvious literate
analogue of \autoref{lem.infty}
(whose proof is based on very rough estimates) and, hence, analogues of
\autoref{th.infty} and \autoref{rem.growth}. Details are left to the reader.


\subsection{Proof of \autoref{prop.52}}\label{proof.52}
Consider the lattice $\bZ_{52}$ spanned (modulo kernel) by the lines
constituting the configuration $Z_{52}$; the polarization
$h\in\bZ_{52}$, $h^2=4$, is recovered as the sum of any four lines constituting a
``plane'' (see~\cite{DIS} for details). It is shown in~\cite{DIS} that
$\rank\bZ_{52}=19$ and $\bZ_{52}$ admits a unique, up to isomorphism
preserving~$h$, embedding into $L=H_2(X)$; the orthogonal complement
$\bZ_{52}^\perp$ is the rank~$3$ lattice~$W$ with the Gram matrix
\[*
\bmatrix0&2&0\\2&0&0\\0&0&24\endbmatrix.
\]
One can easily obtain the following bounds for a nonzero vector
$a\in\bZ_{52}\dual$ such that $2a\in\bZ_{52}$:
if $a\cdot h=0$ or~$2$, then $a^2\le-2$, and
if $a\cdot h=1$, then $a^2\le-1$.

Consider the vector $v_n:=[-n,1,0]\in\bZ_{52}^\perp$ of square $-4n$ and let
$T_n\cong[4n,0,24]$ be its orthogonal complement. For the corresponding
$K3$-surface $X_n:=X(T_n)$, the lattice $\NS(X)$ is the index~$2$ extension of
$\bZ_{52}\oplus\Z v_n$ by a vector $a\oplus\frac12v_n$, where
$a\in\bZ_{52}\dual$, $2a\in\bZ_{52}$, and $a^2=n\bmod2\Z$.
In view of the bounds above, this lattice
satisfies the hypotheses of
\autoref{th.K3}\iref{i.model.birational} and~\iref{i.model.nonsingular} and,
if $n>1$, the corresponding smooth spatial model
contains exactly $52$ lines, \viz. those contained in $\bZ_{52}$.

For the number of models, we observe that, according to
\autoref{prop.extension}, the number of isomorphism classes of extensions
$W\supset T_n\oplus\Z v$, $v^2=-4n$, is at least $2^{\omega(n)-2}$, where
$\omega(n)$ is the number of distinct prime divisors of~$n$.
\qed

\remark
If $n=1$,
two extra lines appear and
the construction used in the proof of \autoref{prop.52} gives us the
quartic $X_{54}$ in~\cite{DIS};
according to~\cite{DIS}, this is the only quartic with more than $52$
lines missing in \autoref{tab.quartics}.
This inclusion $\bZ_{52}\subset\bX_{54}$ of the configurations has not been
observed before.
\endremark

\subsection{Lines defined over $\Q$}\label{s.Q}
A $K3$-surface~$X$ is said to have \emph{Picard rank~$20$ over $\Q$} if $X$
is singular, defined over~$\Q$, and $\NS(X)$ is generated (over~$\Q$) by
divisors defined over~$\Q$. Naturally, spatial models of such surfaces are
the first candidates for quartics containing many lines defined over~$\Q$.

According to M.~Sch\"{u}tt~\cite{Schutt:Q}, a $K3$-surface~$X$ has Picard
rank~$20$ over~$\Q$ if and only if $X=X(T)$ with $T$ primitive and the
discriminant $-\det T$ of class number~$1$. There are $13$ lattices with this
property, with
\[*
\det T\in\bigl\{ 3, 4, 7, 8, 11, 12, 16, 19, 27, 28, 43, 67, 163 \bigr\}.
\]
Furthermore, the known models show that, in each case,
$\NS(X)$ is also generated \emph{over~$\Z$} by
divisors defined over~$\Q$. It follows (M.~Sch\"{u}tt, private communication)
that, in appropriate coordinates in the
projective space, any model of~$X$ is defined over~$\Q$,
and so are all $(-2)$-curves
in~$X$.

\proof[Proof of \autoref{th.Q}]
Comparing the above list of discriminants and \autoref{th.others},
we conclude that
only the last surface, \viz. $X:=X([2,1,82])$, may have smooth spatial
models, and only for these models one may have $\rank\Fano_\Q(X)=20$.
Take for the test lattice
\[*
\bbT:=\bmatrix
      2&0&0&1&0\\
      0&2&0&1&0\\
      0&0&2&1&1\\
      1&1&1&2&0\\
      0&0&1&0&164
\endbmatrix
\]
with $\rt(\bbT)=\bD_4$; it admits isometries only to those Niemeier lattices~$N$
for which $\rt(N)$ has a component of type~$\bD_*$ or~$\bE_*$
(\cf. \autoref{tab.roots}).

Given a root lattice~$R$, define the \emph{minimal dense square}
\[*
\mds(R):=\min\bigl\{a^2\bigm|a\in R\dual,\ \rt(a^\perp\cap R)=0\bigr\}.
\]
Using the description of the irreducible root lattices (see
\autoref{s.root}), one can see that
\[*
\gathered
\mds(\bA_n)=n(n+1)(n+2)/12,\qquad
\mds(\bD_n)=n(n-1)(2n-1)/6,\\
\mds(\bE_6)=78,\qquad
\mds(\bE_7)=399/2,\qquad
\mds(\bE_8)=620.
\endgathered
\]
In the terminology and notation of \autoref{s.many.roots}, if there exists
a dense isometry $\bbT\into N$, we must have $\sum_{k\in I}\mds(R_k)\le163$ for at
least one subset~$I$ obtained from the index set~$\IS$ by removing one
component of type~$\bD_*$ or~$\bE_*$. This observation eliminates all
Niemeier lattices~$N$
with $\rt(N)\ne\bA_7^2\oplus\bD_5^2$,
$\bA_5^4\oplus\bD_4$, or $\bD_4^6$.

For each of the remaining three lattices, there is an essentially unique
isometry $\bD_4\into N$, and its extensions to~$\bbT$ are found by
enumerating the dense vectors in the other components of $\rt(N)$ and
taking into account
the kernel $N/\!\rt(N)$, see, \eg,
\cite[Chaprer 16]{Conway.Sloane}.
Omitting the details and not attempting the complete classification, we
merely state the result: there are over eleven thousands of $\OG(N)$-orbits of
dense isometries, which give rise to $3216$ configurations of lines
distinguishable by simple combinatorial invariants
(rank, pencil structure, and linking structure).

The number $\ls|\Fn_\Q X|$ of lines in the configurations found
(recall that all lines can be assumed defined over~$\Q$)
takes values in the set
$\{26,28,30,31,\ldots,41,42,46\}$.
Ironically, one has $\rank\Fano_\Q(X)=19$ whenever $\ls|\Fn_\Q X|\ge42$.
The complete list of sizes $\ls|\Fn_\Q X|$ of the configurations of rank~$20$
is $\{28,30,31,\ldots,40,41\}$.
\endproof

\remark\label{rem.Q}
The extremal model of $X([2,1,82])$ containing $46$ lines is unique; it
is given by a certain isometry $\bbT\into N(\bA_7^2\oplus\bD_5^2)$.
Using another test lattice, one can obtain the same model from the
isometry $V\into N(\bA_1^{24})$ given by the diagram
\[*
\diagram[Q_{46}]{
\n\ \ \ \ \ \ \ \ \ \ \ \ \ \ \ \n\n\\
\-\-\-\-\-\-\-\-\-\-\-\-\-\-\-\-\\
\ \ \ \ \ \ \ \ \ \ \ \ \=\-\-\-\-\-\-\-\.\.\.\.}
\label{Q46}
\]
(see~\eqref{eq.diagram} for the notation),
extending to a unique $\OG(N)$-orbit.
\endremark

\subsection{Lines in singular quartics}\label{proof.singular}
Let $\Gf\:X:=X(T)\to\Cp3$ be a spatial model as in \autoref{th.K3}, given by
a triple $(S,[\Gg],c)$, so that $S=h^\perp\subset\NS(X)$.
Assume that $\Gf$ is birational, but not smooth, \ie,
$\rt(-S)\ne0$. Then, one of the Weyl chambers $\Delta$ of
the root lattice $\rt(-S)$ is a face of
the K\"{a}hler cone of~$X$. Denoting by $r_1,\ldots,r_n\in S$ the
primitive vectors orthogonal to the facets of~$\Delta$ (these vectors are
roots in $-S$), one can easily see that a $(-2)$-class
$l\in\NS(X)$ such that $l\cdot h=1$ represents an \emph{irreducible}
$(-2)$-curve if and only if $l\cdot r_i\ge0$ for each $i=1,\ldots,n$.
Hence, Statements~\iref{i.model.nonsingular} and~\iref{i.model.lines}
of \autoref{th.K3} take the following form.

\lemma\label{lem.lines}
In the notation introduced above, assuming $\Gf$ birational, one has\rom:
\roster[3]
\item
the classes of the exceptional
divisors contracted by~$\Gf$ are $r_1,\ldots,r_n$\rom;
\item
the straight lines contained in
$\Gf(X)$ are in a one-to-one correspondence with the vectors $a\in S\dual$
representing the class~$\Gg$ and such that
$a^2=-\frac94$ and $a\cdot r_i\ge0$ for each
$i=1,\ldots,n$.
\done
\endroster
\endlemma

In particular, the set of singularities of the image $\Gf(X)$ can be
identified with the ``type'' of the root system $\rt(-S)$: there is one
simple  (\ie, $\bA$--$\bD$--$\bE$ type) singular point for each irreducible
component of $\rt(-S)$ of the same name. Furthermore, since any two Weyl
chambers are in the same orbit of the Weyl group
(and the latter extends to any overlattice containing~$S$), for the purpose of
\emph{counting} lines one can take
for~$\Delta$ any Weyl chamber;
in other words, $r_1,\ldots,r_n$ is any ``standard'' basis
for the root system $\rt(-S)$.

\proof[Proof of \autoref{th.singular}]
Consider the primitive
embedding $V\into N:=N(\bA_1^{24})$ described by the diagram
\[
\diagram[P_{52}]{
\-\-\-\-\-\-\-\-\ \ \n\n\ \ \ \ \ \ \ \ \ \ \d\d\\
\-\-\-\-\ \ \ \ \-\-\-\-\\
\-\-\ \ \=\=\ \ \-\-\ \ \-\-\-\-\-\-\-\-\-\-
}
\]
(see \eqref{eq.diagram} for the notation;
in addition, we use dots $\cdot$ to indicate the components of
$\rt(N)$ that are not in the support of the image,
resulting in extra roots in the orthogonal complement)
and let $S:=-V^\perp$.
A straightforward computation using
\autoref{th.K3} and
\autoref{lem.lines} shows that this lattice~$S$ gives rise to a unique
spatial model $\Gf:X([4,0,12])\to\Cp3$, this model is birational, and its
image has two type~$\bA_1$ singularities (simple nodes) and contains $52$
lines.
\endproof

Since very little is known about lines in singular quartics, we conclude this
section with a short list of quartics with relatively many lines and
singularities. This list is a result of a (non-exhaustive) search of
isometries of test lattices of small discriminant to $N(\bA_1^{24})$,
$N(\bA_2^{12})$, or $N(\bA_3^8)$. Listed below are the combinatorial type
(see~\eqref{eq.diagram}, \eqref{eq.diagram.A2}, \eqref{eq.diagram.A3} for the
notation),
the transcendental lattice~$T$, and the set of singularities of the
corresponding quartic, represented as a root lattice;
the subscript in the notation is the number of
lines.
\begin{alignat}3
&\diagram[P_{48}]{
\-\-\-\-\-\-\-\-\ \ \ \ \ \ \ \ \n\n\ \ \d\d\d\d\\
\ \ \ \ \ \ \ \ \-\-\-\-\-\-\-\-\\
\ \ \ \ \ \ \ \ \-\-\-\-\ \ \ \ \-\-\-\-
}\punct,&\quad&[4,0,8]\punct,&\quad&\bA_1^4,\\\noalign{\allowbreak}
&\diagram[P_{50}]{
\n\ \ \ \n\n\ \ \ \ \ \ \ \ \ \ \ \ \ \ \ \ \ \d\\
\-\-\-\-\-\-\-\-\\
\-\-\-\-\ \ \ \ \-\-\-\-\-\-\-\-\-\-\-\-\.\.\.
}\punct,&&[4,1,16]\punct,&&\bA_1,\\\noalign{\allowbreak}
&\cusptrue\diagram[P_{44}]{
\l\l\ \ \ \ \ \ \ \ \c\d\\
\>\>\>\>\>\>\\
\<\<\<\<\<\<\c\c\c\c
}\punct,&&[8,4,10]\punct,&&\bA_2,\\\noalign{\allowbreak}
&\tacnodetrue\diagram[P_{40}]{
\ \ \ \ \ \C\rr{\ \*\-\*}\d\\
\+\+\+\+\ \ \ \ \\
\-\-\-\-\C
}\punct,&&[6,0,6]\punct,&&\bA_3\oplus\bA_1^2.
\end{alignat}
There are a few other quartics with $48$ lines and one to three simple nodes
and a large number of quartics with fewer lines.

{
\let\.\DOTaccent
\def\cprime{$'$}
\bibliographystyle{amsplain}
\bibliography{degt}
}

\end{document}